\documentclass{article}

\usepackage{amsmath}
\usepackage{amssymb}
\usepackage{authblk}
\usepackage{hyperref}
\usepackage{graphicx}
\usepackage{framed}
\usepackage{setspace}
\usepackage{fullpage}
\usepackage{tikz-cd}

\title{Applied Category Theory in the Wolfram Language using \textsc{Categorica} I: Diagrams, Functors and Fibrations}
\author[1]{Jonathan Gorard}
\affil[1]{Princeton University\protect\\
Princeton, NJ, United States\footnote{\href{mailto:gorard@princeton.edu}{gorard@princeton.edu}}\footnote{Current affiliation. Research and development work was largely performed at, and supported by, the Wolfram Institute.}}

\begin{document}

\maketitle

\begin{abstract}
This article serves as a preliminary introduction to the design of a new, open-source applied and computational category theory framework, named \textsc{Categorica}, built on top of the Wolfram Language. \textsc{Categorica} allows one to configure and manipulate abstract quivers, categories, groupoids, diagrams, functors and natural transformations, and to perform a vast array of automated abstract algebraic computations using (arbitrary combinations of) the above structures; to manipulate and abstractly reason about arbitrary universal properties, including products, coproducts, pullbacks, pushouts, limits and colimits; and to manipulate, visualize and compute with strict (symmetric) monoidal categories, including full support for automated string diagram rewriting and diagrammatic theorem-proving. In so doing, \textsc{Categorica} combines the capabilities of an abstract computer algebra framework (thus allowing one to compute directly with epimorphisms, monomorphisms, retractions, sections, spans, cospans, fibrations, etc.) with those of a powerful automated theorem-proving system (thus allowing one to convert universal properties and other abstract constructions into (higher-order) equational logic statements that can be reasoned about and proved using standard automated theorem-proving methods, as well as to prove category-theoretic statements directly using purely diagrammatic methods). Internally, \textsc{Categorica} relies upon state-of-the-art graph and hypergraph rewriting algorithms for its automated reasoning capabilities, and is able to convert seamlessly between diagrammatic/combinatorial reasoning on labeled graph representations and symbolic/abstract reasoning on underlying algebraic representations of all category-theoretic concepts. In this first of two articles introducing the design of the framework, we shall focus principally upon its handling of quivers, categories, diagrams, groupoids, functors and natural transformations, including demonstrations of both its algebraic manipulation and theorem-proving capabilities in each case.
\end{abstract}

\section{Introduction}

As a field of pure mathematics, category theory emerged in the mid-20th century with the work of Eilenberg and Mac Lane on group theory and algebraic topology\cite{eilenberg}\cite{eilenberg2}, as well as the work of Serre and Grothendieck in the closely related area of homological algebra\cite{serre}\cite{grothendieck}. When conceived as an alternative foundation for mathematics (as exemplified by work on \textit{elementary topos theory} within mathematical logic\cite{mclarty}\cite{johnstone}), the shift from set theory to category theory necessitates a certain shift in philosophical perspective, wherein one transitions away from identifying mathematical objects based upon their internal structure (such as in Zermelo-Fraenkel set theory, where the axiom of extensionality defines set equality based purely upon the set membership relation\cite{kunen}), and towards identifying mathematical objects based upon their relationships to other objects of the same ``type'' (as exemplified by the Yoneda lemma, arguably the most fundamental result in category theory\cite{riehl}). This more \textit{relational} perspective on the nature of structure turns out to be a useful lens through which to view many fields, including many outside of mathematics altogether. For instance, in the foundations of quantum mechanics, moving from set-theoretic to category-theoretic foundations inevitably changes one's view from quantum systems and their states as being the fundamental objects of study (as in the traditional Dirac-von Neumann axioms of quantum mechanics, represented in terms of operator theory on Hilbert space\cite{dirac}\cite{vonneumann}) to quantum processes and their compositions as being the fundamental objects of study (as in the categorical quantum mechanics formalism of Abramsky and Coecke\cite{abramsky}\cite{abramsky2}, later developed into a fully diagrammatic theory of quantum information by Coecke and Duncan\cite{coecke}\cite{coecke2}). Likewise, in computer science, one is able to move from a set-theoretic view in which program state is treated as fundamental (as in the Turing machine picture of computation, and as exemplified by the imperative programming paradigm) to a category-theoretic view in which functions and their compositions are treated as fundamental (as in the ${\lambda}$-calculus picture of computation, and as exemplified by the functional programming paradigm). Such process-theoretic models, in which processes and their algebra of composition are treated on a fundamentally ``higher'' footing than states and their internal structure, have also proved useful and/or instructive for the foundations of natural language processing\cite{coecke3}, cybernetics\cite{capucci}, machine learning\cite{cruttwell}, complex networks and control theory\cite{baez}\cite{baez2}, computational complexity theory\cite{gorard8} (based on prior work on computational irreducibility\cite{gorard9}), database systems\cite{spivak}, and many other domains of knowledge. This general process of applying foundational methods and concepts from category theory to fields outside of pure mathematics has come to be known as \textit{applied category theory}\cite{fong}.

With the advent and maturation of so many fruitful domains of applied category theory, many software tools have now been developed to facilitate working with category-theoretic data structures in an automated or semi-automated fashion, such as \textsc{\href{https://algebraicjulia.github.io/Catlab.jl/dev/}{Catlab.jl}}\cite{halter}, built on top of the Julia programming language, which uses the formalism of generalized algebraic theories and dependent types to facilitate automated algebraic manipulation of certain fundamental data structures such as operads and symmetric monoidal categories. There are also many examples of diagrammatic proof assistants, such as the highly general \textsc{Quantomatic}\cite{kissinger} and \textsc{Cartographer}\cite{sobocinski} frameworks, or the more specialized \textsc{PyZX} framework\cite{kissinger2}, which facilitate the automated rewriting of string diagrams in symmetric monoidal category theory (with, in the latter case, a particular emphasis on diagrammatic quantum information theory). There are even projects, such as the ANR \textsc{\href{https://coreact.wiki}{CoREACT}} project, which aim to formalize certain key aspects of applied category theory, such as compositional/categorical rewriting theory, within existing proof assistant frameworks such as \textsc{Coq}. The aim of the present article is to introduce another such software framework, known as \textsc{Categorica}, which is written principally in the Wolfram Language and is designed to be fully integrated into the Mathematica software system, and which seeks to combine many of the computational algebraic capabilities of existing frameworks such as \textsc{\href{https://algebraicjulia.github.io/Catlab.jl/dev/}{Catlab}} with many of the diagrammatic theorem-proving capabilities of frameworks such as \textsc{Quantomatic}. At its core, the \textsc{Categorica} framework consists of a unified collection of advanced symbolic and diagrammatic algorithms for efficiently converting between purely algebraic representations of categories, diagrams, functors and other key category-theoretic constructions (by employing an entirely presentation-theoretic view of categories, with arrows within an underlying quiver acting as generators, composition acting as the fundamental binary operation, and algebraic equivalences between objects and morphisms acting as relations) and purely graph-theoretic/combinatorial representations of the very same structures (by employing a description of quivers, categories and functors in terms of graphs, hypergraphs and combinatorial rewriting systems). In this manner, \textsc{Categorica} is able to combine the features of a diagrammatic theorem-prover, a higher-order (equational) logic theorem prover, a (hyper)graph rewriting framework, and an abstract computer algebra system, all in an entirely seamless fashion, automatically adopting the most appropriate algorithmic approach and most efficient concrete data structures for any particular problem, and then converting the result back into the desired abstract representation at the end.

The first in this series of two articles introducing the design of \textsc{Categorica} will focus primarily on the core algebraic structures of the framework (with a particular emphasis on \textsc{Categorica}'s representation and handling of quivers, categories, functors and natural transformations), some simple diagrammatic theorem-proving capabilities (such as \textsc{Categorica}'s ability to derive and prove necessary and sufficient algebraic conditions for diagrams to commute, for categories to be groupoids, for functors to be faithful, etc.), as well as some more advanced abstract algebraic capabilities (such as the computation of retractions and sections, epimorphisms and monomorphisms, initial and terminal objects, constant and coconstant morphisms, injectivity/surjectivity of functors on both objects and morphisms, and even some simple cases of Grothendieck fibrations). The forthcoming second article in the series will focus upon the extension of these capabilities to the handling of universal properties (especially products and coproducts, pullbacks and pushouts, and limits and colimits), including illustrations both of how to prove theorems using them, and of how to perform algebraic computations involving them, as well as the applications of these methods to the handling of (strict, symmetric) monoidal categories, including capabilities for string diagram representation, manipulation and rewriting. \textsc{Categorica}'s core algorithmic framework is based on the hypergraph rewriting/Wolfram model formalism outlined within \cite{gorard}\cite{gorard2}\cite{gorard3}, which may be augmented with basic techniques from the theory of graph grammars, and algebraic/compositional graph rewriting theory\cite{ehrig}\cite{ehrig2} (and, in particular, the theory of \textit{double-pushout} rewriting\cite{habel}) in order to develop a highly efficient (hyper)graphical/diagrammatic theorem-proving system\cite{gorard4}\cite{gorard5}\cite{gorard6}. We begin in Section \ref{sec:Section1} by introducing the fundamental data structure underlying all abstract category-theoretic objects in \textsc{Categorica}, namely the abstract \textit{quiver} (a directed multigraph whose vertices are \textit{objects} and whose edges are \textit{arrows}), and we show how an \texttt{AbstractQuiver} object in \textsc{Categorica} may be used in order to ``freely generate'' a corresponding \texttt{AbstractCategory} object by allowing for certain arrows to be composed (at which point they become \textit{morphisms}), as well as introducing identity morphisms for each object, in a manner that is consistent with both associativity and identity axioms. We demonstrate how \textsc{Categorica} is able to keep track of all necessary algebraic equivalences between morphisms (as necessitated by the underlying axioms of category theory) automatically, as well as how new algebraic relations between objects and arrows/morphisms may be introduced, in order to construct more general examples \textit{non-free} categories. We also highlight a very simple application of \textsc{Categorica}'s automated theorem-proving capabilities, by proving that certain categorical diagrams \textit{commute} (i.e. that, within certain categories, all directed paths yield the same morphism up to algebraic equivalence), and by automatically computing both necessary and sufficient algebraic conditions to force non-commuting categorical diagrams to commute, using the \texttt{AbstractCategory} function.

In Section \ref{sec:Section2}, we proceed to show some examples of simple algebraic computations that can be performed using the \texttt{AbstractCategory} function, such as the computation of monomorphisms, epimorphisms and bimorphisms (i.e. the category-theoretic generalizations of injective, surjective and bijective functions in mathematical analysis), and the computation of sections, retractions and isomorphisms (i.e. the category-theoretic generalizations of left-invertible, right-invertible and invertible elements in abstract algebra). We introduce a special kind of category known as a \textit{groupoid}, in which all morphisms are isomorphisms, and we show how \textsc{Categorica}'s automated theorem-proving capabilities can again be used in order to prove that certain \texttt{AbstractCategory} objects are or are not groupoids, and indeed to compute necessary and sufficient algebraic conditions (whenever they exist) to force non-groupoidal \texttt{AbstractCategory} objects to become groupoidal. The relationships between monomorphisms vs. epimorphisms and retractions vs. sections are both examples of a far more general category-theoretic notion of \textit{duality}: the \textit{dual} of a particular construction can be obtained by simply reversing the direction of morphisms, and therefore reversing the order in which morphisms are composed. \textsc{Categorica} has in-built functionality (via the \textit{``DualCategory''} property in particular) for making the systematic exploration of category-theoretic dualities completely straightforward, and we illustrate this by investigating some other common examples of dual constructions using \textsc{Categorica}, including initial vs. terminal objects (i.e. the category-theoretic generalizations of bottom and top elements in partially ordered sets), strict initial vs. strict terminal objects (i.e. initial and terminal objects whose incoming/outgoing morphisms are all isomorphisms), and constant vs. coconstant morphisms (i.e. the category-theoretic generalizations of constant functions and zero-maps in mathematical analysis). In Section \ref{sec:Section3}, we move on to explore functors (i.e. homomorphisms/structure-preserving maps between categories) using the \texttt{AbstractFunctor} function, demonstrating \textsc{Categorica}'s capabilities for handling both covariant and contravariant functors (i.e. functors in which morphism directions are preserved and reversed, respectively). We illustrate how \textsc{Categorica} is able to distinguish between different notions of injectivity, surjectivity and bijectivity for \texttt{AbstractFunctor} objects, including injectivity, surjectivity and bijectivity on objects; essential injectivity, essential surjectivity and essential bijectivity (i.e. injectivity, surjectivity and bijectivity on objects, but only up to isomorphism); and faithfulness, fullness and full faithfulness (i.e. injectivity, surjectivity and bijectivity on \textit{morphisms}). We also show \textsc{Categorica} may be used to investigate an important class of functors known as Grothendieck fibrations (i.e. the category-theoretic generalizations of fiber bundles in topology) from a total category to a base category, including special cases such as discrete fibrations (i.e. Grothendieck fibrations in which the fiber categories are all discrete categories containing only objects and identity morphisms).

Finally, in Section \ref{sec:Section4}, we introduce the \texttt{AbstractNaturalTransformation} function, allowing one to represent arbitrary natural transformations between \texttt{AbstractFunctor} objects in \textsc{Categorica}, and showcase its ability to compute (and prove) both necessary and sufficient algebraic conditions to force transformations between functors to be natural. We demonstrate furthermore how the \texttt{AbstractNaturalTransformation} framework may be used to detect natural isomorphisms between both objects and functors in a highly general way. We also include a short note regarding some of the internal algorithmic details of the \textsc{Categorica} system, and, in particular, its use of double-pushout rewriting formalism and compositional graph rewriting theory to reduce abstract algebraic and diagrammatic reasoning problems to concrete (hyper)graph rewriting problems. We conclude in Section \ref{sec:Section5} with a summary of directions for future research and development. The majority of the core \textsc{Categorica} functionality presented within this article, and its forthcoming companion article, has been fully documented and exposed via the \textit{Wolfram Function Repository}, including functions such as \href{https://resources.wolframcloud.com/FunctionRepository/resources/AbstractQuiver/}{\texttt{AbstractQuiver}}, \href{https://resources.wolframcloud.com/FunctionRepository/resources/AbstractCategory/}{\texttt{AbstractCategory}}, \href{https://resources.wolframcloud.com/FunctionRepository/resources/AbstractFunctor/}{\texttt{AbstractFunctor}}, \href{https://resources.wolframcloud.com/FunctionRepository/resources/AbstractProduct/}{\texttt{AbstractProduct}}, \href{https://resources.wolframcloud.com/FunctionRepository/resources/AbstractCoproduct/}{\texttt{AbstractCoproduct}}, \href{https://resources.wolframcloud.com/FunctionRepository/resources/AbstractPullback/}{\texttt{AbstractPullback}}, \href{https://resources.wolframcloud.com/FunctionRepository/resources/AbstractPushout/}{\texttt{AbstractPushout}} and \href{https://resources.wolframcloud.com/FunctionRepository/resources/AbstractStrictMonoidalCategory/}{\texttt{AbstractStrictMonoidalCategory}}. However, there is also a significant amount of functionality that has not yet been documented and/or exposed, but which is nevertheless available through the \textsc{Categorica} \href{http://github.com/JonathanGorard/Categorica}{GitHub Repository} (with the overall framework currently consisting of over 20,000 lines of symbolic, high-level Wolfram Language code). All explicit examples of categories ${\mathcal{C}}$ presented within this article are \textit{small}, in the sense that the object sets ${\mathrm{ob} \left( \mathcal{C} \right)}$ and morphism sets ${\mathrm{hom} \left( \mathcal{C} \right)}$ are both actually sets, rather than proper classes. Indeed, the only categories currently supported explicitly by \textsc{Categorica} are strictly \textit{finite}, which enables us to bypass any considerations of set-theoretic foundations.

\section{Quivers, Categories and Diagrams}
\label{sec:Section1}

\begin{figure}[ht]
\centering
\begin{framed}
\includegraphics[width=0.595\textwidth]{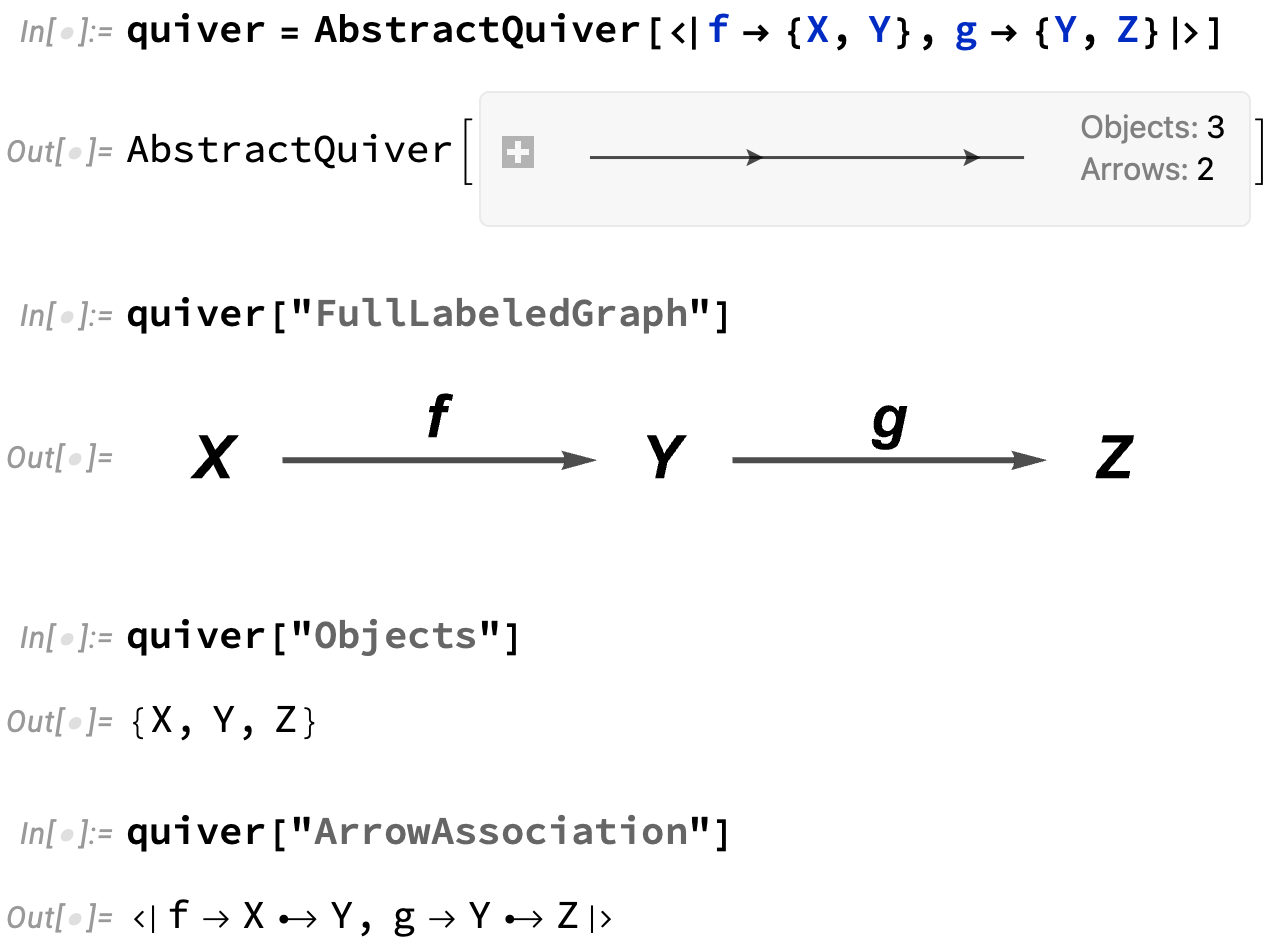}
\vrule
\includegraphics[width=0.395\textwidth]{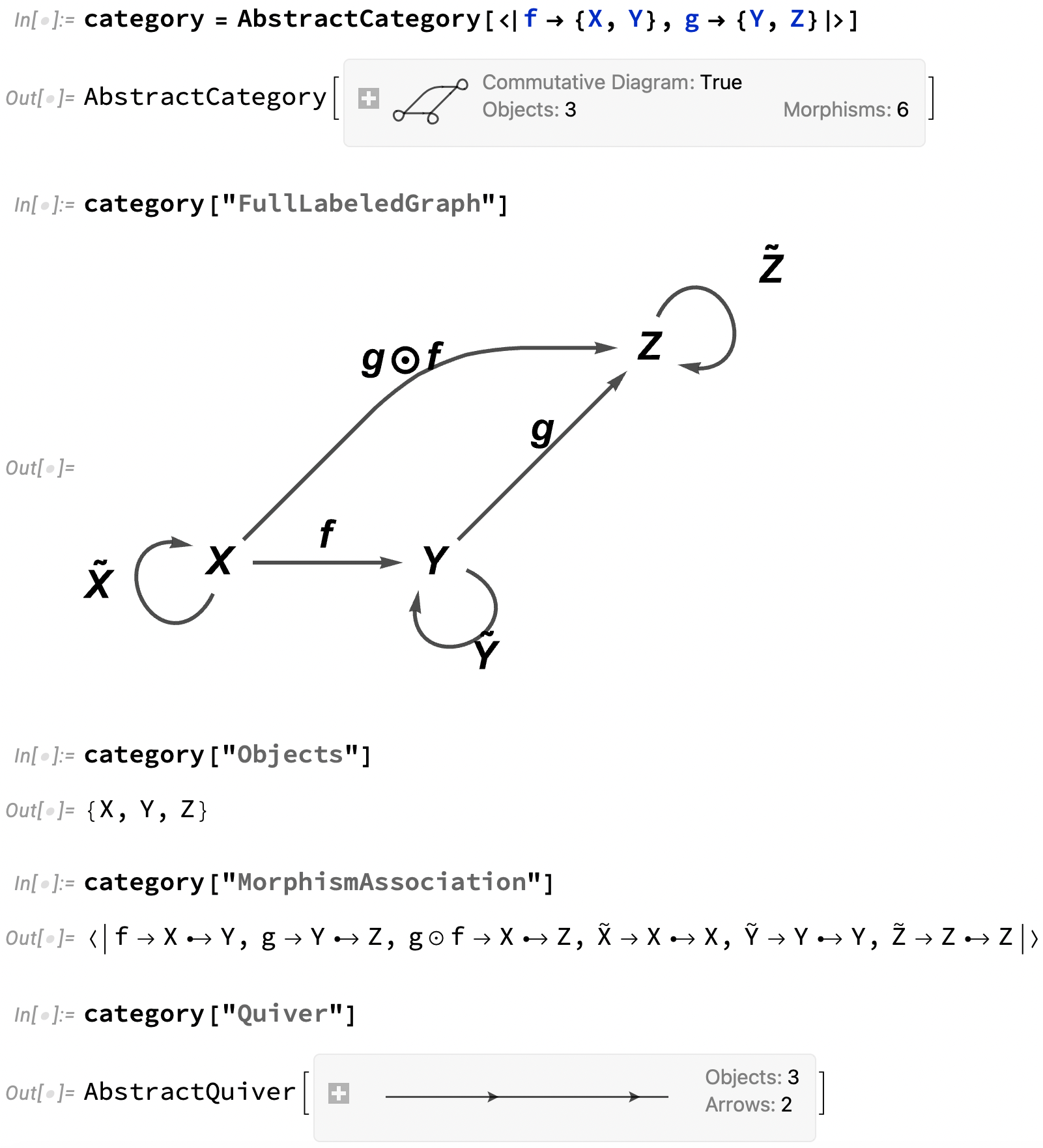}
\end{framed}
\caption{On the left, the \texttt{AbstractQuiver} object for a simple (three-object, two-arrow) quiver. On the right, the \texttt{AbstractCategory} object for the three-object, six-morphism category that is freely generated by this quiver.}
\label{fig:Figure1}
\end{figure}

Every category represented within the \textsc{Categorica} framework has, as its underlying ``skeleton'', a directed multigraph called a \textit{quiver} (i.e. a collection of vertices called \textit{objects}, connected by directed edges called \textit{arrows}, which can have arbitrary multiplicity). One important terminological convention that is enforced consistently throughout the \textsc{Categorica} framework is that the directed edges in a quiver are always known as \textit{arrows}, whereas the directed edges in the corresponding category that is generated by that quiver (potentially freely, potentially with additional algebraic structure) are always known as \textit{morphisms}. By and large, most conventions regarding the notation and nomenclature used within the \textsc{Categorica} framework have been chosen to reflect the conventions set within Mac Lane's classic treatise on the subject\cite{maclane}. If ${\mathcal{Q}}$ is an arbitrary quiver, then we shall denote the set of objects/nodes in ${\mathcal{Q}}$ by ${\mathrm{ob} \left( \mathcal{Q} \right)}$, and the set of directed edges/arrows in ${\mathcal{Q}}$ by ${\mathrm{arr} \left( \mathcal{Q} \right)}$. The rules for how to generate a category ${\mathcal{C}}$ (whose set of objects is denoted ${\mathrm{ob} \left( \mathcal{C} \right)}$ and whose set of morphisms is denoted ${\mathrm{hom} \left( \mathcal{C} \right)}$, by analogy to sets of homomorphisms in abstract algebra) from the underlying quiver ${\mathcal{Q}}$ are then very straightforward: firstly, for every object $X$ or arrow ${f : X \to Y}$ in the quiver ${\mathcal{Q}}$, there is a corresponding object $X$ or morphism ${f : X \to Y}$ in the category ${\mathcal{C}}$:

\begin{equation}
\forall X \in \mathrm{ob} \left( \mathcal{Q} \right), \qquad \exists X \in \mathrm{ob} \left( \mathcal{C} \right), \qquad \text{ and } \qquad \forall \left( f : X \to Y \right) \in \mathrm{arr} \left( \mathcal{Q} \right), \qquad \exists \left( f : X \to Y \right) \in \mathrm{hom} \left( \mathcal{C} \right);
\end{equation}
secondly, for every object $X$ in the quiver ${\mathcal{Q}}$, there exists a corresponding morphism ${id_X : X \to X}$ in the category ${\mathcal{C}}$ mapping that object to itself:

\begin{equation}
\forall X \in \mathrm{ob} \left( \mathcal{Q} \right), \qquad \exists \left( id_X : X \to X \right) \in \mathrm{hom} \left( \mathcal{C} \right);
\end{equation}
and, thirdly, for every pair of morphisms ${f : X \to Y}$ and ${g : Y \to Z}$ in the category ${\mathcal{C}}$ where the \textit{codomain} of the first morphism matches the \textit{domain} of the second, i.e. ${\mathrm{codom} \left( f : X \to Y \right) = \mathrm{dom} \left( g: Y \to Z \right)}$ (in this case, for instance, the domain of the arrow/morphism ${f : X \to Y}$ is $X$, i.e. ${\mathrm{dom} \left( f : X \to Y \right) = X}$, and its codomain is $Y$, i.e. ${\mathrm{codom} \left( f : X \to Y \right) = Y}$), there exists a third morphism ${g \circ f : X \to Z}$ mapping from the domain of the first morphism to the codomain of the second:

\begin{equation}
\forall \left( f : X \to Y \right), \left( g : Y \to Z \right) \in \mathrm{hom} \left( \mathcal{C} \right), \qquad \exists \left( g \circ f : X \to Z \right) \in \mathrm{hom} \left( \mathcal{C} \right),
\end{equation}
or, slightly more abstractly (omitting the arbitrary object labels $X$, $Y$ and $Z$):

\begin{multline}
\forall f, g \in \mathrm{hom} \left( \mathcal{C} \right), \qquad \text{ such that } \qquad \mathrm{codom} \left( f \right) = \mathrm{dom} \left( g \right),\\
\exists \left( g \circ f : \mathrm{dom} \left( f \right) \to \mathrm{codom} \left( g \right) \right) \in \mathrm{hom} \left( \mathcal{C} \right).
\end{multline}
The resulting category ${\mathcal{C}}$ is known as the \textit{free category} generated by the quiver ${\mathcal{Q}}$ (with the arrows of ${\mathcal{Q}}$ acting as the generators, in the algebraic sense, of the morphisms of ${\mathcal{C}}$). This construction can be illustrated diagrammatically by means of the following minimal example:

\begin{equation}
\begin{tikzcd}
& Y \arrow[dr, "g"] &\\
X \arrow[ur, "f"] & & Z
\end{tikzcd} \qquad \mapsto \qquad
\begin{tikzcd}
& Y \arrow[dr, "g"] \arrow[loop above, "id_Y"] &\\
X \arrow[ur, "f"] \arrow[rr, swap, "g \circ f"] \arrow[loop left, "id_X"] & & Z \arrow[loop right, "id_Z"]
\end{tikzcd}.
\end{equation}
As shown in Figure \ref{fig:Figure1}, the \texttt{AbstractQuiver} and \texttt{AbstractCategory} functions in \textsc{Categorica} allow one to represent arbitrary (finite) quivers and categories within the language, with both quivers and free categories being represented internally as lists of objects and associations of arrows/morphisms between those objects (which are then trivially interconvertible with the corresponding combinatorial representations of these structures as labeled directed graphs); every \texttt{AbstractCategory} object comes equipped with an associated \texttt{AbstractQuiver} object representing its underlying ``skeletal'' structure.

\begin{figure}[ht]
\centering
\begin{framed}
\includegraphics[width=0.495\textwidth]{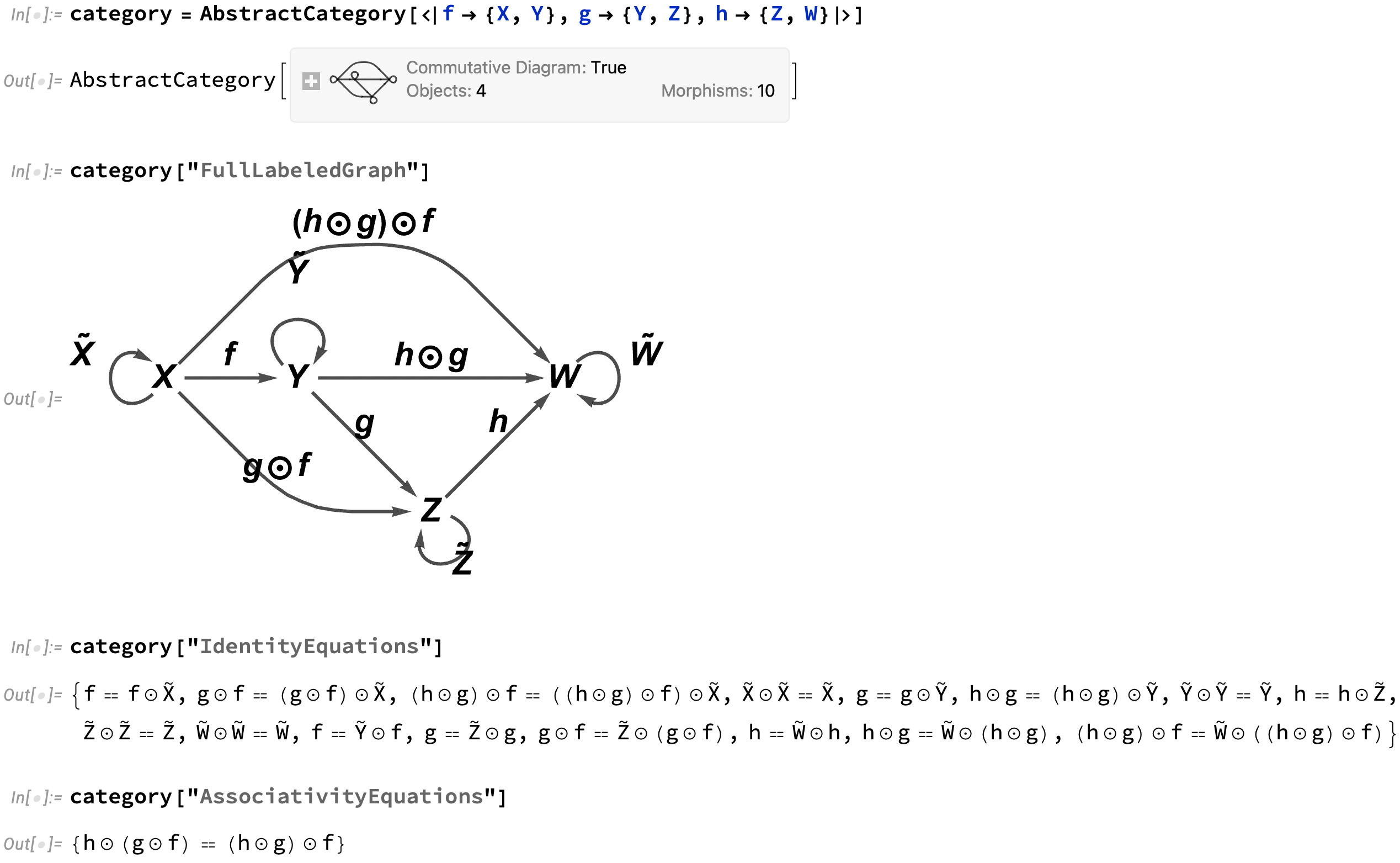}
\vrule
\includegraphics[width=0.495\textwidth]{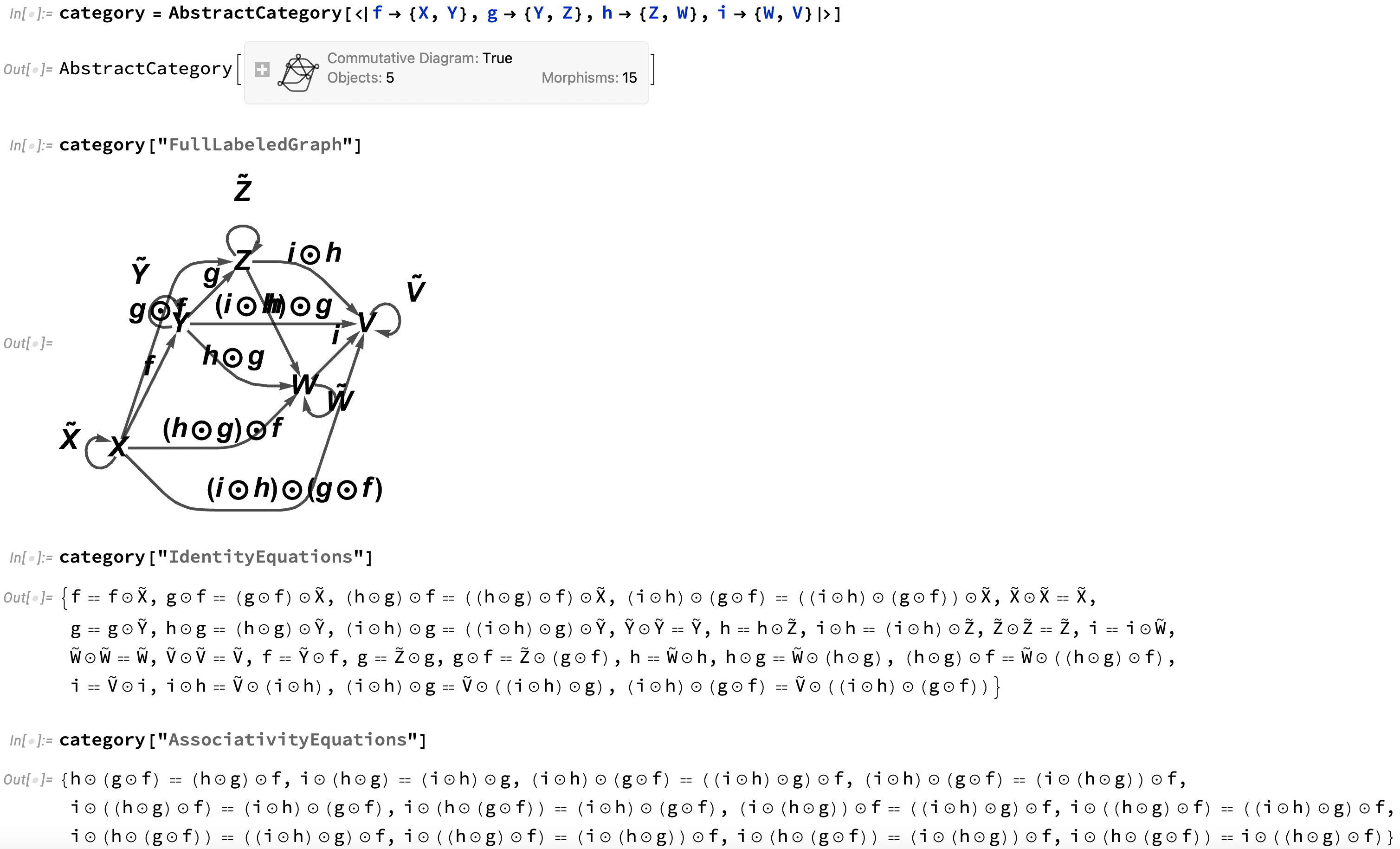}
\end{framed}
\caption{On the left, the \texttt{AbstractCategory} object for a simple (four-object, ten-morphism) category, showing all algebraic equivalences between morphisms that must hold by virtue of the identity and associativity axioms of category theory. On the right, the \texttt{AbstractCategory} object for a slightly larger (five-object, fifteen-morphism) category, again showing all algebraic equivalences between morphisms that must hold by virtue of the identity and associativity axioms of category theory.}
\label{fig:Figure2}
\end{figure}

Within the construction of the free category ${\mathcal{C}}$ described above, the morphisms ${id_X : X \to X}$ that were introduced for each object $X$ are known as \textit{identity} morphisms on $X$, and the abstract binary operation:

\begin{equation}
\circ : \mathrm{hom} \left( \mathcal{C} \right) \times \mathrm{hom} \left( \mathcal{C} \right) \to \mathrm{hom} \left( \mathcal{C} \right),
\end{equation}
that appears in the definition of the morphism ${g \circ f : X \to Z}$ is known as the \textit{composition} operation (such that the resulting morphism ${g \circ f : X \to Z}$ may be referred to as the \textit{composition} of morphisms $f : X \to Y$ and $g : Y \to Z$). The identity morphisms ${id_X : X \to X}$ must satisfy the requisite axioms to act as both left and right identities under the operation of morphism composition, such that if ${f : X \to Y}$ is a morphism in the category ${\mathcal{C}}$, then ${f \circ id_X : X \to Y}$ should be the same as ${f : X \to Y}$ (i.e. ${id_X}$ acts as a right identity on $f$):

\begin{equation}
\forall \left( f: X \to Y \right) \in \mathrm{hom} \left( \mathcal{C} \right), \qquad \left( f \circ id_X : X \to Y \right) = \left( f : X \to Y \right),
\end{equation}
or, illustrated diagrammatically:

\begin{equation}
\begin{tikzcd}
X \arrow[rr, bend left, "f"] \arrow[rr, bend right, swap, "f \circ id_X"] \arrow[loop left, "id_X"] & & Y
\end{tikzcd} \qquad \mapsto \qquad
\begin{tikzcd}
X \arrow[rr, "f = f \circ id_X"] \arrow[loop left, "id_X"] & & Y
\end{tikzcd},
\end{equation}
and, similarly, ${id_Y \circ f : X \to Y}$ should be the same as ${f : X \to Y}$ (i.e. ${id_Y}$ acts as a left identity on $f$):

\begin{equation}
\forall \left( f: X \to Y \right) \in \mathrm{hom} \left( \mathcal{C} \right), \qquad \left( id_Y \circ f : X \to Y \right) = \left( f : X \to Y \right),
\end{equation}
or, illustrated diagrammatically:

\begin{equation}
\begin{tikzcd}
X \arrow[rr, bend left, "f"] \arrow[rr, bend right, swap, "id_Y \circ f"] & & Y \arrow[loop right, "id_Y"]
\end{tikzcd} \qquad \mapsto \qquad
\begin{tikzcd}
X \arrow[rr, "f = id_Y \circ f"] & & Y \arrow[loop right, "id_Y"]
\end{tikzcd}.
\end{equation}
Additionally, the operation of morphism composition itself must satisfy the axiom of associativity, such that if ${f : X \to Y}$, ${g : Y \to Z}$ and ${h : Z \to W}$ form a set of three composable morphisms in the category ${\mathcal{C}}$, then the composition ${\left( h \circ g \right) \circ f : X \to W}$ should be the same as the composition ${h \circ \left( g \circ f \right) : X \to W}$:

\begin{multline}
\forall \left( f : X \to Y \right), \left( g : Y \to Z \right), \left( h : Z \to W \right) \in \mathrm{hom} \left( \mathcal{C} \right),\\
\left( \left( h \circ g \right) \circ f : X \to W \right) = \left( h \circ \left( g \circ f \right) : X \to W \right),
\end{multline}
or, illustrated diagrammatically:

\begin{equation}
\begin{tikzcd}
& Y \arrow[r, "g"] & Z \arrow[dr, "h"] &\\
X \arrow[ur, "f"] \arrow[rrr, "\left( h \circ g \right) \circ f"] \arrow[rrr, bend right, swap, "h \circ \left( g \circ f \right)"] & & & W
\end{tikzcd} \qquad \mapsto \qquad
\begin{tikzcd}
& Y \arrow[r, "g"] & Z \arrow[dr, "h"] &\\
X \arrow[ur, "f"] \arrow[rrr, swap, "\left( h \circ g \right) \circ f = h \circ \left( g \circ f \right)"] & & & W
\end{tikzcd}.
\end{equation}
\textsc{Categorica} automatically keeps track of all algebraic equivalences between morphisms that must be imposed in order to maintain consistency with these identity and associativity axioms, as shown in Figure \ref{fig:Figure2} for the case of two relatively simple \texttt{AbstractCategory} objects. Note that the \texttt{CircleDot[ ${\dots}$ ]} and \texttt{OverTilde[ ${\dots}$ ]} representations of the composition operation ${\circ}$ and the identity morphisms ${id_X}$, respectively, are simply the defaults chosen by \textsc{Categorica}, and can be overridden simply by passing additional arguments to \texttt{AbstractCategory}.

\begin{figure}[ht]
\centering
\begin{framed}
\includegraphics[width=0.445\textwidth]{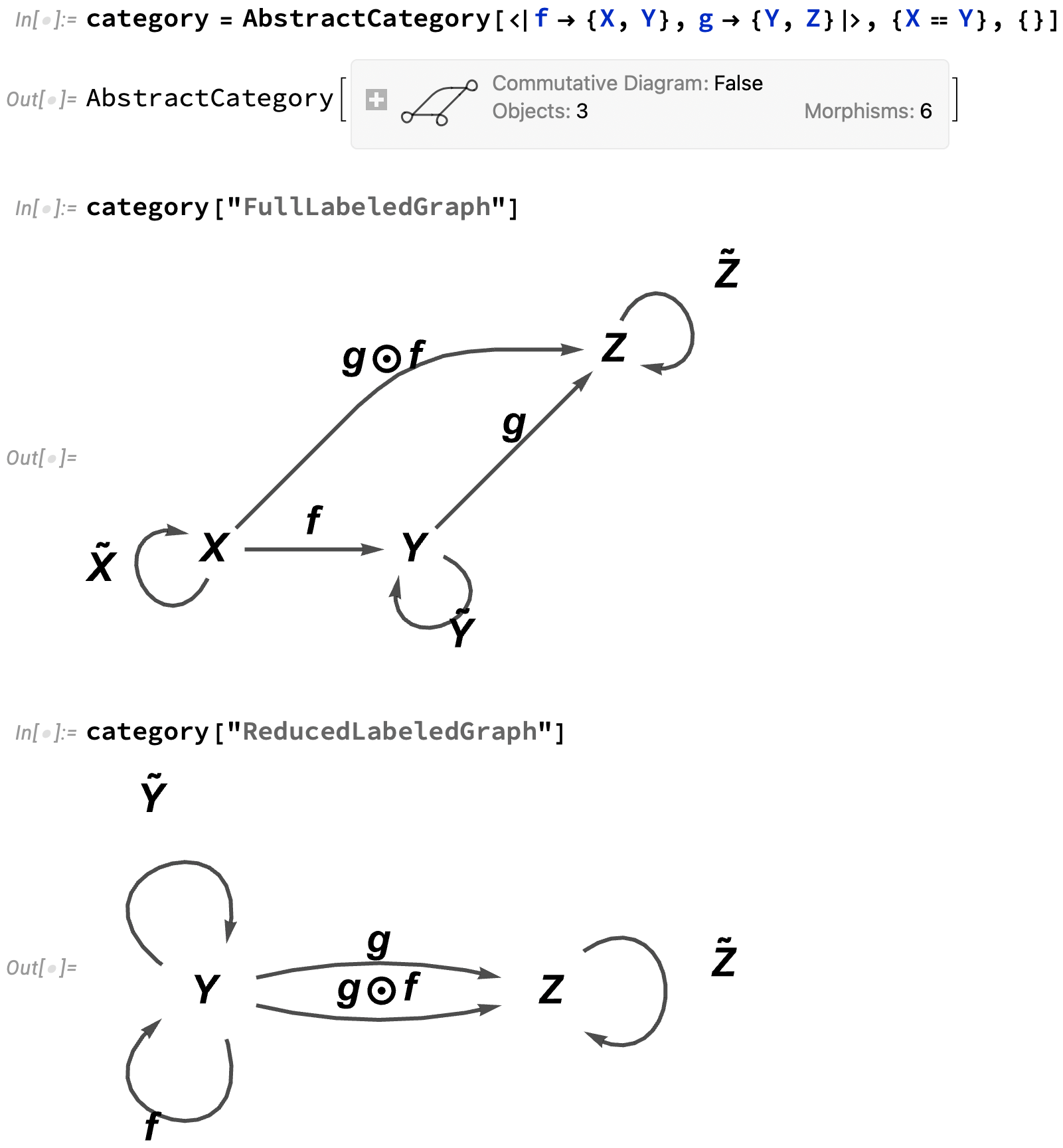}
\vrule
\includegraphics[width=0.545\textwidth]{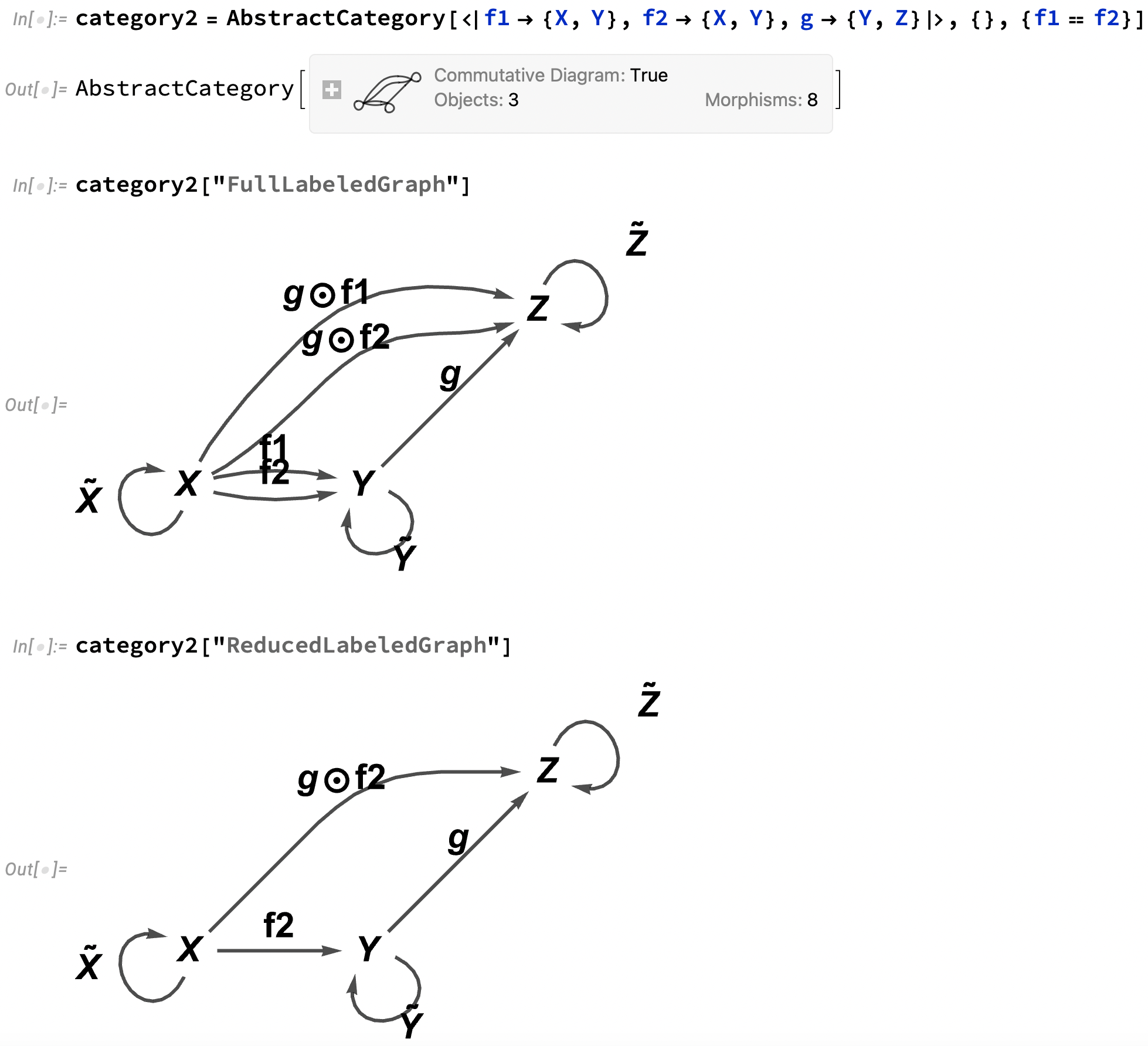}
\end{framed}
\caption{On the left, the \texttt{AbstractCategory} object for a simple category with the object equivalence ${X = Y}$ specified, showing the labeled graph representations of the category both with and without algebraic equivalences imposed. On the right, the \texttt{AbstractCategory} object for a simple category with the morphism equivalence ${\left( f_1 : X \to Y \right) = \left( f_2 : X \to Y \right)}$ specified, showing the labeled graph representations of the category both with and without algebraic equivalences imposed.}
\label{fig:Figure3}
\end{figure}

Although we have considered only \textit{free} categories so far, it is also possible to construct non-free categories which possess additional algebraic structure in the form of additional equivalences between their objects and morphisms (beyond simply the minimal algebraic equivalences necessitated by the axioms of category theory). For instance, we could consider imposing the equivalence ${X = Y}$ between objects $X$ and $Y$ in the following simple category:

\begin{equation}
\begin{tikzcd}
& Y \arrow[dr, "g"] &\\
X \arrow[ur, "f"] \arrow[rr, swap, "g \circ f"] & & Y
\end{tikzcd} \qquad \mapsto \qquad
\begin{tikzcd}
X = Y \arrow[rr, bend left, "g"] \arrow[rr, bend right, swap, "g \circ f"] \arrow[loop left, "f"] & & Z
\end{tikzcd},
\end{equation}
or the equivalence ${\left( f_1 : X \to Y \right) = \left( f_2 : X \to Y \right)}$ between morphisms ${f_1 : X \to Y}$ and ${f_2 : X \to Y}$ in the following, slightly different, simple category:

\begin{equation}
\begin{tikzcd}
& Y \arrow[dr, "g"] &\\
X \arrow[ur, swap, "f_1"] \arrow[ur, bend left, "f_2"] \arrow[rr, swap, "g \circ f_1"] \arrow[rr, swap, bend right, "g \circ f_2"] & & Z
\end{tikzcd} \qquad \mapsto \qquad
\begin{tikzcd}
& Y \arrow[dr, "g"] &\\
X \arrow[ur, "f_1 = f_2"] \arrow[rr, swap, "g \circ f_1 = g \circ f_2"] & & Z
\end{tikzcd}.
\end{equation}
Note that we are only permitted to impose algebraic equivalences between morphisms whose domains and codomains both match, i.e., abstractly, using the functions ${\mathrm{dom} : \mathrm{hom} \left( \mathcal{C} \right) \to \mathrm{ob} \left( \mathcal{C} \right)}$ and ${\mathrm{codom} : \mathrm{hom} \left( \mathcal{C} \right) \to \mathrm{ob} \left( \mathcal{C} \right)}$ introduced above, we have:

\begin{equation}
\forall f, g \in \mathrm{hom} \left( \mathcal{C} \right), \qquad \implies \qquad \mathrm{dom} \left( f \right) = \mathrm{dom} \left( g \right), \qquad \text{ and } \qquad \mathrm{codom} \left( f \right) = \mathrm{codom} \left( g \right).
\end{equation}
Note also that an algebraic equivalence between morphisms, such as ${\left( f_1 : X \to Y \right) = \left( f_2 : X \to Y \right)}$ in the above, will, in general, automatically imply other algebraic equivalences between morphisms, such as:

\begin{equation}
\left( g \circ f_1 : X \to Z \right) = \left( g \circ f_2 : X \to Z \right),
\end{equation}
as shown in the corresponding diagram, although in this particular case the converse does not necessarily hold, as we shall see in more detail later. In Figure \ref{fig:Figure3}, we see the \texttt{AbstractCategory} objects corresponding to the two examples presented above (with object equivalence ${X = Y}$ and morphism equivalence ${\left( f_1 : X \to Y \right) = \left( f_2 : X \to Y \right)}$ specified, respectively), showing the directed graph representations of the corresponding categories both with and without algebraic equivalences imposed. One of the general design principles of \textsc{Categorica} is that the keyword \textit{``Reduced''} within property names such as \textit{``ReducedLabeledGraph''} is used to indicate that all known algebraic equivalences should be applied (as opposed to the keyword \textit{``Full''} within property names such as \textit{``FullLabeledGraph''}, which is used to indicate that no algebraic equivalences should be applied, and that the full, ``maximally-unreduced'' algebraic structure should be presented instead). There is also a related keyword \textit{``Simple''} (as in the property name \textit{``SimpleLabeledGraph''}), which can be used to remove all self-loops and multi-edges from the corresponding directed graph representation of any quiver, category, functor, etc. (for instance in order to yield a more presentable form of a particular diagrammatic representation, of the kind that might be shown within an academic paper). Although this automatic simplification feature is frequently very useful when performing practical diagrammatic manipulations, since the purpose of the present article is to be completely explicit and to illustrate (in as much detail as necessary) the underlying design and algorithmic underpinnings of the \textsc{Categorica} framework, we shall not make use of it here.

\begin{figure}[ht]
\centering
\begin{framed}
\includegraphics[width=0.445\textwidth]{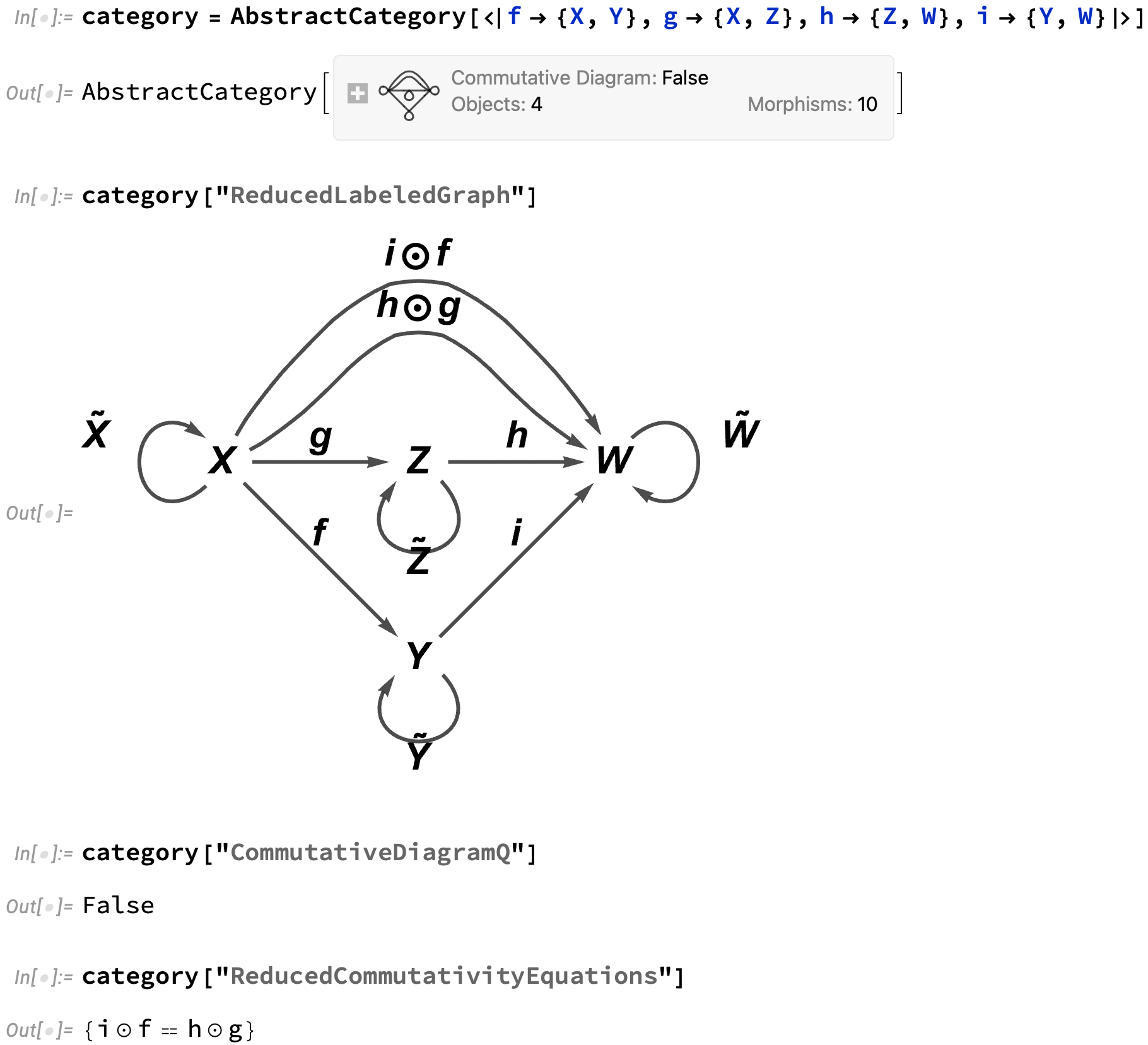}
\vrule
\includegraphics[width=0.545\textwidth]{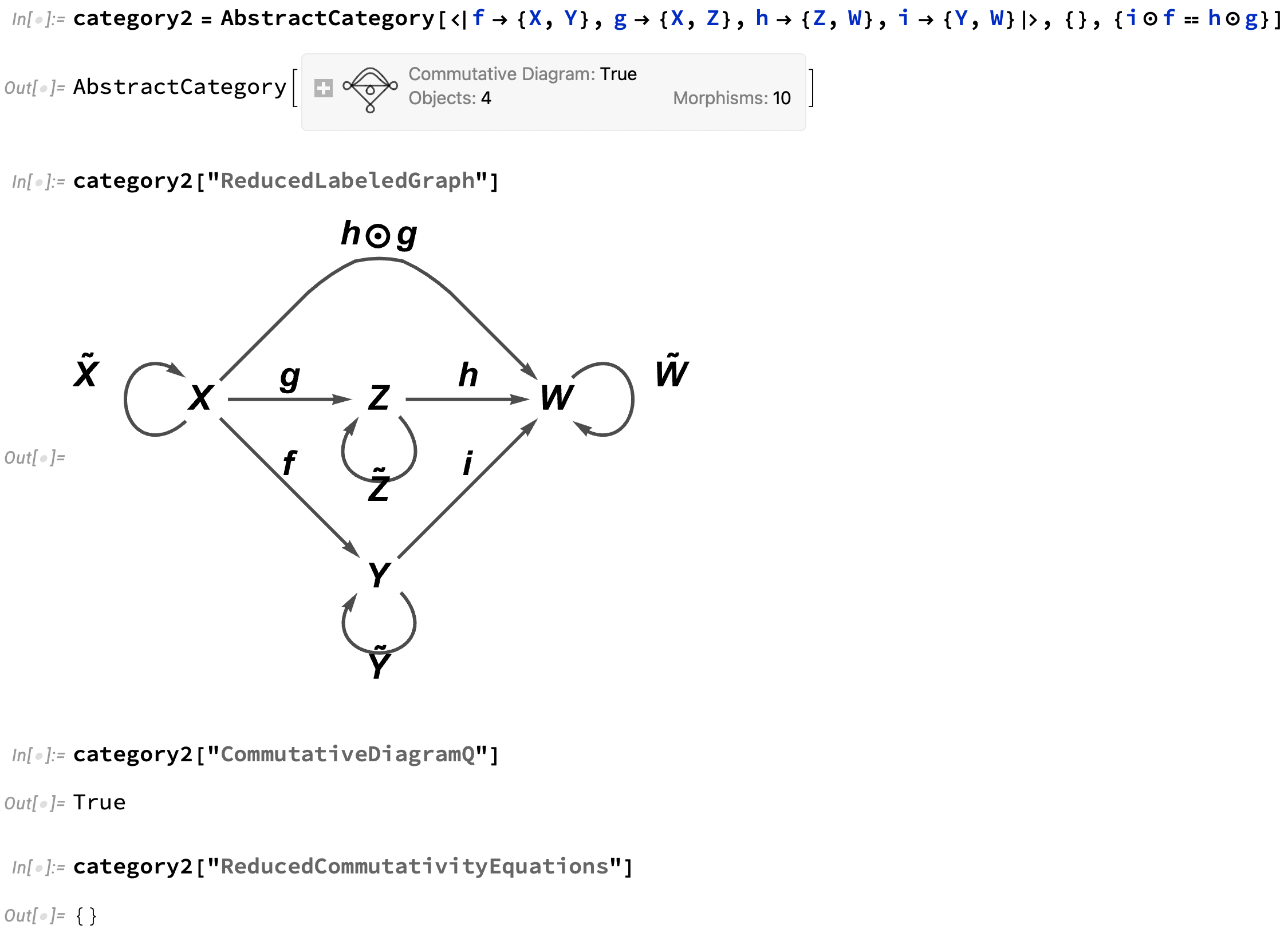}
\end{framed}
\caption{On the left, the \texttt{AbstractCategory} object corresponding to a simple (not yet commutative) square diagram, illustrating that the diagram is not commutative and showing that the morphism equivalence ${i \circ f = h \circ g}$ is the minimal algebraic condition necessary to force the diagram to commute. On the right, the \texttt{AbstractCategory} object for the commutative case of this same square diagram, with the morphism equivalence ${i \circ f = h \circ g}$ imposed, demonstrating that this equivalence is indeed sufficient to force the diagram to commute.}
\label{fig:Figure4}
\end{figure}

\begin{figure}[ht]
\centering
\begin{framed}
\includegraphics[width=0.495\textwidth]{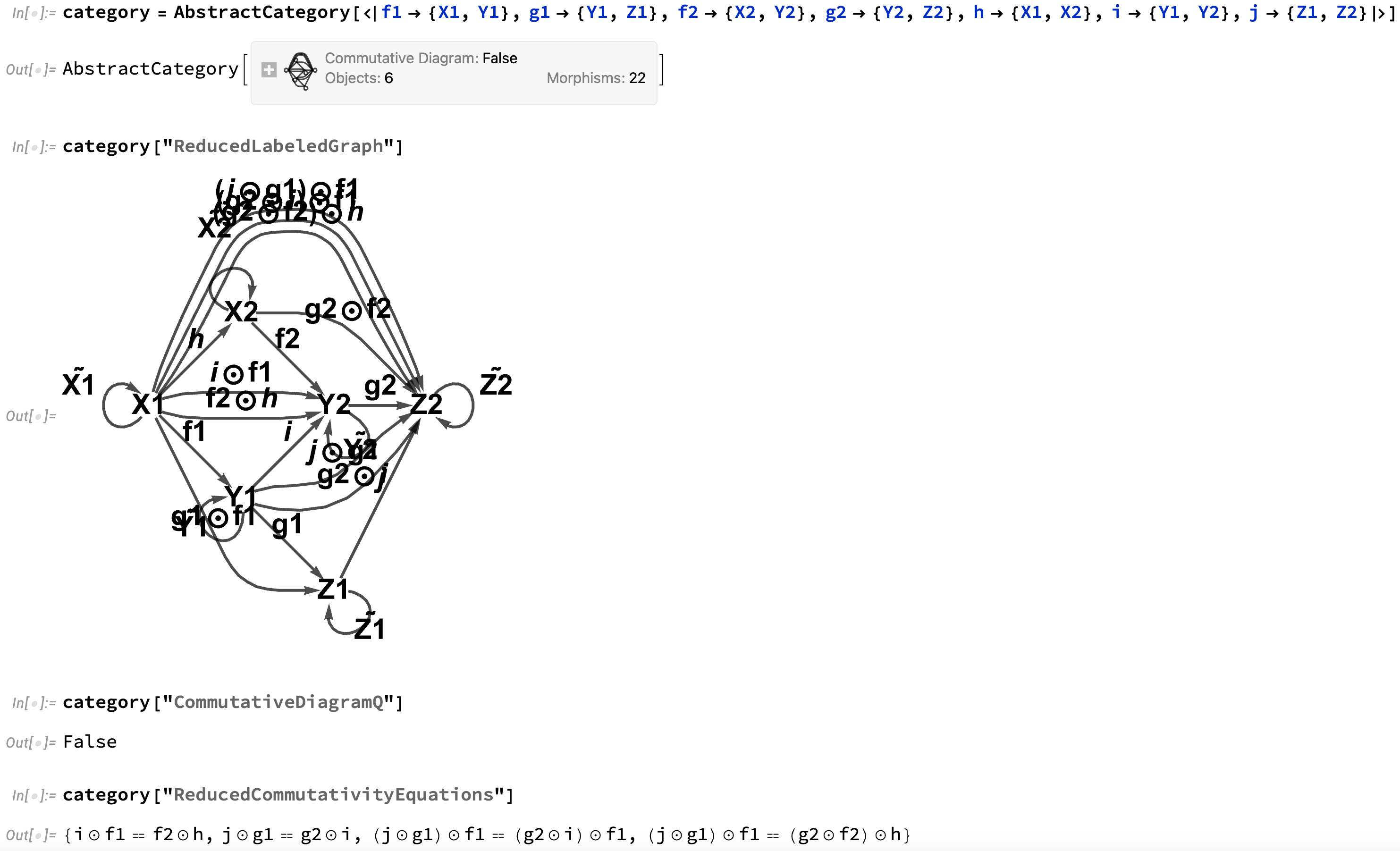}
\vrule
\includegraphics[width=0.495\textwidth]{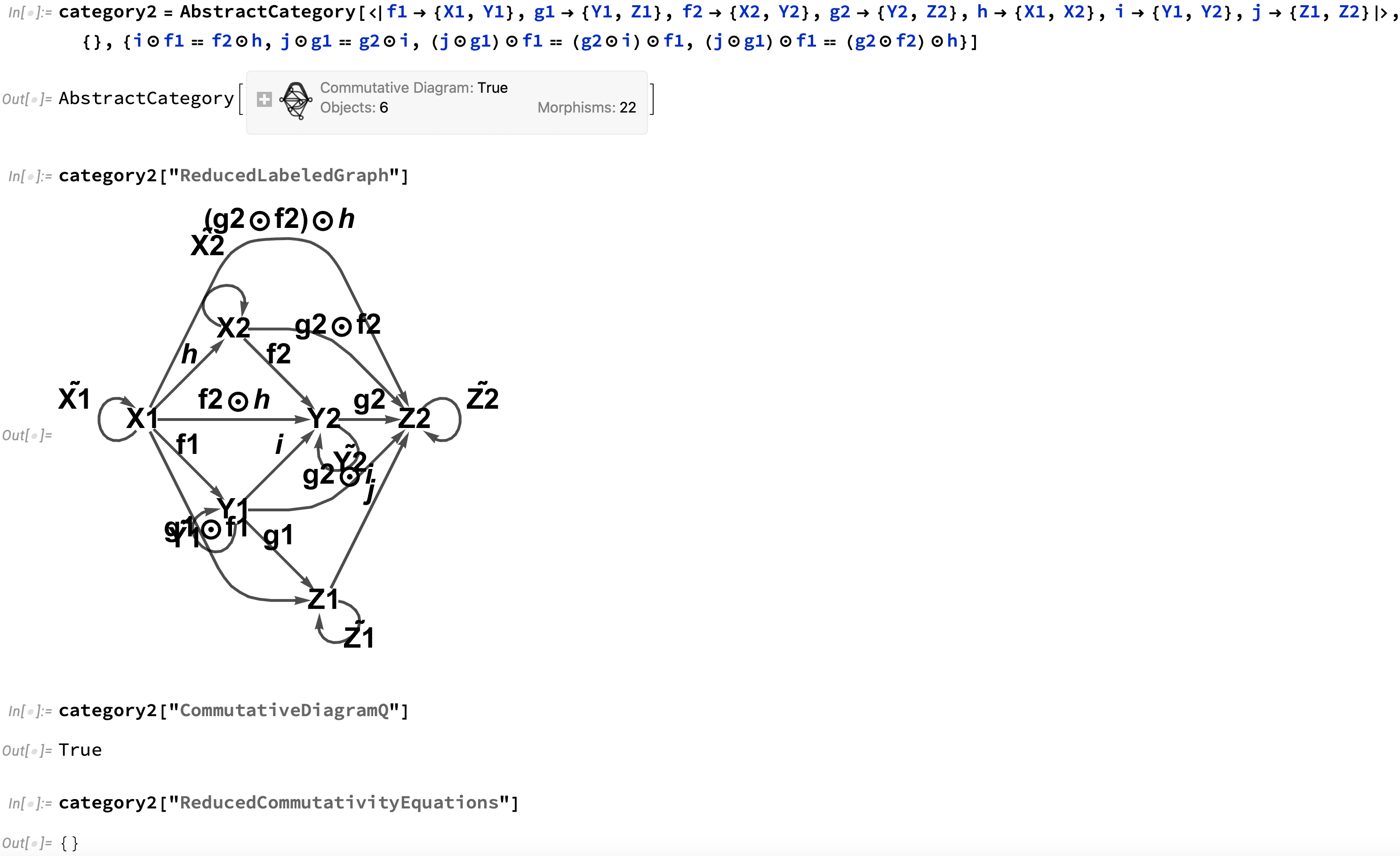}
\end{framed}
\caption{On the left, the \texttt{AbstractCategory} object corresponding to a slightly more complex (not yet commutative) oblong diagram, illustrating that the diagram is not yet commutative and computing the minimum set of morphism equivalences necessary to force the diagram to commute. On the right, the \texttt{AbstractCategory} object for the commutative case of this same oblong diagram, with the aforementioned morphism equivalences imposed, demonstrating that these equivalences are indeed sufficient to force the diagram to commute.}
\label{fig:Figure5}
\end{figure}

One of the most common applications of non-free categories and algebraic equivalences between morphisms, occurring throughout many fields of mathematics (and especially in abstract and homological algebra), is the formal treatment of \textit{commutative diagrams}. Indeed, it is often said that commutative diagrams play the same role in category theory and abstract/homological algebra that equations play in more elementary algebra. For instance, when one says that the following diagram \textit{commutes}:

\begin{equation}
\begin{tikzcd}
X \arrow[rr, "f"] \arrow[dd, swap, "g"] & & Y \arrow[dd, "i"]\\ \\
Z \arrow[rr, swap, "h"] & & W
\end{tikzcd},
\end{equation}
one is effectively imposing the following equivalence between morphisms:

\begin{equation}
\left( h \circ g : X \to W \right) = \left( i \circ f : X \to W \right),
\end{equation}
or, in diagrammatic form, one has:

\begin{equation}
\begin{tikzcd}
X \arrow[rr, "f"] \arrow[dd, swap, "g"] \arrow[ddrr, bend left = 10, "h \circ g"] \arrow[ddrr, bend right = 10, swap, "i \circ f"] & & Y \arrow[dd, "i"] \\ \\
Z \arrow[rr, swap, "h"] & & W
\end{tikzcd} \qquad \mapsto \qquad
\begin{tikzcd}
X \arrow[rr, "f"] \arrow[dd, swap, "g"] \arrow[ddrr, "h \circ g = i \circ f"] & & Y \arrow[dd, "i"]\\ \\
Z \arrow[rr, swap, "h"] & & W
\end{tikzcd}.
\end{equation}
This condition of commutativity of diagrams may thus be characterized purely combinatorially, namely as the condition that all directed paths through the labeled graph representation of the associated category yield the same morphism (up to algebraic equivalence). \textsc{Categorica} is hence able to fall back to using purely graph-theoretic search algorithms in order to compute the minimum set of morphism equivalences necessary to force this diagram to commute, as well as to prove that these equivalences are indeed sufficient, as shown in Figure \ref{fig:Figure4}. However, although the case of a commutative square is very straightforward, it does not take long before this method of \textit{diagram chasing} becomes essentially unmanageable for a human algebraist to enact in full detail; for instance, even upon considering the next obvious case, namely a commutative \textit{oblong} of the form:

\begin{equation}
\begin{tikzcd}
X_1 \arrow[rr, "f_1"] \arrow[dd, swap, "h"] & & Y_1 \arrow[rr, "g_1"] \arrow[dd, "i"] & & Z_1 \arrow[dd, "j"]\\ \\
X_2 \arrow[rr, swap, "f_2"] & & Y_2 \arrow[rr, swap, "g_2"] & & Z_2
\end{tikzcd},
\end{equation}
things have already become quite a bit trickier to analyze than they were for the commutative square, since one not only needs to force the two interior squares to commute:

\begin{equation}
\begin{tikzcd}
X_1 \arrow[rr, "f_1"] \arrow[dd, swap, "h"] \arrow[ddrr, bend left = 10, "i \circ f_1"] \arrow[ddrr, bend right = 10, swap, "f_2 \circ h"] & & Y_1 \arrow[rr, "g_1"] \arrow[dd, "i"] \arrow[ddrr, bend left = 10, "j \circ g_1"] \arrow[ddrr, bend right = 10, swap, "g_2 \circ i"] & & Z_1 \arrow[dd, "j"]\\ \\
X_2 \arrow[rr, swap, "f_2"] & & Y_2 \arrow[rr, swap, "g_2"] & & Z_2
\end{tikzcd} \qquad \mapsto \qquad
\begin{tikzcd}
X_1 \arrow[rr, "f_1"] \arrow[dd, swap, "h"] \arrow[ddrr, "\substack{i \circ f_1\\ = f_2 \circ h}"] & & Y_1 \arrow[rr, "g_1"] \arrow[dd, "i"] \arrow[ddrr, "\substack{j \circ g_1\\ = g_2 \circ i}"] & & Z_1 \arrow[dd, "j"]\\ \\
X_2 \arrow[rr, swap, "f_2"] & & Y_2 \arrow[rr, swap, "g_2"] & & Z_2
\end{tikzcd},
\end{equation}
i.e:

\begin{equation}
\left( i \circ f_1 : X_1 \to Y_2 \right) = \left( f_2 \circ h : X_1 \to Y_2 \right), \qquad \text{ and } \qquad \left( j \circ g_1 : Y_1 \to Z_2 \right) = \left( g_2 \circ i : Y_1 \to Z_2 \right),
\end{equation}
but one must also somehow take care of the commutativity of the outer rectangle:

\begin{equation}
\begin{tikzcd}
X_1 \arrow[rr, "f_1"] \arrow[dd, swap, "h"] \arrow[ddrrrr, bend left = 10, "\left( j \circ g_1 \right) \circ f_1"] \arrow[ddrrrr, bend right = 10, swap, "\left( g_2 \circ f_2 \right) \circ h"] & & Y_1 \arrow[rr, "g_1"] & & Z_1 \arrow[dd, "j"]\\ \\
X_2 \arrow[rr, swap, "f_2"] & & Y_2 \arrow[rr, swap, "g_2"] & & Z_2
\end{tikzcd} \qquad \mapsto \qquad
\begin{tikzcd}
X_1 \arrow[rr, "f_1"] \arrow[dd, swap, "h"] \arrow[ddrrrr, "\substack{\left( j \circ g_1 \right) \circ f_1\\ = \left( g_2 \circ f_2 \right) \circ h}"] & & Y_1 \arrow[rr, "g_1"] & & Z_1 \arrow[dd, "j"]\\ \\
X_2 \arrow[rr, swap, "f_2"] & & Y_2 \arrow[rr, swap, "g_2"] & & Z_2
\end{tikzcd},
\end{equation}
i.e:

\begin{equation}
\left( \left( j \circ g_1 \right) \circ f_1 : X_1 \to Z_2 \right) = \left( \left( g_2 \circ f_2 \right) \circ h : X_1 \to Z_2 \right).
\end{equation}
Nevertheless, despite the additional complexity, \textsc{Categorica} is able to handle this case (and, indeed, the case of much larger and more complicated diagrams) in much the same way, as demonstrated in Figure \ref{fig:Figure5}. A fully rigorous definition and treatment of categorical diagrams necessitates the introduction of a certain functor defined over \textit{index categories} (or \textit{schemes}), which we shall revisit following our introduction to \textsc{Categorica}'s handling of abstract functor objects in Section 
\ref{sec:Section3}. Note that \texttt{AbstractQuiver} objects in \textsc{Categorica} may also carry algebraic equivalence information on both objects and arrows: these equivalences are then translated into corresponding algebraic equivalences on objects and \textit{morphisms} upon promotion of the quiver to a full \texttt{AbstractCategory} object via the pipeline outlined above.

\clearpage

\section{Monos, Epis, Retractions and Sections: The Case of Groupoids}
\label{sec:Section2}

\begin{figure}[ht]
\centering
\begin{framed}
\includegraphics[width=0.495\textwidth]{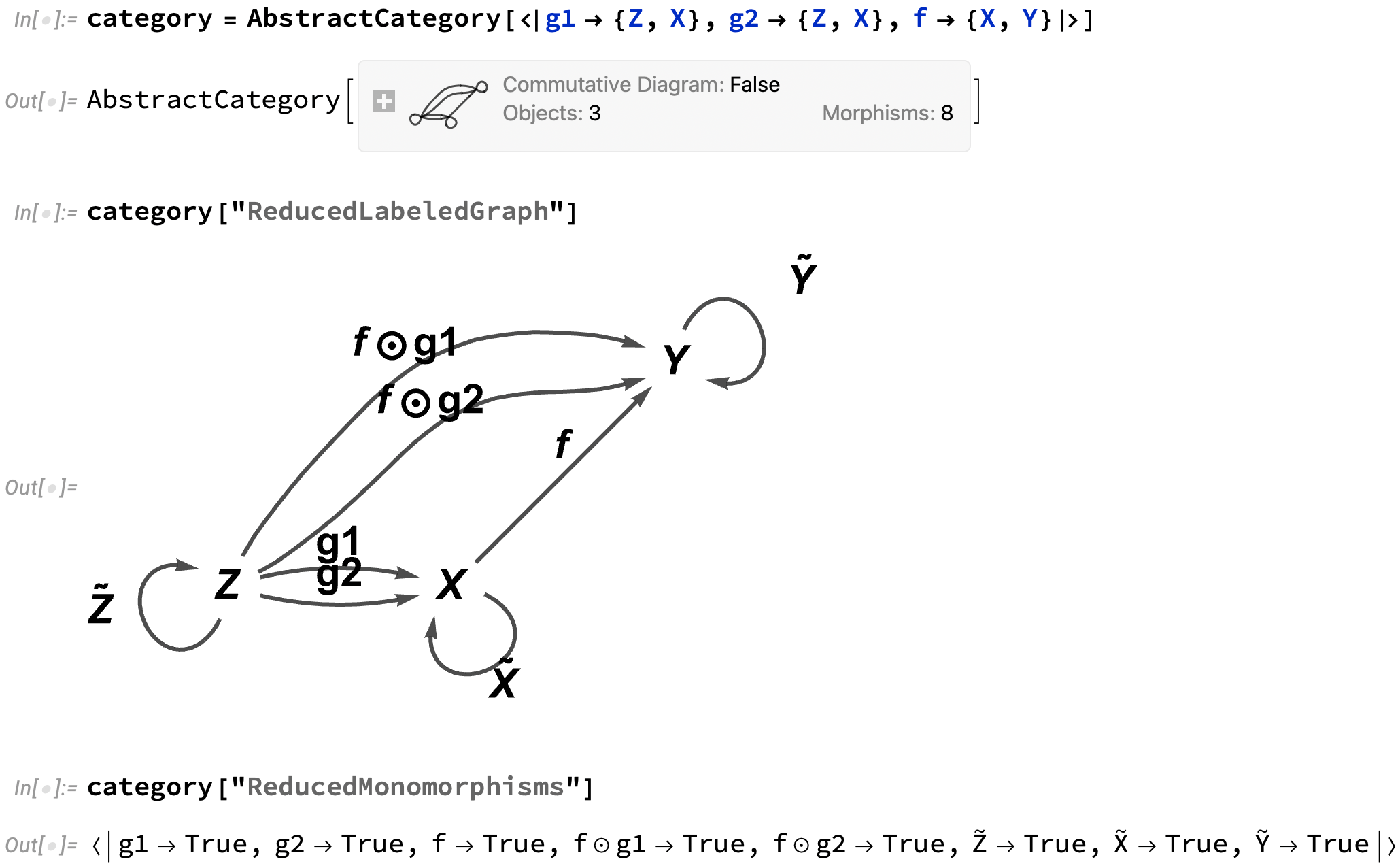}
\vrule
\includegraphics[width=0.495\textwidth]{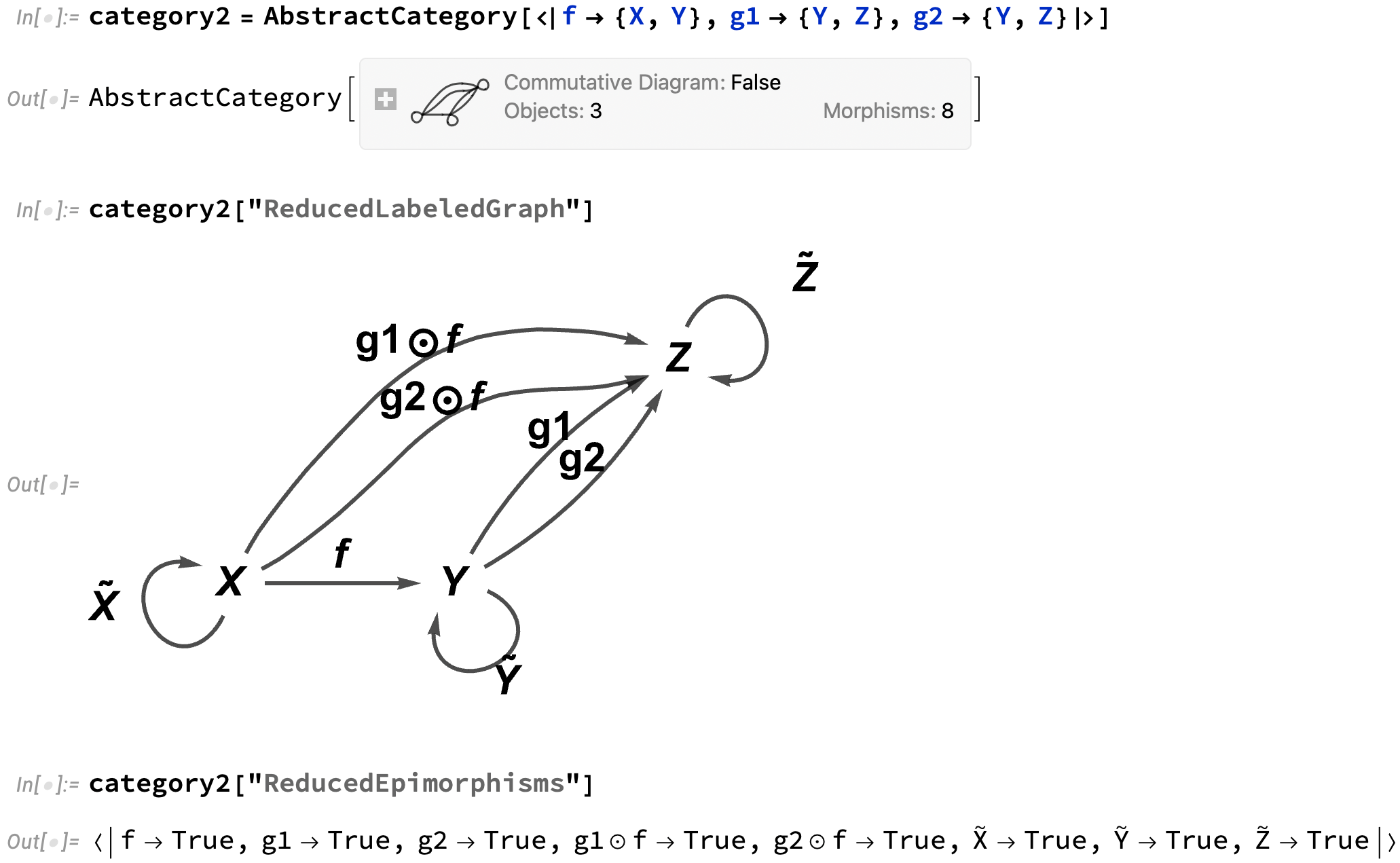}
\end{framed}
\caption{On the left, the \texttt{AbstractCategory} object corresponding to the basic diagrammatic setup of a monomorphism, showing that all morphisms are initially monomorphisms. On the right, the \texttt{AbstractCategory} object corresponding to the basic diagrammatic setup of an epimorphism, showing that all morphisms are initially epimorphisms.}
\label{fig:Figure6}
\end{figure}

\begin{figure}[ht]
\centering
\begin{framed}
\includegraphics[width=0.495\textwidth]{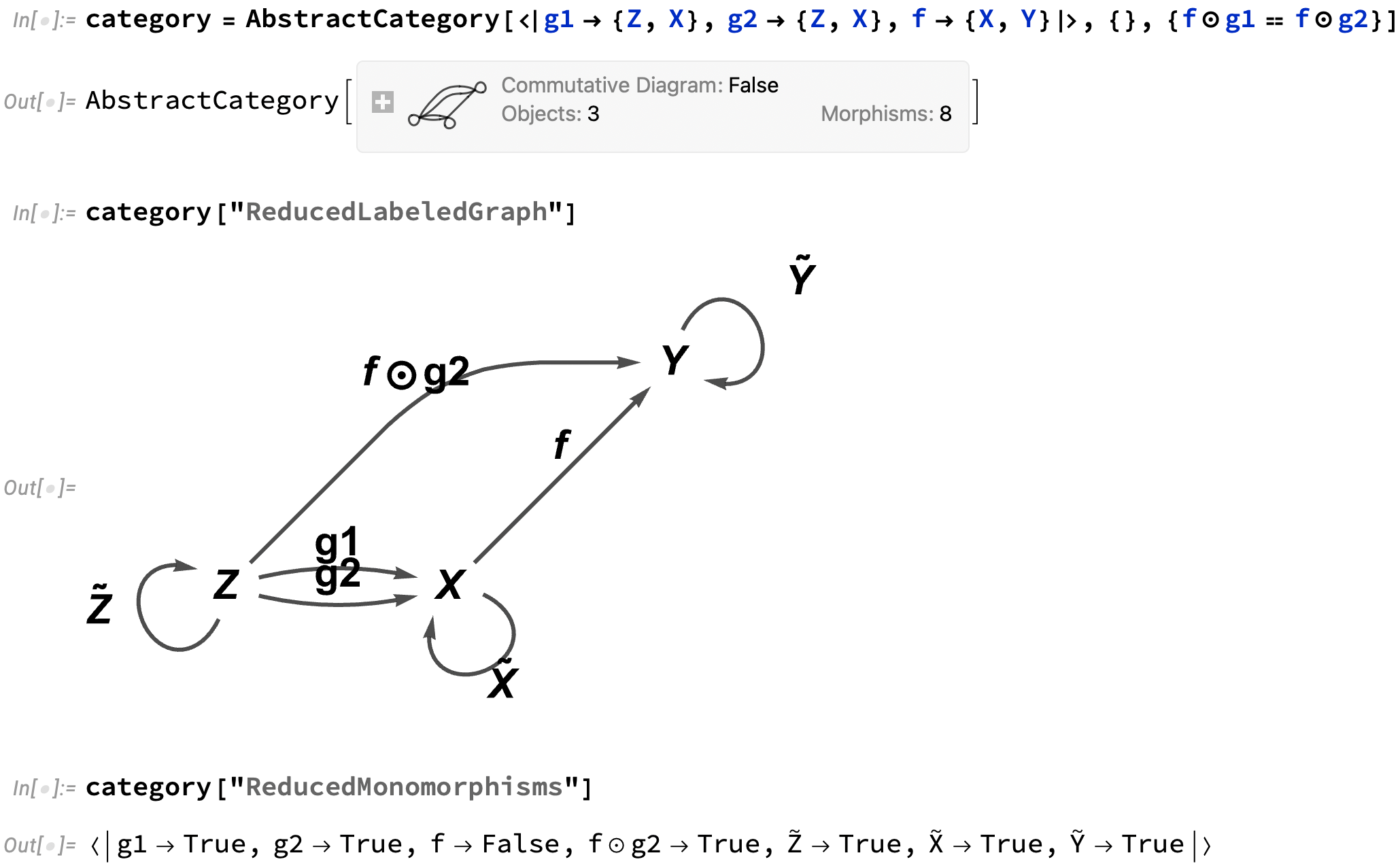}
\vrule
\includegraphics[width=0.495\textwidth]{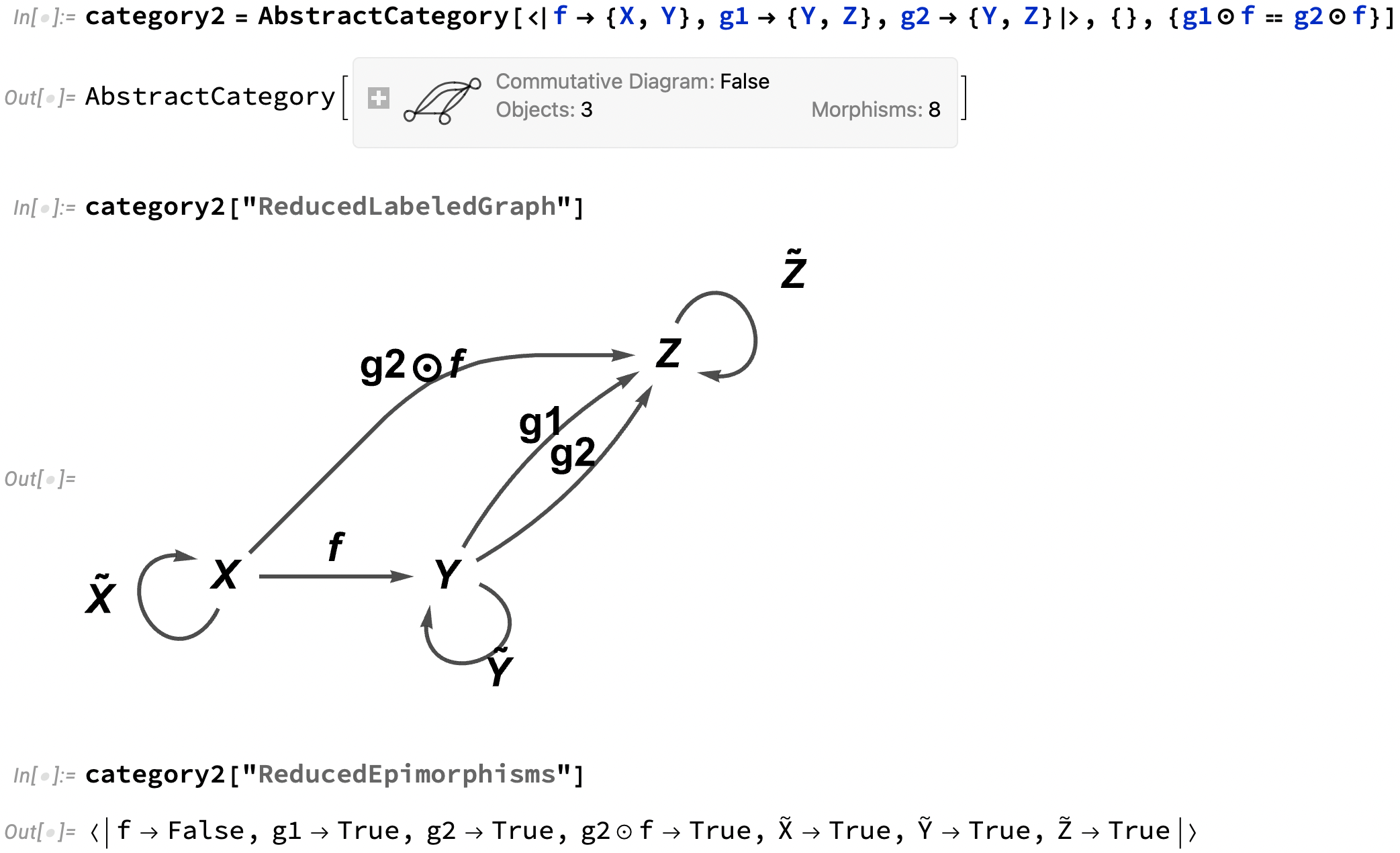}
\end{framed}
\caption{On the left, the \texttt{AbstractCategory} object corresponding to the basic diagrammatic setup of a monomorphism, with the additional algebraic equivalence ${f \circ g_1 = f \circ g_2}$ imposed, showing that morphism $f$ has now ceased to be a monomorphism. On the right, the \texttt{AbstractCategory} object corresponding to the basic diagrammatic setup of an epimorphism, with the additional algebraic equivalence ${g_1 \circ f = g_2 \circ f}$ imposed, showing that morphism $f$ has now ceased to be an epimorphism.}
\label{fig:Figure7}
\end{figure}

\begin{figure}[ht]
\centering
\begin{framed}
\includegraphics[width=0.495\textwidth]{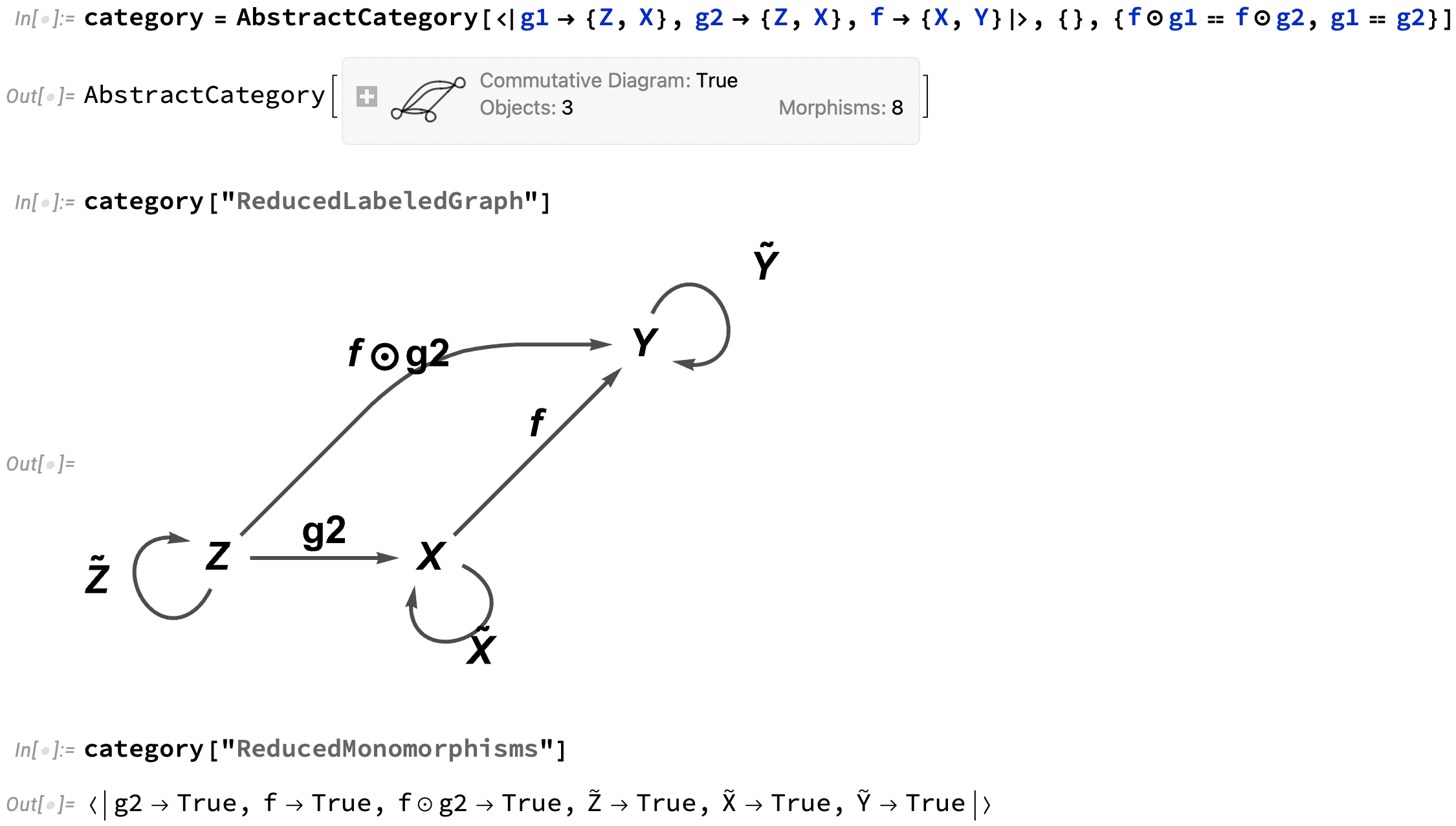}
\vrule
\includegraphics[width=0.495\textwidth]{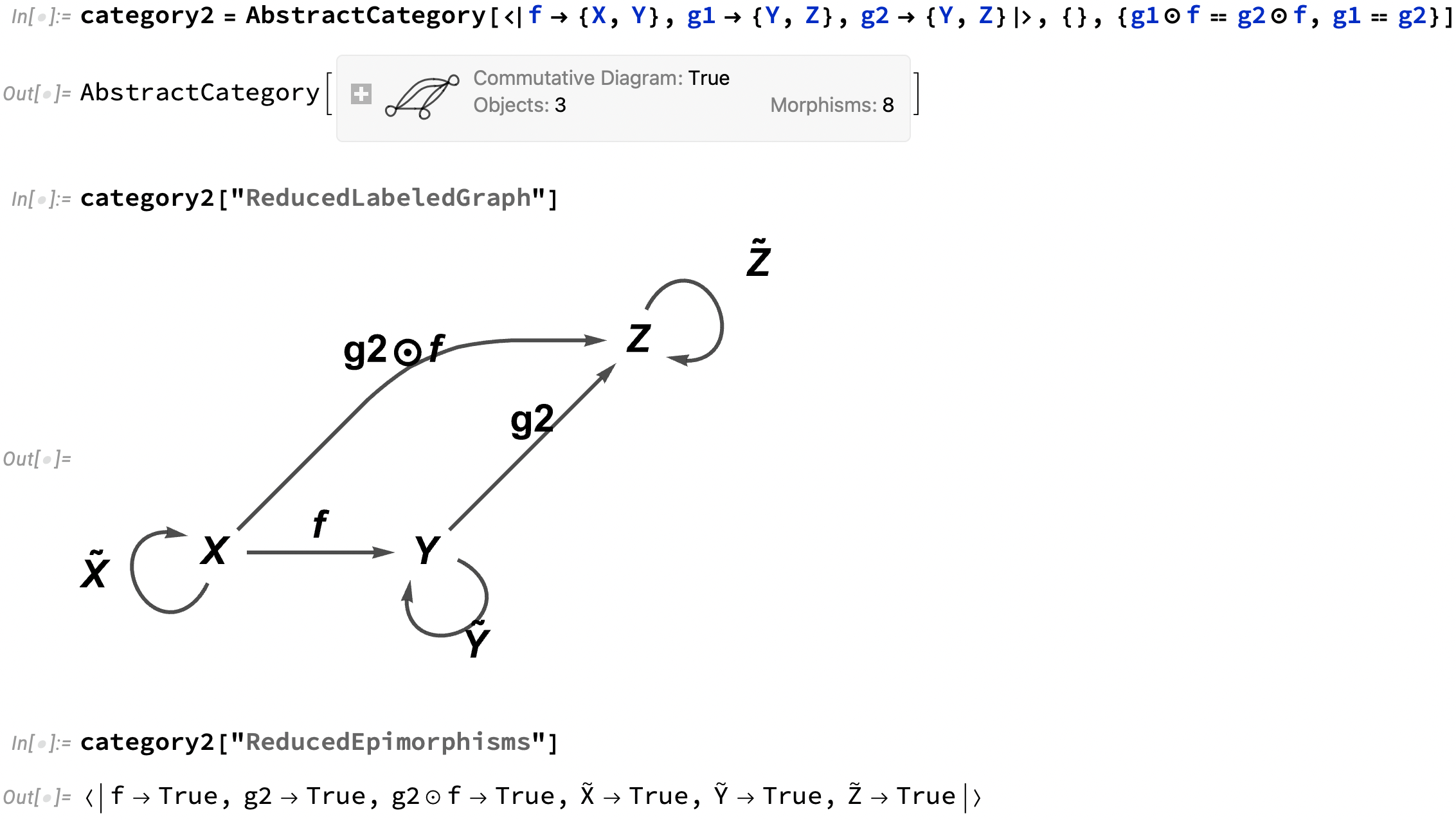}
\end{framed}
\caption{On the left, the \texttt{AbstractCategory} object corresponding to the basic diagrammatic setup of a monomorphism, with the additional algebraic equivalences ${f \circ g_1 = f \circ g_2}$ \textit{and} ${g_1 = g_2}$ both imposed, showing that morphism $f$ has now returned to being a monomorphism. On the right, the \texttt{AbstractCategory} object corresponding to the basic diagrammatic setup of an epimorphism, with the additional algebraic equivalences ${g_1 \circ f = g_2 \circ f}$ \textit{and} ${g_1 = g_2}$ both imposed, showing that morphism $f$ has now returned to being an epimorphism.}
\label{fig:Figure8}
\end{figure}

The category-theoretic analog of an injective function in mathematical analysis is a \textit{monomorphism}\cite{borceux}\cite{bergman}: a left-cancellative morphism. Specifically, if the morphism ${f : X \to Y}$ in the category ${\mathcal{C}}$ is such that, for all objects $Z$ and all pairs of morphisms ${g_1 : Z \to X}$ and ${g_2 : Z \to X}$ such that ${\left( f \circ g_1 : Z \to Y \right) = \left( f \circ g_2 : Z \to Y \right)}$, one is able to ``cancel on the left'' by ${f : X \to Y}$ to obtain that ${\left( g_1 : Z \to X \right ) = \left( g_2 : Z \to X \right)}$, then ${f : X \to Y}$ is a monomorphism:

\begin{multline}
\forall \left( f : X \to Y \right) \in \mathrm{hom} \left( \mathcal{C} \right), \qquad \left( f : X \to Y \right) \text{ is a monomorphism }\\
\iff \qquad \forall Z \in \mathrm{ob} \left( \mathcal{C} \right), \qquad \forall \left( g_1 : Z \to X \right), \left( g_2 : Z \to X \right) \in \mathrm{hom} \left( \mathcal{C} \right),\\
\left( f \circ g_1 : Z \to Y \right) = \left( f \circ g_2 : Z \to Y \right) \qquad \implies \qquad \left( g_1 : Z \to X \right) = \left( g_2 : Z \to X \right),
\end{multline}
or, illustrated diagrammatically, one has that if:

\begin{equation}
\begin{tikzcd}
& X \arrow[dr, "f"] &\\
Z \arrow[ur, swap, "g_1"] \arrow[ur, bend left, "g_2"] \arrow[rr, swap, "f \circ g_1"] \arrow[rr, bend right, swap, "f \circ g_2"] & & Y
\end{tikzcd} \qquad \mapsto \qquad
\begin{tikzcd}
& X \arrow[dr, "f"] &\\
Z \arrow[ur, swap, "g_1"] \arrow[ur, bend left, "g_2"] \arrow[rr, swap, "f \circ g_1 = f \circ g_2"] & & Y
\end{tikzcd},
\end{equation}
then one also necessarily has:

\begin{equation}
\begin{tikzcd}
& X \arrow[dr, "f"] &\\
Z \arrow[ur, "g_1 = g_2"] \arrow[rr, swap, "f \circ g_1 = f \circ g_2"] & & Y
\end{tikzcd}.
\end{equation}
The corresponding \textit{dual} notion to that of a monomorphism (i.e. the construction obtained by reversing the direction of the morphisms, and hence reversing the order of morphism composition, in the diagrams shown above) is that of an \textit{epimorphism}: a right-cancellative morphism, and hence the category-theoretic analog of a surjective function in  analysis. Specifically, if the morphism ${f : X \to Y}$ in the category ${\mathcal{C}}$ is such that, for all objects $Z$ and all pairs of morphisms ${g_1 : Y \to Z}$ and ${g_2 : Y \to Z}$ such that ${\left( g_1 \circ f : X \to Z \right) = \left( g_2 \circ f : X \to Z \right)}$, one is able to ``cancel on the right'' by ${f : X \to Y}$ to obtain that ${\left( g_1 : Y \to Z \right) = \left( g_2 : Y \to Z \right)}$, then ${f : X \to Y}$ is an epimorphism:

\begin{multline}
\forall \left( f : X \to Y \right) \in \mathrm{hom} \left( \mathcal{C} \right), \qquad \left( f : X \to Y \right) \text{ is an epimorphism }\\
\iff \qquad \forall Z \in \mathrm{ob} \left( \mathcal{C} \right), \qquad \forall \left( g_1 : Y \to Z \right), \left( g_2 : Y \to Z \right) \in \mathrm{hom} \left( \mathcal{C} \right),\\
\left( g_1 \circ f : X \to Z \right) = \left( g_2 \circ f : X \to Z \right) \qquad \implies \qquad \left( g_1 : Y \to Z \right) = \left( g_2 : Y \to Z \right),
\end{multline}
or, illustrated diagrammatically, one has that if:

\begin{equation}
\begin{tikzcd}
& Y \arrow[dr, swap, "g_1"] \arrow[dr, bend left, "g_2"] &\\
X \arrow[ur, "f"] \arrow[rr, swap, "g_1 \circ f"] \arrow[rr, bend right, swap, "g_2 \circ f"] & & Z
\end{tikzcd} \qquad \mapsto \qquad
\begin{tikzcd}
& Y \arrow[dr, swap, "g_1"] \arrow[dr, bend left, "g_2"] &\\
X \arrow[ur, "f"] \arrow[rr, swap, "g_1 \circ f = g_2 \circ f"] & & Z
\end{tikzcd},
\end{equation}
then one also necessarily has:

\begin{equation}
\begin{tikzcd}
& Y \arrow[dr, "g_1 = g_2"] &\\
X \arrow[ur, "f"] \arrow[rr, swap, "g_1 \circ f = g_2 \circ f"] & & Z
\end{tikzcd}.
\end{equation}
Figure \ref{fig:Figure6} illustrates the basic diagrammatic setup for both monomorphisms and epimorphisms (commonly abbreviated simply to \textit{monos} and \textit{epis} in the category-theoretic literature) in \textsc{Categorica}, and demonstrates in these two cases that all morphisms initially, i.e. in the absence of any further algebraic equivalences, correspond to monomorphisms and epimorphisms, respectively. Figure \ref{fig:Figure7} shows that, by imposing the algebraic equivalence ${\left( f \circ g_1 : Z \to Y \right) = \left( f \circ g_2 : Z \to Y \right)}$ in the former case and ${\left( g_1 \circ f : X \to Z \right) = \left( g_2 \circ f : X \to Z \right)}$ in the latter case, one is able to force the morphism ${f : X \to Y}$ to cease to be a monomorphism in the former example, and to cease to be an epimorphism in the latter example. Finally, Figure \ref{fig:Figure8} demonstrates that one is able to restore the status of morphism ${f : X \to Y}$ as a monomorphism (in the former case) or an epimorphism (in the latter case) by imposing the additional algebraic equivalence ${\left( g_1 : Z \to X \right) = \left( g_2 : Z \to X \right)}$ (for the monomorphism case) or ${\left( g_1 : Y \to Z \right) = \left( g_2 : Y \to Z \right)}$ (for the epimorphism case). Note that any morphism that is both a monomorphism and an epimorphism is known as a \textit{bimorphism} (the category-theoretic analog of a bijective function), and \textsc{Categorica} contains in-built functionality for detecting and manipulating bimorphisms in much the same way.

\begin{figure}[ht]
\centering
\begin{framed}
\includegraphics[width=0.495\textwidth]{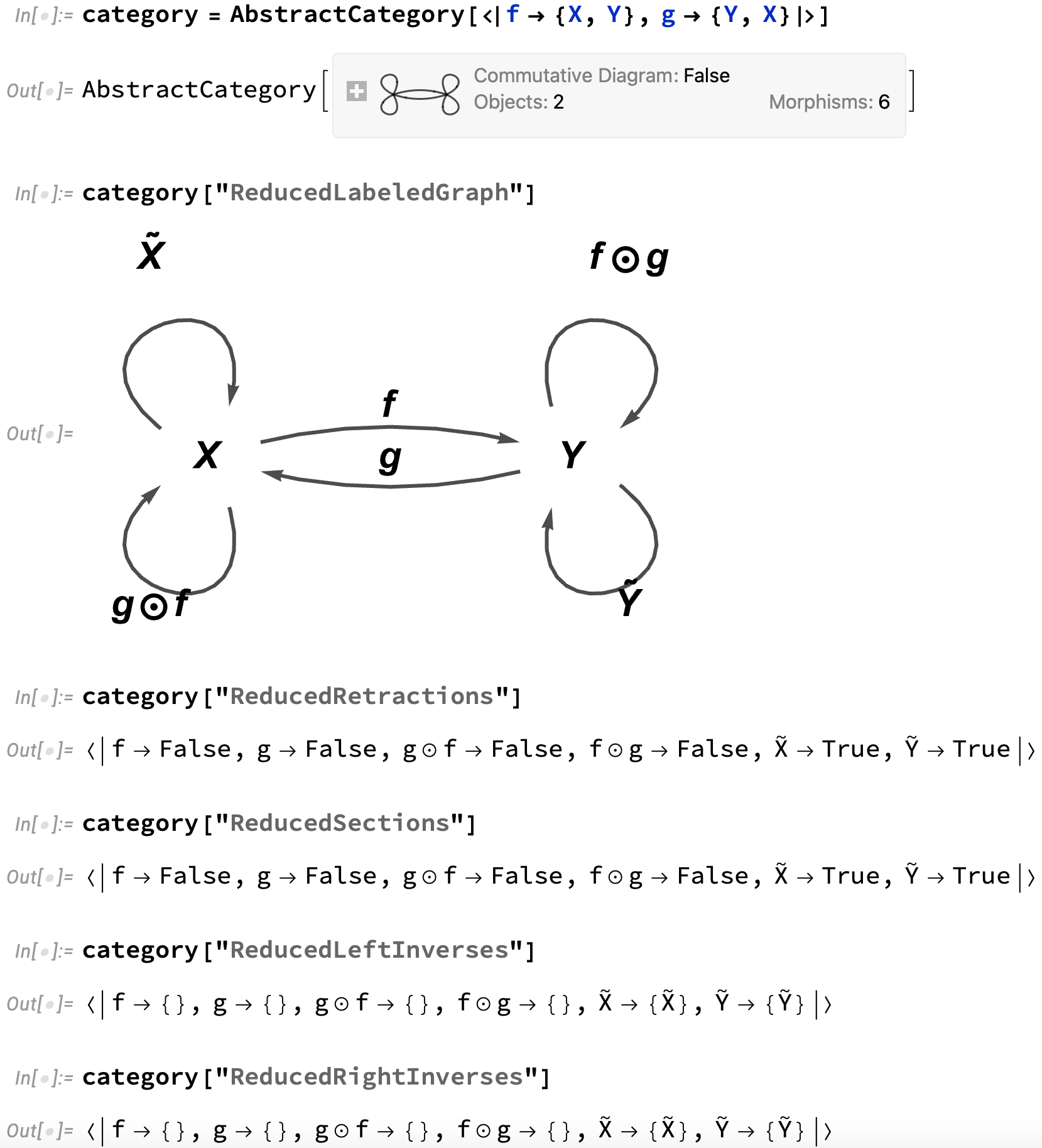}
\vrule
\includegraphics[width=0.495\textwidth]{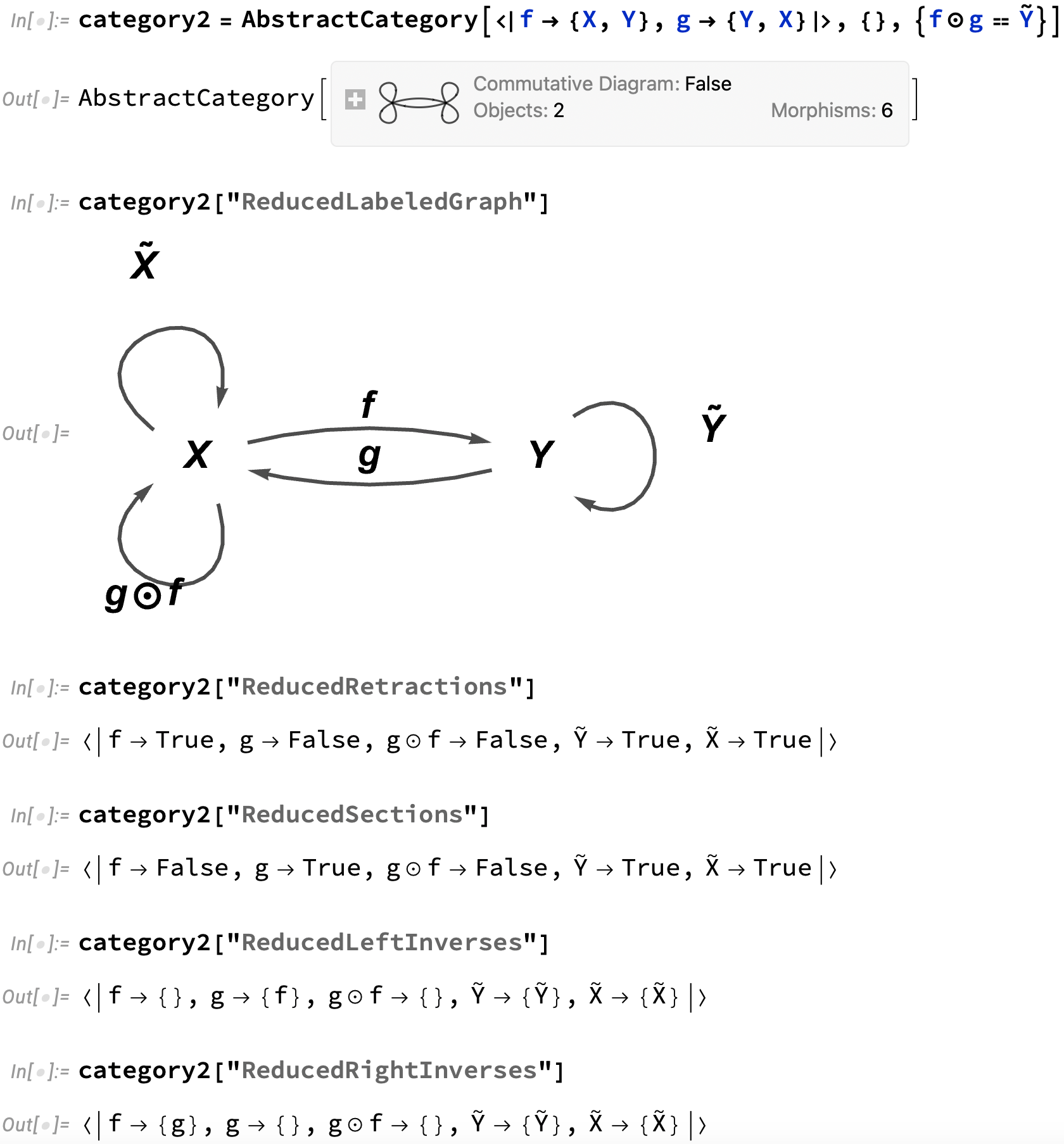}
\end{framed}
\caption{On the left, the \texttt{AbstractCategory} object corresponding to the basic diagrammatic setup of a retraction/section, showing that initially no morphisms (except for identity morphisms) are either retractions or sections, and therefore that no morphisms (except for identity morphisms) possess either left or right inverses. On the right, the \texttt{AbstractCategory} object corresponding to the basic diagrammatic setup of a retraction/section, with the additional algebraic equivalence ${f \circ g = id_Y}$ imposed, showing that morphism $f$ has now become a retraction, and morphism $g$ has now become a section, with $f$ being the left inverse of $g$ and $g$ being the right inverse of $f$.}
\label{fig:Figure9}
\end{figure}

An important special case of monomorphisms are \textit{sections}\cite{mitchell}\cite{maclane}: morphisms that possess left inverses (and hence which are, themselves, right inverses of some other morphism). Likewise, an important special case of epimorphisms are \textit{retractions}: dual to sections, these are morphisms that possess right inverses (and hence which are, themselves, left inverses of some other morphism). Specifically, if the morphisms ${f : X \to Y}$ and ${g : Y \to X}$ in the category ${\mathcal{C}}$ are such that ${f \circ g : Y \to Y}$ is the identity morphism on $Y$ (i.e. ${\left( f \circ g : Y \to Y \right) = \left( id_Y : Y \to Y \right)}$), then ${f: X \to Y}$ is a retraction of ${g : Y \to X}$, and ${g : Y \to X}$ is a section of ${f : X \to Y}$:

\begin{multline}
\forall \left( f : X \to Y \right) \in \mathrm{hom} \left( \mathcal{C} \right), \qquad \left( f : X \to Y \right) \text{ is a retraction }\\
\iff \qquad \exists \left( g : Y \to X \right) \in \mathrm{hom} \left( \mathcal{C} \right), \qquad \text{ such that } \qquad \left( f \circ g : Y \to Y \right) = \left( id_Y : Y \to Y \right),
\end{multline}
or, illustrated diagrammatically, a retraction is a morphism ${f : X \to Y}$ such that one has:

\begin{equation}
\begin{tikzcd}
X \arrow[rr, bend left, "f"] \arrow[loop above, "id_X"] \arrow[loop below, "g \circ f"] & & Y \arrow[ll, bend left, "g"] \arrow[loop above, "id_Y"] \arrow[loop below, "f \circ g"]
\end{tikzcd} \qquad \mapsto \qquad
\begin{tikzcd}
X \arrow[rr, bend left, "f"] \arrow[loop above, "id_X"] \arrow[loop below, "g \circ f"] & & Y \arrow[ll, bend left, "g"] \arrow[loop right, "f \circ g = id_Y"]
\end{tikzcd}.
\end{equation}
Similarly, if ${g \circ f : X \to X}$ is the identity morphism on $X$ (i.e. ${\left( g \circ f : X \to X \right) = \left( id_X : X \to X \right)}$) then ${f : X \to Y}$ is a section of ${g : Y \to X}$, and ${g : Y \to X}$ is a retraction of ${f : X \to Y}$:

\begin{multline}
\forall \left( f : X \to Y \right) \in \mathrm{hom} \left( \mathcal{C} \right), \qquad \left( f : X \to Y \right) \text{ is a section }\\
\iff \qquad \exists \left( g : Y \to X \right) \in \mathrm{hom} \left( \mathcal{C} \right), \qquad \text{ such that } \qquad \left( g \circ f : X \to X \right) = \left( id_X : X \to X \right),
\end{multline}
or, illustrated diagrammatically, a section is a morphism ${f : X \to Y}$ such that one has:

\begin{equation}
\begin{tikzcd}
X \arrow[rr, bend left, "f"] \arrow[loop above, "id_X"] \arrow[loop below, "g \circ f"] & & Y \arrow[ll, bend left, "g"] \arrow[loop above, "id_Y"] \arrow[loop below, "f \circ g"]
\end{tikzcd} \qquad \mapsto \qquad
\begin{tikzcd}
X \arrow[rr, bend left, "f"] \arrow[loop left, "g \circ f = id_X"] & & Y \arrow[ll, bend left, "g"] \arrow[loop above, "id_Y"] \arrow[loop below, "f \circ g"]
\end{tikzcd}.
\end{equation}
These category-theoretic notions of retractions and sections are so-named because they naturally generalize the corresponding notions in topology, wherein one considers \textit{retractions} of topological spaces into subspaces to be continuous maps that preserve all points in the subspace (and for which the corresponding inclusion maps of the subspaces into the original spaces, such that the compositions of the two maps always reduce to the identity maps on the subspaces, would be \textit{sections})\cite{borsuk}\cite{eilenberg3}. Figure \ref{fig:Figure9} illustrates the basic diagrammatic setup for both retractions and sections in \textsc{Categorica}, and demonstrates that initially (in the absence of any further algebraic equivalences) no morphisms, with the exception of identity morphisms, are either retractions or sections, and therefore that no morphisms, with the exception of identity morphisms, possess either left or right inverses. However, by imposing the algebraic equivalence ${\left( f \circ g : Y \to Y \right) = \left( id_Y : Y \to Y \right)}$, one is able to force the morphism ${f : X \to Y}$ to be a retraction (and hence to possess the morphism ${g : Y \to X}$ as its right inverse), as well as to force the morphism ${g : Y \to X}$ to be a section (and hence to possess the morphism ${f : X \to Y}$ as its left inverse). It is therefore straightforward to prove, using \textsc{Categorica}'s diagrammatic theorem-proving capabilities, that every section is necessarily a monomorphism (i.e. that the existence of a left inverse necessarily implies left-cancellativity) and that every retraction is necessarily an epimorphism (i.e. that the existence of a right inverse necessarily implies right-cancellativity).

\begin{figure}[ht]
\centering
\begin{framed}
\includegraphics[width=0.445\textwidth]{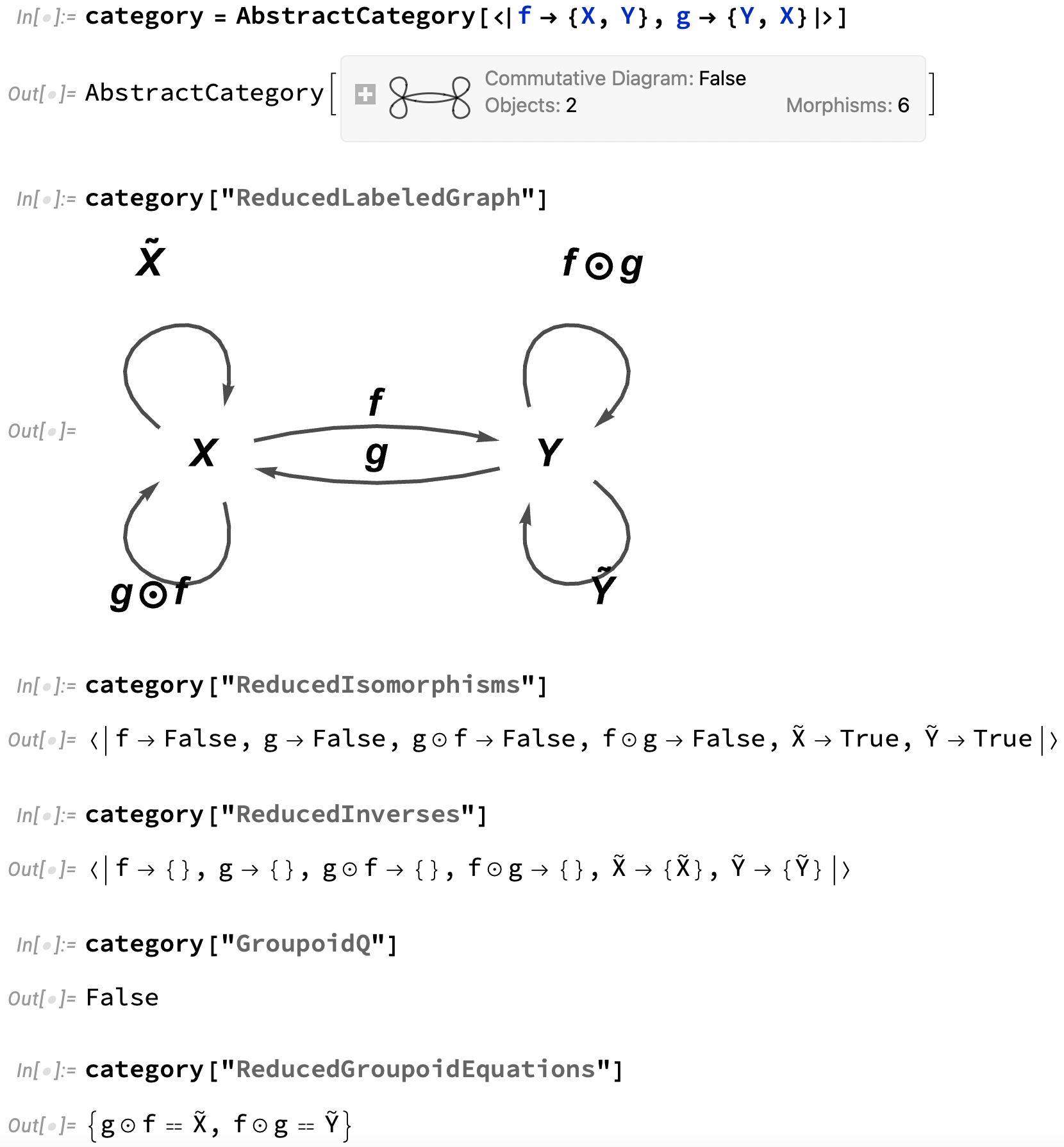}
\vrule
\includegraphics[width=0.545\textwidth]{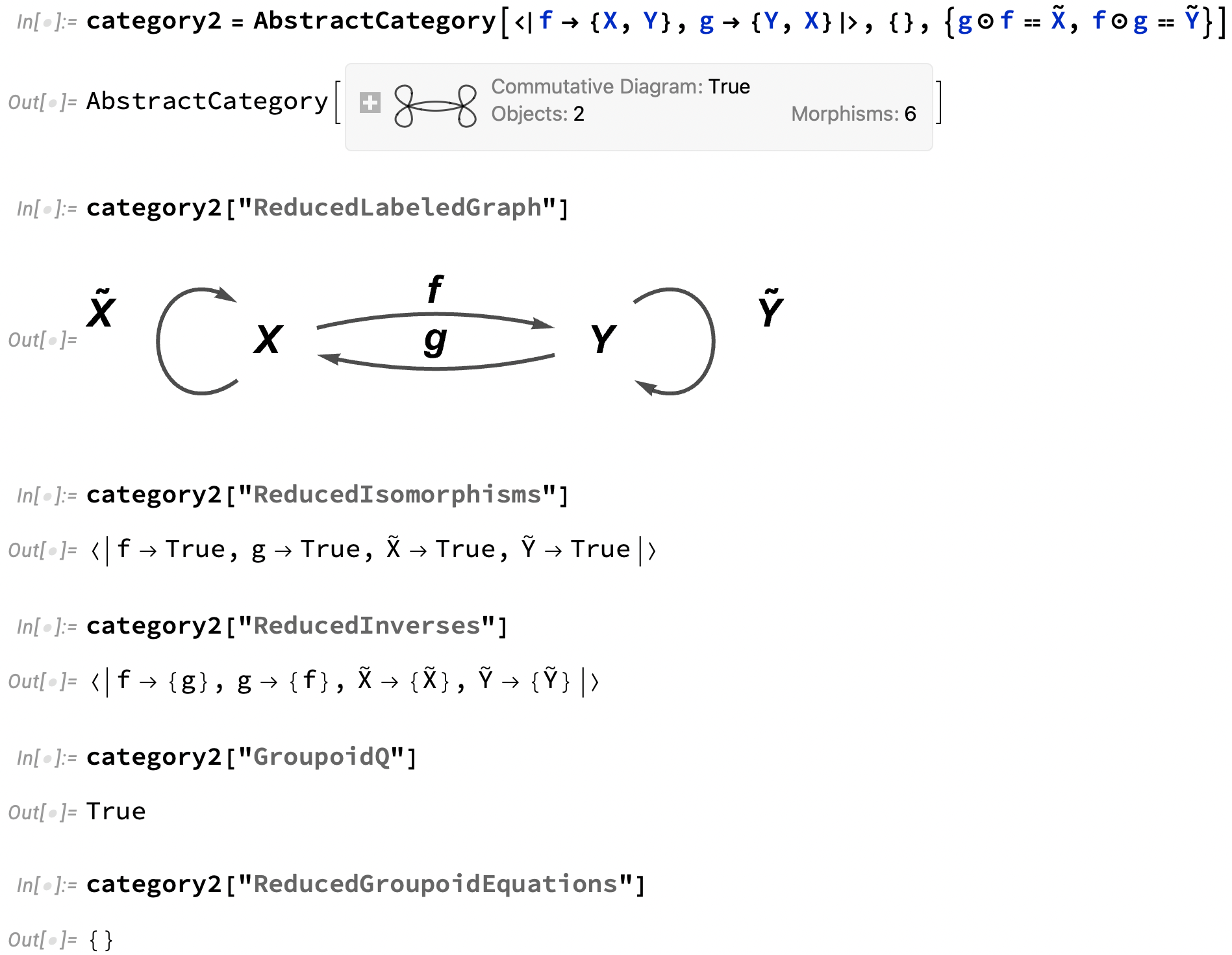}
\end{framed}
\caption{On the left, the \texttt{AbstractCategory} object corresponding to a simple (not yet groupoidal) category, showing that initially no morphisms (except for identity morphisms) are isomorphisms, and therefore that no morphisms (except for identity isomorphisms) possess inverses, hence illustrating that the category is not a groupoid and showing that the morphism equivalences ${g \circ f = id_X}$ and ${f \circ g = id_Y}$ are the minimal algebraic conditions necessary to force the category to be groupoidal. On the right, the \texttt{AbstractCategory} object for the groupoidal case of the same category, with the morphism equivalences ${g \circ f = id_X}$ and ${f \circ g = id_Y}$ imposed, showing that the morphisms $f$ and $g$ have now become isomorphisms, with $f$ being the inverse of $g$ and $g$ being the inverse of $f$, hence demonstrating that these equivalences are indeed sufficient to force the category to be groupoidal.}
\label{fig:Figure10}
\end{figure}

\begin{figure}[ht]
\centering
\begin{framed}
\includegraphics[width=0.445\textwidth]{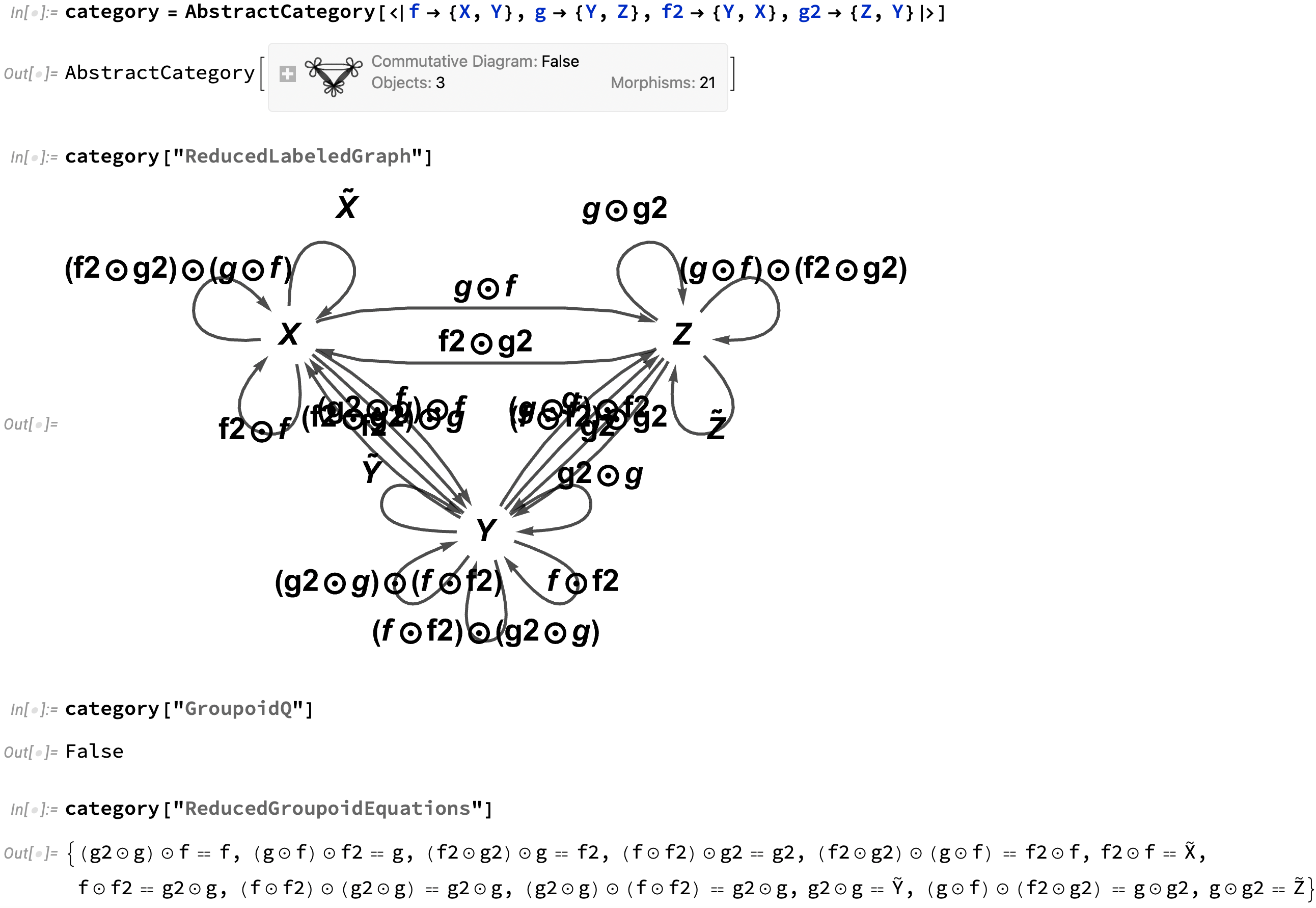}
\vrule
\includegraphics[width=0.545\textwidth]{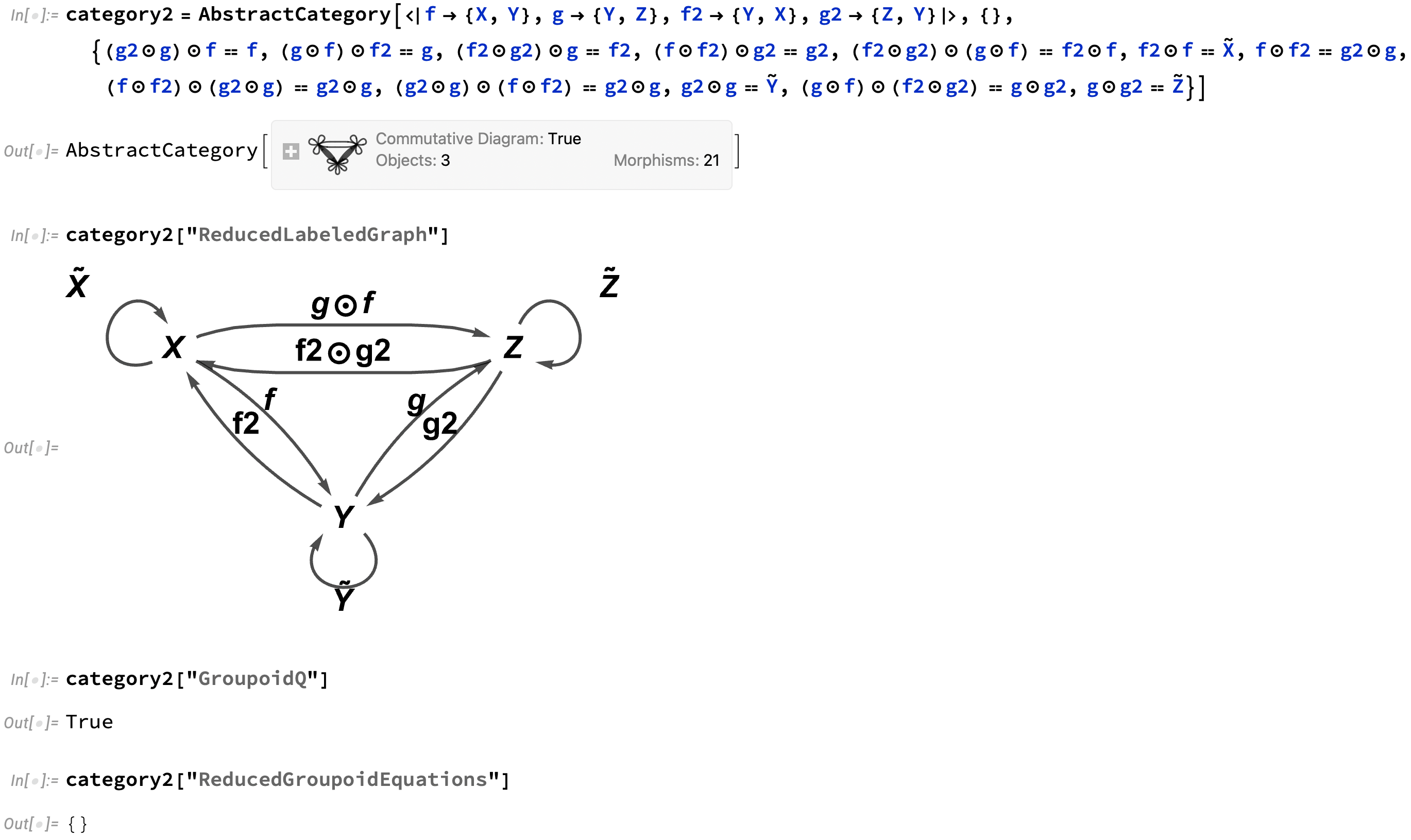}
\end{framed}
\caption{On the left, the \texttt{AbstractCategory} object corresponding to a slightly more complex (not yet groupoidal) category, illustrating that the category is not a groupoid and computing the minimum set of morphism equivalences necessary to force the category to be groupoidal. On the right, the \texttt{AbstractCategory} object for the groupoidal case of the same category, with the aforementioned morphism equivalences imposed, demonstrating that these equivalences are indeed sufficient to force the category to be groupoidal.}
\label{fig:Figure11}
\end{figure}

A morphism that acts as both a retraction and a section (and therefore which both is, and also possesses, a left and right inverse) is known as an \textit{isomorphism}. Specifically, if the morphism ${f : X \to Y}$ in the category ${\mathcal{C}}$ is such that there exists another morphism ${g : Y \to X}$ such that both ${g \circ f : X \to X}$ is the identity morphism on $X$ (i.e. ${\left( g \circ f : X \to X \right) = \left( id_X : X \to X \right)}$) and ${f \circ g : Y \to Y}$ is the identity morphism on $Y$ (i.e. ${\left( f \circ g : Y \to Y \right) = \left( id_Y : Y \to Y \right)}$, then ${f : X \to Y}$ is an isomorphism (as, correspondingly, is ${g : Y \to X}$):

\begin{multline}
\forall \left( f : X \to Y \right) \in \mathrm{hom} \left( \mathcal{C} \right), \qquad \left( f : X \to Y \right) \text{ is an isomorphism }\\
\iff \qquad \exists \left( g : Y \to X \right) \in \mathrm{hom} \left( \mathcal{C} \right), \qquad \text{ such that } \qquad \left( g \circ f : X \to X \right) = \left( id_X : X \to X \right)\\
\text{ and } \qquad \left( f \circ g : Y \to Y \right) = \left( id_Y : Y \to Y \right),
\end{multline}
or, illustrated diagrammatically, an isomorphism is a morphism ${f : X \to Y}$ such that one has:

\begin{equation}
\begin{tikzcd}
X \arrow[rr, bend left, "f"] \arrow[loop above, "id_X"] \arrow[loop below, "g \circ f"] & & Y \arrow[ll, bend left, "g"] \arrow[loop above, "id_Y"] \arrow[loop below, "f \circ g"]
\end{tikzcd} \qquad \mapsto \qquad
\begin{tikzcd}
X \arrow[rr, bend left, "f"] \arrow[loop left, "g \circ f = id_X"] & & Y \arrow[ll, bend left, "g"] \arrow[loop right, "f \circ g = id_Y"]
\end{tikzcd}.
\end{equation}
A category in which all morphisms are isomorphisms is known as a \textit{groupoid}, since the morphisms of a groupoid consisting of a single object trivially form a group under the operation of morphism composition\cite{dicks} (with the associativity and identity axioms deriving from the underlying axioms of category theory, and the inverse axiom deriving from the existence of both left and right inverses for each morphism), and therefore the morphisms of a groupoid consisting of multiple objects may be thought of as forming an abstract generalization of a group in which the group operation is \textit{partial}, i.e. only well-defined for particular pairs of elements. Groupoids are particularly widely-studied in algebraic topology and homotopy theory, wherein the \textit{fundamental groupoid} of a topological space generalizes the more traditional fundamental group to the case where one does not necessarily fix a single distinguished base point\cite{brown}. Much like in the case of commutative diagrams discussed previously, \textsc{Categorica} is able to use purely graph-theoretic algorithms to compute the minimum set of morphism equivalences necessary to force the category shown above to be a groupoid, as well as to prove that these equivalences are indeed sufficient to do so, as shown in Figure \ref{fig:Figure10}. However, once the initial category is even marginally more complex than this, the task of determining the minimum set of algebraic conditions necessary for the category to be groupoidal quickly evolves to be highly non-trivial; for instance, for the case of a three-object category of the form:

\begin{equation}
\begin{tikzcd}
& Y \arrow[dr, swap, "g"] \arrow[dl, bend right, swap, "f^{-1}"] &\\
X \arrow[ur, swap, "f"] & & Z \arrow[ul, bend right, swap, "g^{-1}"]
\end{tikzcd},
\end{equation}
one must consider not only the conditions necessary to force the morphisms ${f : X \to Y}$ and ${g : Y \to Z}$ (and hence also the morphisms ${f^{-1} : Y \to X}$ and ${g^{-1} : Z \to Y}$) to be isomorphisms, namely:

\begin{equation}
\begin{tikzcd}
& Y \arrow[dr, swap, "g"] \arrow[dl, bend right, swap, "f^{-1}"] \arrow[loop above, "f \circ f^{-1}"] \arrow[loop below, "id_Y"] &\\
X \arrow[ur, swap, "f"] \arrow[loop left, "id_X"] \arrow[loop below, "f^{-1} \circ f"] & & Z \arrow[ul, bend right, swap, "g^{-1}"]
\end{tikzcd} \qquad \mapsto \qquad
\begin{tikzcd}
& Y \arrow[dr, swap, "g"] \arrow[dl, bend right, swap, "f^{-1}"] \arrow[loop above, "f \circ f^{-1} = id_Y"] &\\
X \arrow[ur, swap, "f"] \arrow[loop left, "f^{-1} \circ f = id_X"] & & Z \arrow[ul, bend right, swap, "g^{-1}"]
\end{tikzcd},
\end{equation}
i.e:

\begin{equation}
\left( f^{-1} \circ f : X \to X \right) = \left( id_X : X \to X \right), \qquad \text{ and } \qquad \left( f \circ f^{-1} : Y \to Y \right) = \left( id_Y : Y \to Y \right),
\end{equation}
and:

\begin{equation}
\begin{tikzcd}
& Y \arrow[dr, swap, "g"] \arrow[dl, bend right, swap, "f^{-1}"] \arrow[loop above, "g^{-1} \circ g"] \arrow[loop below, "id_Y"] &\\
X \arrow[ur, swap, "f"] & & Z \arrow[ul, bend right, swap, "g^{-1}"] \arrow[loop right, "id_Z"] \arrow[loop below, "g \circ g^{-1}"]
\end{tikzcd} \qquad \mapsto \qquad
\begin{tikzcd}
& Y \arrow[dr, swap, "g"] \arrow[dl, bend right, swap, "f^{-1}"] \arrow[loop above, "g^{-1} \circ g = id_Y"] &\\
X \arrow[ur, swap, "f"] & & Z \arrow[ul, bend right, swap, "g^{-1}"] \arrow[loop right, "g \circ g^{-1} = id_Z"]
\end{tikzcd},
\end{equation}
i.e:

\begin{equation}
\left( g^{-1} \circ g : Y \to Y \right) = \left( id_Y : Y \to Y \right), \qquad \text{ and } \qquad \left( g \circ g^{-1} : Z \to Z \right) = \left( id_Z : Z \to Z \right),
\end{equation}
respectively, but also the conditions necessary to force the resulting pair of composite morphisms ${g \circ f : X \to Z}$ and ${f^{-1} \circ g^{-1} : Z \to X}$ to be isomorphisms as well, namely to collapse:

\begin{equation}
\begin{tikzcd}
& Y \arrow[dr, swap, "g"] \arrow[dl, bend right, swap, "f^{-1}"] &\\
X \arrow[ur, swap, "f"] \arrow[rr, "g \circ f"] \arrow[loop left, "\left( f^{-1} \circ g^{-1} \right) \circ \left( g \circ f \right)"] \arrow[loop below, "id_X"] & & Z \arrow[ul, bend right, swap, "g^{-1}"] \arrow[ll, bend left, "f^{-1} \circ g^{-1}"] \arrow[loop right, "\left( g \circ f \right) \circ \left( f^{-1} \circ g^{-1} \right)"] \arrow[loop below, "id_Z"]
\end{tikzcd},
\end{equation}
down to:

\begin{equation}
\mapsto \qquad \begin{tikzcd}
& Y \arrow[dr, swap, "g"] \arrow[dl, bend right, swap, "f^{-1}"] &\\
X \arrow[ur, swap, "f"] \arrow[rr, "g \circ f"] \arrow[loop left, "\left( f^{-1} \circ g^{-1} \right) \circ \left( g \circ f \right) = id_X"] & & Z \arrow[ul, bend right, swap, "g^{-1}"] \arrow[ll, bend left, "f^{-1} \circ g^{-1}"] \arrow[loop right, "\left( g \circ f \right) \circ \left( f^{-1} \circ g^{-1} \right) = id_Z"]
\end{tikzcd},
\end{equation}
in addition to any relevant permutations thereof. Nevertheless, the generality of \textsc{Categorica}'s algebraic reasoning algorithms ensures that it is able to handle such cases in exactly the same way, as demonstrated in Figure \ref{fig:Figure11}.

\begin{figure}[ht]
\centering
\begin{framed}
\includegraphics[width=0.495\textwidth]{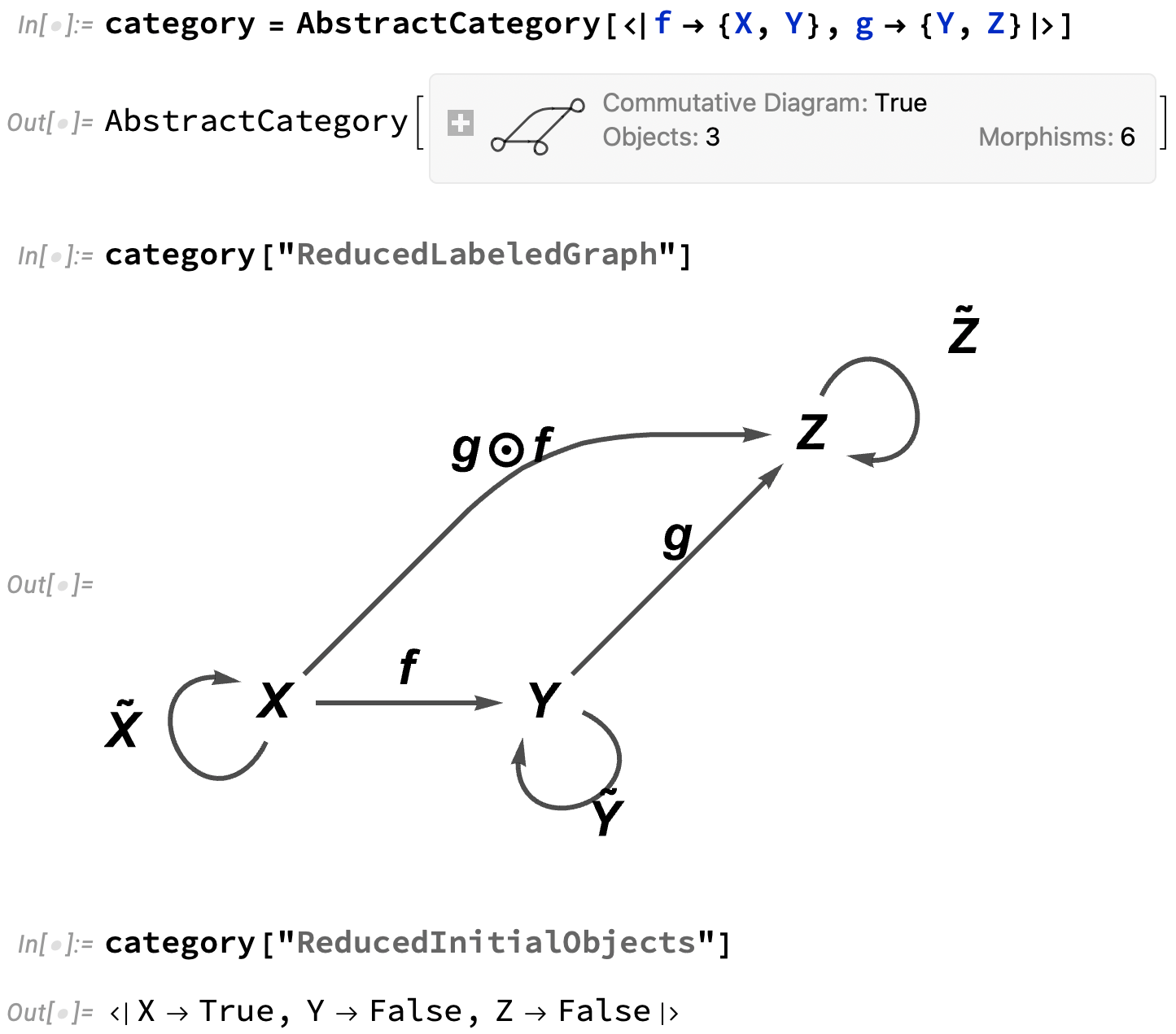}
\vrule
\includegraphics[width=0.495\textwidth]{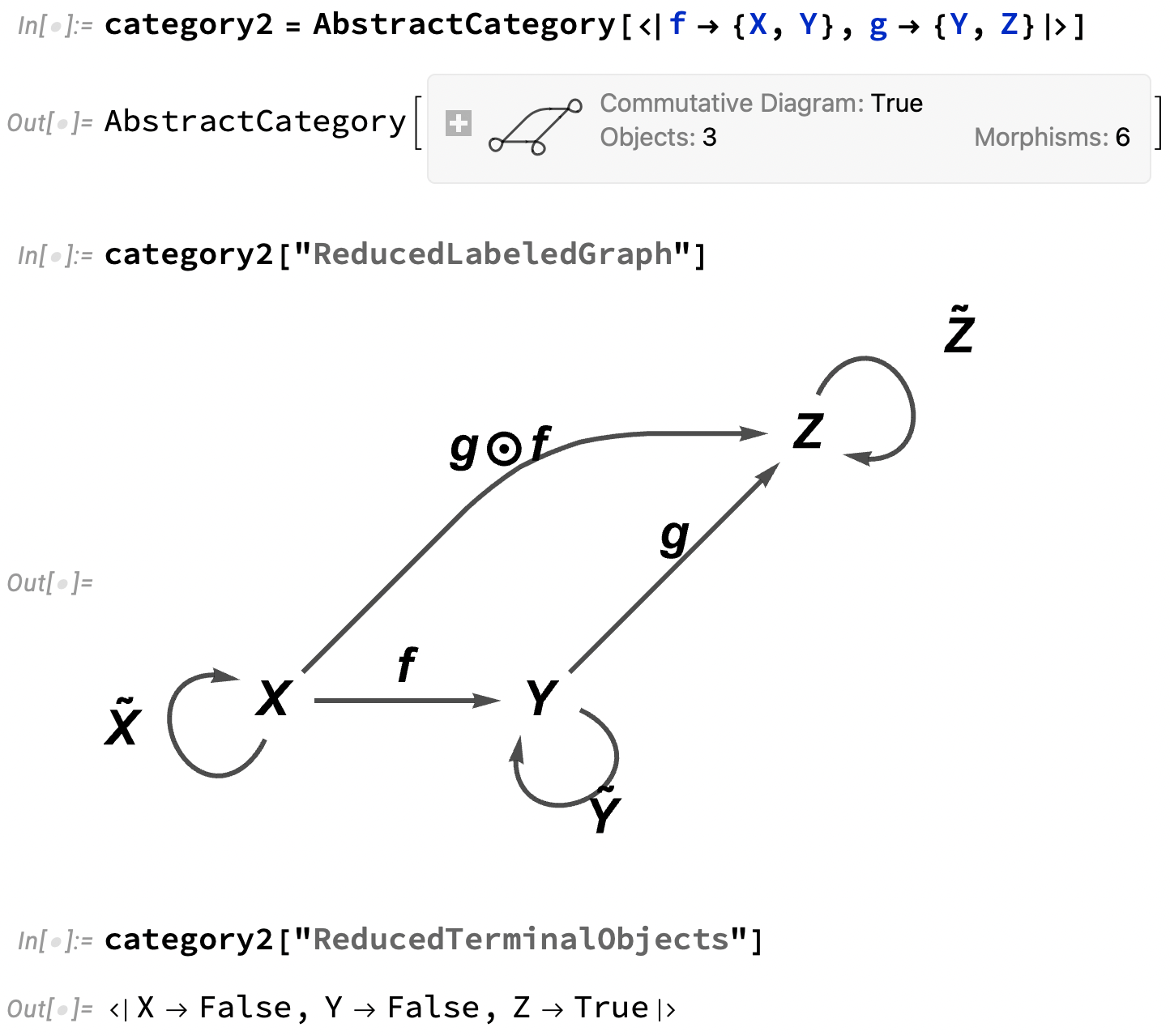}
\end{framed}
\caption{On the left, the \texttt{AbstractCategory} object corresponding to a simple triangular diagram, showing that object $X$ is an initial object. On the right, the same \texttt{AbstractCategory} object corresponding to the same simple triangular diagram, showing that object $Z$ is a terminal object.}
\label{fig:Figure20}
\end{figure}

Dual constructions of this general kind, such as monomorphisms vs. epimorphisms, retractions vs. sections, etc., are ubiquitous throughout category theory, and can be investigated systematically using the \textit{``DualCategory''} property of \texttt{AbstractCategory} objects in the \textsc{Categorica} framework (e.g. a morphism that registers as a monomorphism within a particular \texttt{AbstractCategory} object will register as an epimorphism in the corresponding \texttt{AbstractCategory} object returned by the \textit{``DualCategory''} property, etc.). A very straightforward example of a dual construction is that of \textit{initial} vs. \textit{terminal} objects; an initial object $X$ in the category ${\mathcal{C}}$ is an object such that, for every object $P$ (including $X$ itself) in ${\mathcal{C}}$, there exists a unique \textit{outgoing} morphism ${f : X \to P}$:

\begin{equation}
\forall X \in \mathrm{ob} \left( \mathcal{C} \right), \qquad X \text{ is an initial object } \qquad \iff \qquad \forall P \in \mathrm{ob} \left( \mathcal{C} \right), \qquad \exists! \left( f : X \to P \right) \in \mathrm{hom} \left( \mathcal{C} \right),
\end{equation}
or, illustrated diagrammatically:

\begin{equation}
\begin{tikzcd}
X \arrow[rr, dashed, "\exists! f"] & & \forall P
\end{tikzcd}.
\end{equation}
Dually, a terminal object $X$ in the category ${\mathcal{C}}$ is an object such that, for every object $P$ (including $X$ itself) in ${\mathcal{C}}$, there exists a unique \textit{incoming} morphism ${f : P \to X}$:

\begin{multline}
\forall X \in \mathrm{ob} \left( \mathcal{C} \right), \qquad X \text{ is a terminal object }\\
\iff \qquad \forall P \in \mathrm{ob} \left( \mathcal{C} \right), \qquad \exists! \left( f : P \to X \right) \in \mathrm{hom} \left( \mathcal{C} \right),
\end{multline}
or, illustrated diagrammatically:

\begin{equation}
\begin{tikzcd}
\forall P \arrow[rr, dashed, "\exists! f"] & & X
\end{tikzcd}.
\end{equation}
Initial and terminal objects generalize many key construction in pure mathematics\cite{pedicchio}\cite{maclane}, such as the bottom/minimal and top/maximal elements ${\bot}$ and ${\top}$ of partially-ordered sets (since, if one considers a category ${\mathcal{C}}$ constructed from a given poset ${\mathcal{P}}$ whose objects are elements of ${\mathcal{P}}$ and where the morphism ${f : X \to Y}$ exists in ${\mathcal{C}}$ if and only if ${X \leq Y}$ in ${\mathcal{P}}$, then the bottom/minimal element ${\bot}$ is an initial object and the top/maximal element ${\top}$ is a terminal object in the category ${\mathcal{C}}$), and the empty and point spaces ${\varnothing}$ and ${*}$ in topology (since, if one considers the category ${\mathbf{Top}}$ whose objects are topological spaces and whose morphisms are continuous functions between them, then the empty space ${\varnothing}$ is an initial object in ${\mathbf{Top}}$ because there exists a unique continuous function mapping the empty space to any other topological space, and the point space ${*}$ is a terminal object in ${\mathbf{Top}}$ because there exists a unique continuous function mapping any topological space to the point space). If we take the simple triangular diagram case that we have investigated previously:

\begin{equation}
\begin{tikzcd}
& Y \arrow[dr, "g"] \arrow[loop above, "id_Y"] &\\
X \arrow[ur, "f"] \arrow[rr, swap, "g \circ f"] \arrow[loop left, "id_X"] & & Z \arrow[loop right, "id_Z"]
\end{tikzcd},
\end{equation}
then it is easy to see that $X$ here is an initial object, and, dually, that $Z$ is a final object, as illustrated by the elementary \textsc{Categorica} implementation shown in Figure \ref{fig:Figure20}. However, a slightly more subtle refinement of the concept of an ordinary initial object is that of a \textit{strict} initial object, namely an initial object $X$ for which all incoming morphisms ${f : Q \to X}$ (where $Q$ is an arbitrary object in the category ${\mathcal{C}}$) must be isomorphisms:

\begin{multline}
\forall X \in \mathrm{ob} \left( \mathcal{C} \right), \qquad X \text{ is a strict initial object }\\
\iff \qquad \forall P \in \mathrm{ob} \left( \mathcal{C} \right), \qquad \exists! \left( f : X \to P \right) \in \mathrm{hom} \left( \mathcal{C} \right), \qquad \text{ and }\\
\forall Q \in \mathrm{ob} \left( \mathcal{C} \right), \qquad \forall \left( g : Q \to X \right) \in \mathrm{hom} \left( \mathcal{C} \right), \qquad \exists \left( g^{-1} : X \to Q \right) \in \mathrm{hom} \left( \mathcal{C} \right), \qquad \text{ such that},\\
\left( g \circ g^{-1} : X \to X \right) = \left( id_X : X \to X \right), \qquad \text{ and } \qquad \left( g^{-1} \circ g : Q \to Q \right) = \left( id_Q : Q \to Q \right),
\end{multline}
or, illustrated diagrammatically:

\begin{equation}
\begin{tikzcd}
\forall Q \arrow[rr, bend left, "\forall g"] \arrow[loop left, "g^{-1} \circ g = id_Q"] & & X \arrow[ll, bend left, dashed, "\exists g^{-1}"] \arrow[rr, dashed, "\exists! f"] \arrow[loop above, "g \circ g^{-1} = id_X"] & & \forall P
\end{tikzcd}.
\end{equation}
Dually, a \textit{strict} terminal object is a terminal object $X$ for which all outgoing morphisms ${f : X \to Q}$ (where $Q$ is again an arbitrary object in the category ${\mathcal{C}}$) must be isomorphisms:

\begin{multline}
\forall X \in \mathrm{ob} \left( \mathcal{C} \right), \qquad X \text{ is a strict terminal object },\\
\iff \qquad \forall P \in \mathrm{ob} \left( \mathcal{C} \right), \qquad \exists! \left( f : P \to X \right) \in \mathrm{hom} \left( \mathcal{C} \right), \qquad \text{ and }\\
\forall Q \in \mathrm{ob} \left( \mathcal{C} \right), \qquad \forall \left( g : X \to Q \right) \in \mathrm{hom} \left( \mathcal{C} \right), \qquad \exists \left( g^{-1} : Q \to X \right) \in \mathrm{hom} \left( \mathcal{C} \right), \qquad \text{ such that }\\
\left( g \circ g^{-1} : Q \to Q \right) = \left( id_Q : Q \to Q \right), \qquad \text{ and } \qquad \left( g^{-1} \circ g : X \to X \right) = \left( id_X : X \to X \right),
\end{multline}
or, illustrated diagrammatically:

\begin{equation}
\begin{tikzcd}
\forall P \arrow[rr, dashed, "\exists! f"] & & X \arrow[rr, bend left, "\forall g"] \arrow[loop above, "g^{-1} \circ g = id_X"] & & \forall Q \arrow[ll, bend left, dashed, "\exists g^{-1}"] \arrow[loop right, "g \circ g^{-1} = id_Q"]
\end{tikzcd}.
\end{equation}

\begin{figure}[ht]
\centering
\begin{framed}
\includegraphics[width=0.495\textwidth]{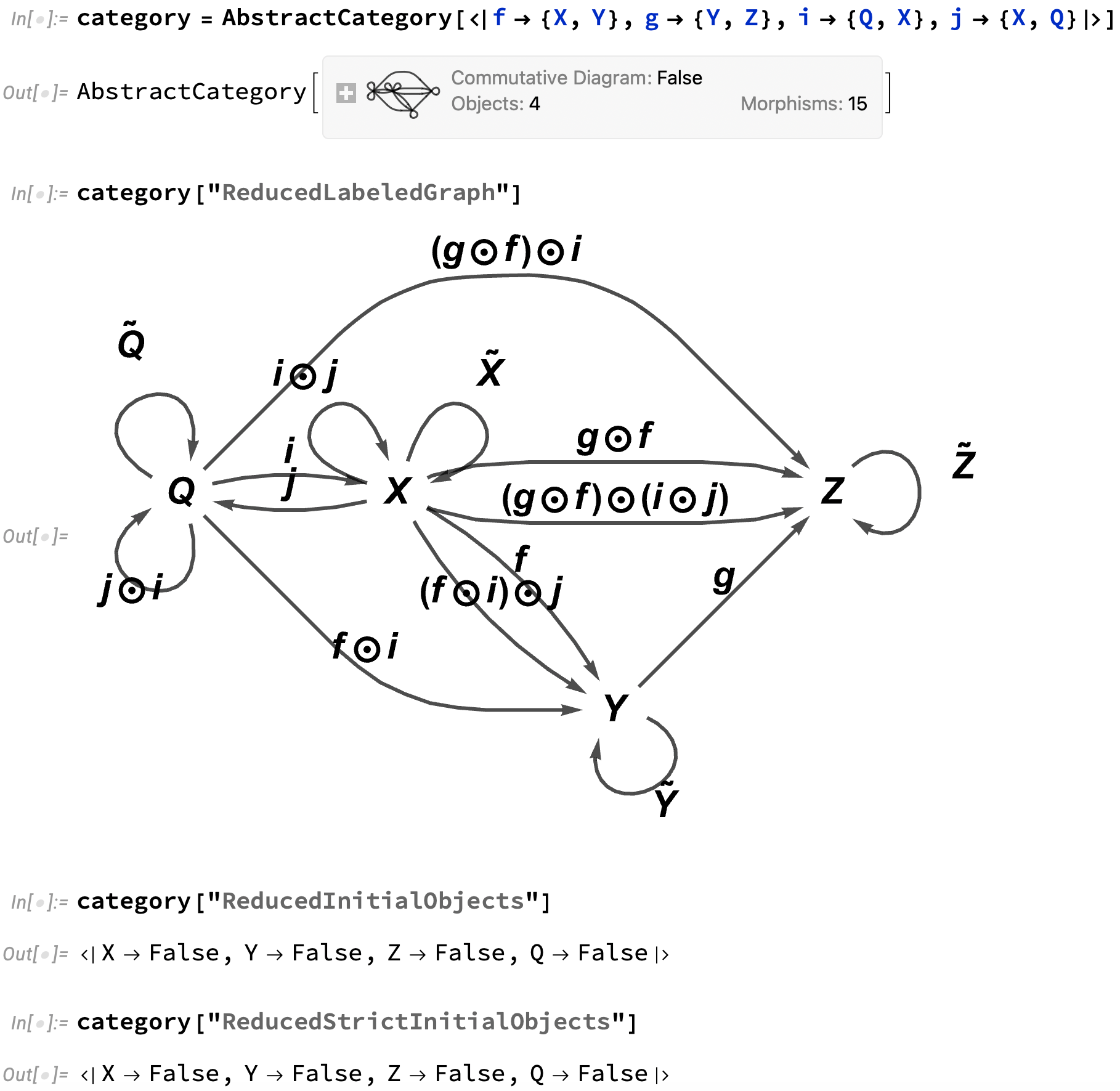}
\vrule
\includegraphics[width=0.495\textwidth]{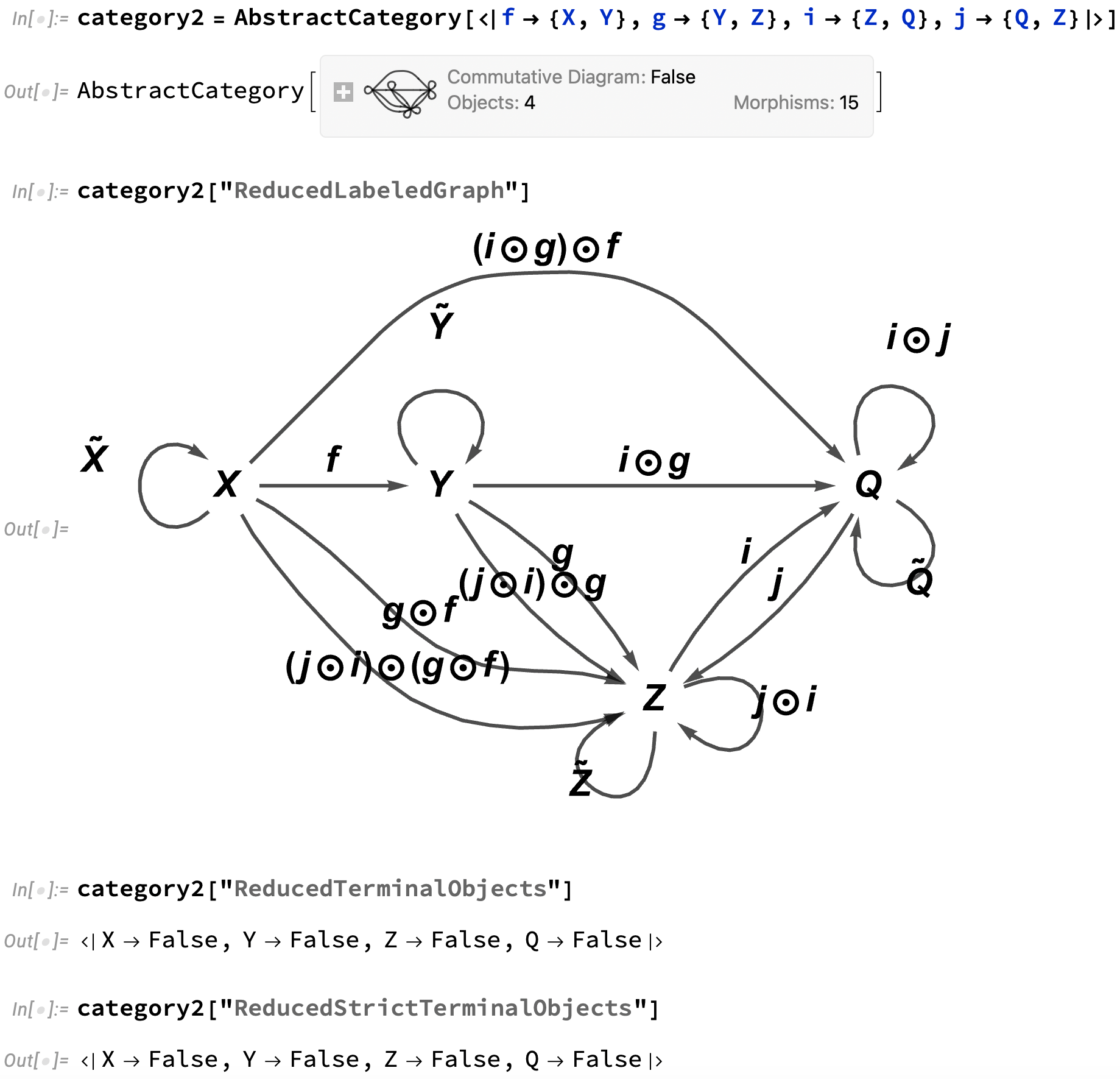}
\end{framed}
\caption{On the left, the \texttt{AbstractCategory} object corresponding to the simple triangular diagram from before, with new incoming and outgoing morphisms $i$ and $j$ added between the objects $X$ and $Q$, showing that $X$ has now ceased to be a (strict) initial object. On the right, the \texttt{AbstractCategory} object corresponding to the simple triangular diagram from before, with new incoming and outgoing morphisms $i$ and $j$ added between the objects $Z$ and $Q$, showing that $Z$ has now ceased to be a (strict) terminal object.}
\label{fig:Figure21}
\end{figure}

\begin{figure}[ht]
\centering
\begin{framed}
\includegraphics[width=0.495\textwidth]{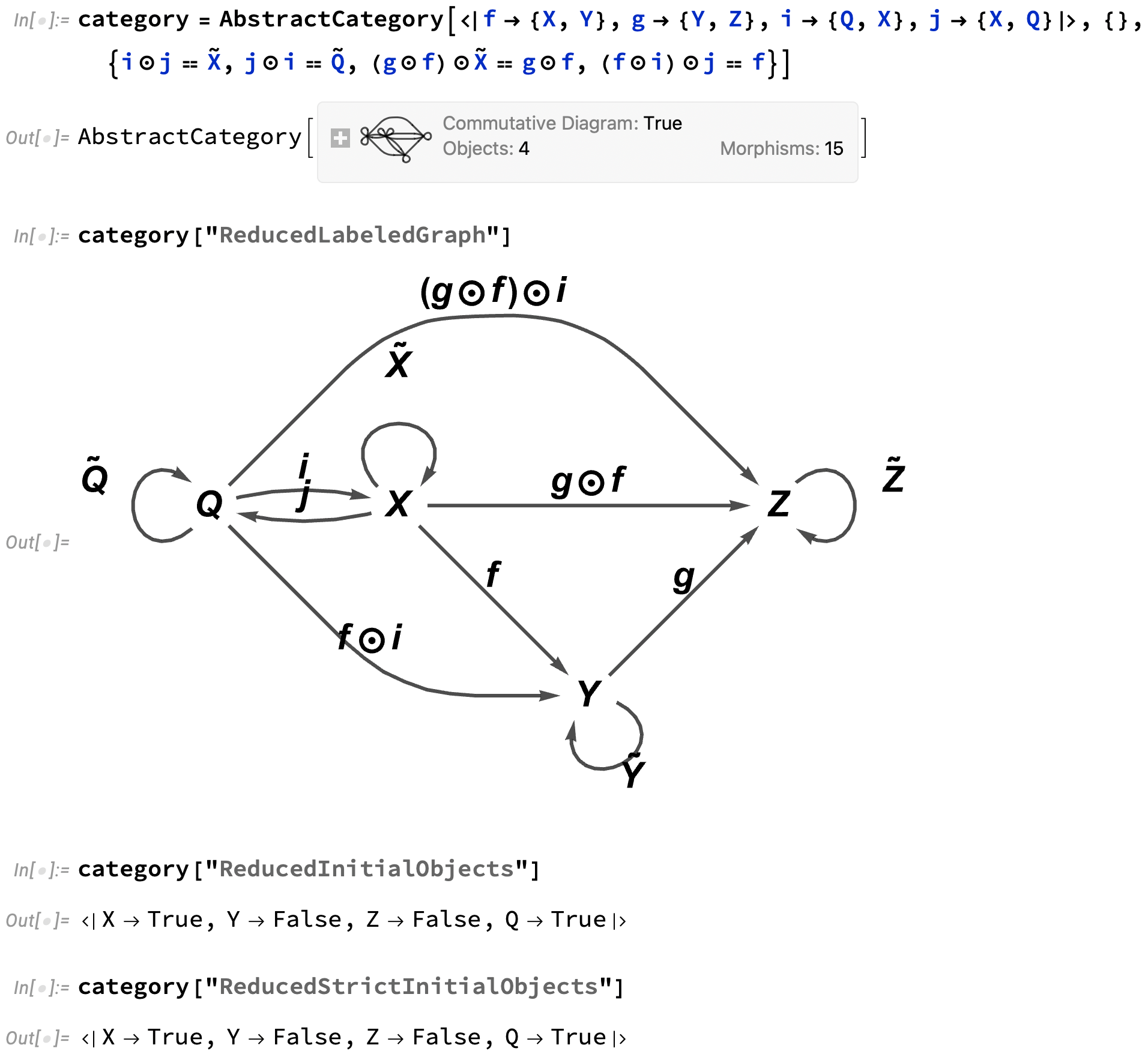}
\vrule
\includegraphics[width=0.495\textwidth]{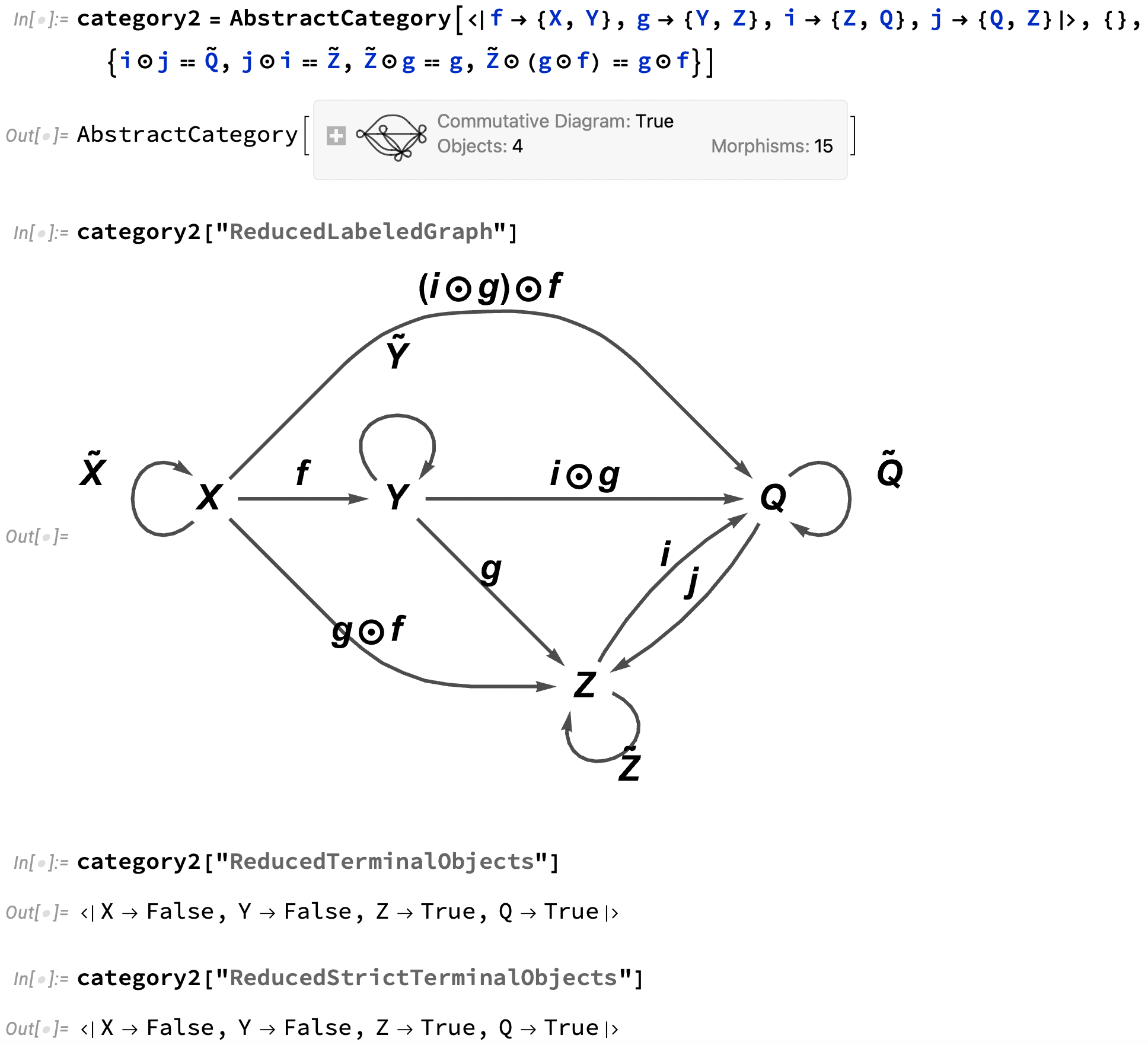}
\end{framed}
\caption{On the left, the \texttt{AbstractCategory} object corresponding to the simple triangular diagram from before, with new incoming and outgoing morphisms $i$ and $j$ added between the objects $X$ and $Q$, along with the morphism equivalences ${i \circ j = id_X}$ and ${j \circ i = id_Q}$ (as well as their corollaries), showing that $X$ is now a strict initial object. On the right, the \texttt{AbstractCategory} object corresponding to the simple triangular diagram from before, with incoming and outgoing morphisms $i$ and $j$ added between the objects $Z$ and $Q$, along with the morphism equivalences ${i \circ j = id_Q}$ and ${j \circ i = id_Z}$ (as well as their corollaries), showing that $Z$ is now a strict terminal object.}
\label{fig:Figure22}
\end{figure}

Returning again to the simple triangular diagram above, and introducing either a new incoming morphism ${h : Q \to X}$ to the initial object $X$, or a new outgoing morphism ${h : Z \to Q}$ from the terminal object $Z$ (along with corresponding morphisms ${h^{-1} : X \to Q}$ or ${h^{-1}: Q \to Z}$ in the reverse directions, respectively), i.e:

\begin{equation}
\begin{tikzcd}
& & & & & & Y \arrow[ddrr, "g"] \arrow[loop above, "id_Y"] & &\\ \\
Q \arrow[rrrr, swap, "h"] \arrow[uurrrrrr, bend left, "f \circ h"] \arrow[rrrrrrrr, bend right, swap, "\left( g \circ f \right) \circ h"] \arrow[loop above, "id_Q"] \arrow[loop left, "h^{-1} \circ h"] & & & & X \arrow[uurr, swap, "f"] \arrow[uurr, bend left, "\left( f \circ h \right) \circ h^{-1}"] \arrow[rrrr, "g \circ f"] \arrow[rrrr, bend right, "\left( g \circ f \right) \circ \left( h \circ h^{-1} \right)"] \arrow[llll, bend right, swap, "h^{-1}"] \arrow[loop above, "id_X"] \arrow[loop below, "h \circ h^{-1}"] & & & & Z \arrow[loop right, "id_Z"]
\end{tikzcd},
\end{equation}
or:

\begin{equation}
\begin{tikzcd}
& & Y \arrow[ddrr, swap, "g"] \arrow[ddrr, bend left, "\left( h^{-1} \circ h \right) \circ g"] \arrow[ddrrrrrr, bend left, "h \circ g"] \arrow[loop above, "id_Y"] & & & & & &\\ \\
X \arrow[uurr, "f"] \arrow[rrrr, "g \circ f"] \arrow[rrrr, bend right, "\left( h^{-1} \circ h \right) \circ \left( g \circ f \right)"] \arrow[rrrrrrrr, bend right, swap, "\left( h \circ g \right) \circ f"] \arrow[loop left, "id_X"] & & & & Z \arrow[rrrr, swap, "h"] \arrow[loop above, "id_Z"] \arrow[loop below, "h^{-1} \circ h"] & & & & Q \arrow[llll, bend right, swap, "h^{-1}"] \arrow[loop above, "id_Q"] \arrow[loop right, "h \circ h^{-1}"]
\end{tikzcd},
\end{equation}
respectively, we see that, in the first case, $X$ ceases to be a strict initial object (since the incoming morphism ${h : Q \to X}$ is not an isomorphism), and in the second case, $Z$ ceases to be a strict terminal object (since the outgoing morphism ${h : Z \to Q}$ is also not an isomorphism), as demonstrated in Figure \ref{fig:Figure21}. Indeed, $X$ and $Z$ are no longer initial/terminal objects even in the non-strict sense, due to the existence of multiple (currently algebraically inequivalent) outgoing morphisms:

\begin{multline}
\left( f : X \to Y \right) \neq \left( \left( f \circ h \right) \circ h^{-1} : X \to Y \right),\\
\text{ and } \qquad \left( g \circ f : X \to Z \right) \neq \left( \left( g \circ f \right) \circ \left( h \circ h^{-1} \right) : X \to Z \right),
\end{multline}
from $X$ to $Y$ and from $X$ to $Z$ in the former case, and multiple (currently algebraically inequivalent) incoming morphisms:

\begin{multline}
\left( g : Y \to Z \right) \neq \left( \left( h^{-1} \circ h \right) \circ g : Y \to Z \right)\\
\text{ and } \qquad \left( g \circ f : X \to Z \right) \neq \left( \left( h^{-1} \circ h \right) \circ \left( g \circ f \right) : X \to Z \right),
\end{multline}
from $Y$ to $Z$ and from $X$ to $Z$ in the latter case. However, if we now impose the pair of algebraic equivalences:

\begin{equation}
\left( h \circ h^{-1} : X \to X \right) = \left( id_X : X \to X \right), \qquad \text{ and } \qquad \left( h^{-1} \circ h : Q \to Q \right) = \left( id_Q : Q \to Q \right),
\end{equation}
thereby ensuring that:

\begin{multline}
\left( \left( f \circ h \right) \circ h^{-1} : X \to Y \right) = \left( f : X \to Y \right),\\
\text{ and } \qquad \left( \left( g \circ f \right) \circ \left( h \circ h^{-1} \right) : X \to Z \right) = \left( g \circ f : X \to Z \right),
\end{multline}
in the former case, and the pair of algebraic equivalences:

\begin{equation}
\left( h^{-1} \circ h : Z \to Z \right) = \left( id_Z : Z \to Z \right), \qquad \text{ and } \qquad \left( h \circ h^{-1} : Q \to Q \right) = \left( id_Q : Q \to Q \right),
\end{equation}
thereby ensuring that:

\begin{multline}
\left( \left( h^{-1} \circ h \right) \circ g : Y \to Z \right) = \left( g : Y \to Z \right),\\
\text{ and } \qquad \left( \left( h^{-1} \circ h \right) \circ \left( g \circ f \right) : X \to Z \right) = \left( g \circ f : X \to Z \right),
\end{multline}
in the latter case, then morphisms ${h : Q \to X}$ and ${h : Z \to Q}$ (as well as their corresponding reverse morphisms ${h^{-1} : X \to Q}$ and ${h^{-1} : Q \to Z}$, respectively) now become isomorphisms, collapsing the above diagrams to:

\begin{equation}
\begin{tikzcd}
& & & & & & Y \arrow[ddrr, "g"] \arrow[loop above, "id_Y"] & &\\ \\
Q \arrow[rrrr, swap, "h"] \arrow[uurrrrrr, bend left, "f \circ h"] \arrow[rrrrrrrr, bend right, swap, "\left( g \circ f \right) \circ h"] \arrow[loop left, "h^{-1} \circ h = id_Q"] & & & & X \arrow[uurr, "\left( f \circ h \right) \circ h^{-1} = f"] \arrow[rrrr, swap, "\left( g \circ f \right) \circ \left( h \circ h^{-1} \right) = g \circ f"] \arrow[llll, bend right, swap, "h^{-1}"] \arrow[loop below, "h \circ h^{-1} = id_X"] & & & & Z \arrow[loop right, "id_Z"]
\end{tikzcd},
\end{equation}
and:

\begin{equation}
\begin{tikzcd}
& & Y \arrow[ddrr, "\left( h^{-1} \circ h \right) \circ g = g"] \arrow[ddrrrrrr, bend left, "h \circ g"] \arrow[loop above, "id_Y"] & & & & & &\\ \\
X \arrow[uurr, "f"] \arrow[rrrr, swap, "\left( h^{-1} \circ h \right) \circ \left( g \circ f \right) = g \circ f"] \arrow[rrrrrrrr, bend right, swap, "\left( h \circ g \right) \circ f"] \arrow[loop left, "id_X"] & & & & Z \arrow[rrrr, swap, "h"] \arrow[loop below, "h^{-1} \circ h = id_Z"] & & & & Q \arrow[llll, bend right, "h^{-1}"] \arrow[loop right, "h \circ h^{-1} = id_Q"]
\end{tikzcd},
\end{equation}
respectively, and therefore causing $X$ and $Z$ to be strict initial and terminal objects, respectively, as illustrated in Figure \ref{fig:Figure22}. Any object that is both an initial object and a terminal object simultaneously is known as a \textit{zero object} (and, by extension, any object that is both a strict initial object and a strict terminal object simultaneously is known as a \textit{strict zero object}), and shares certain structural features in common with the point space object ${*}$ in the category ${\mathbf{Top_{*}}}$ of \textit{pointed} topological spaces (i.e. topological spaces with a distinguished base point). Categories that contain zero objects are consequently referred to as \textit{pointed categories}. \textsc{Categorica} also has specialized functionality for handling zero objects, strict zero objects and pointed categories built into the \texttt{AbstractCategory} function.

\begin{figure}[ht]
\centering
\begin{framed}
\includegraphics[width=0.495\textwidth]{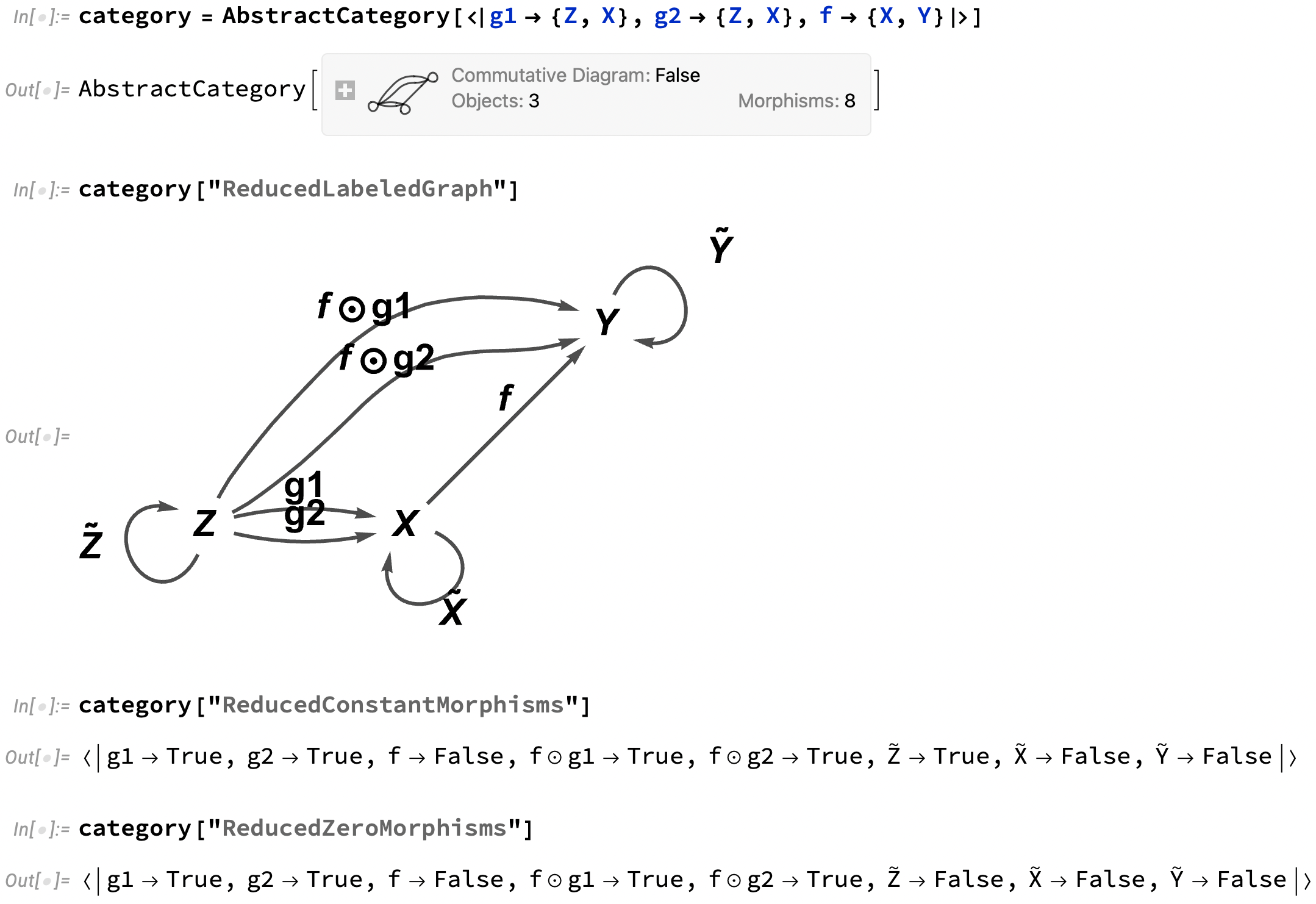}
\vrule
\includegraphics[width=0.495\textwidth]{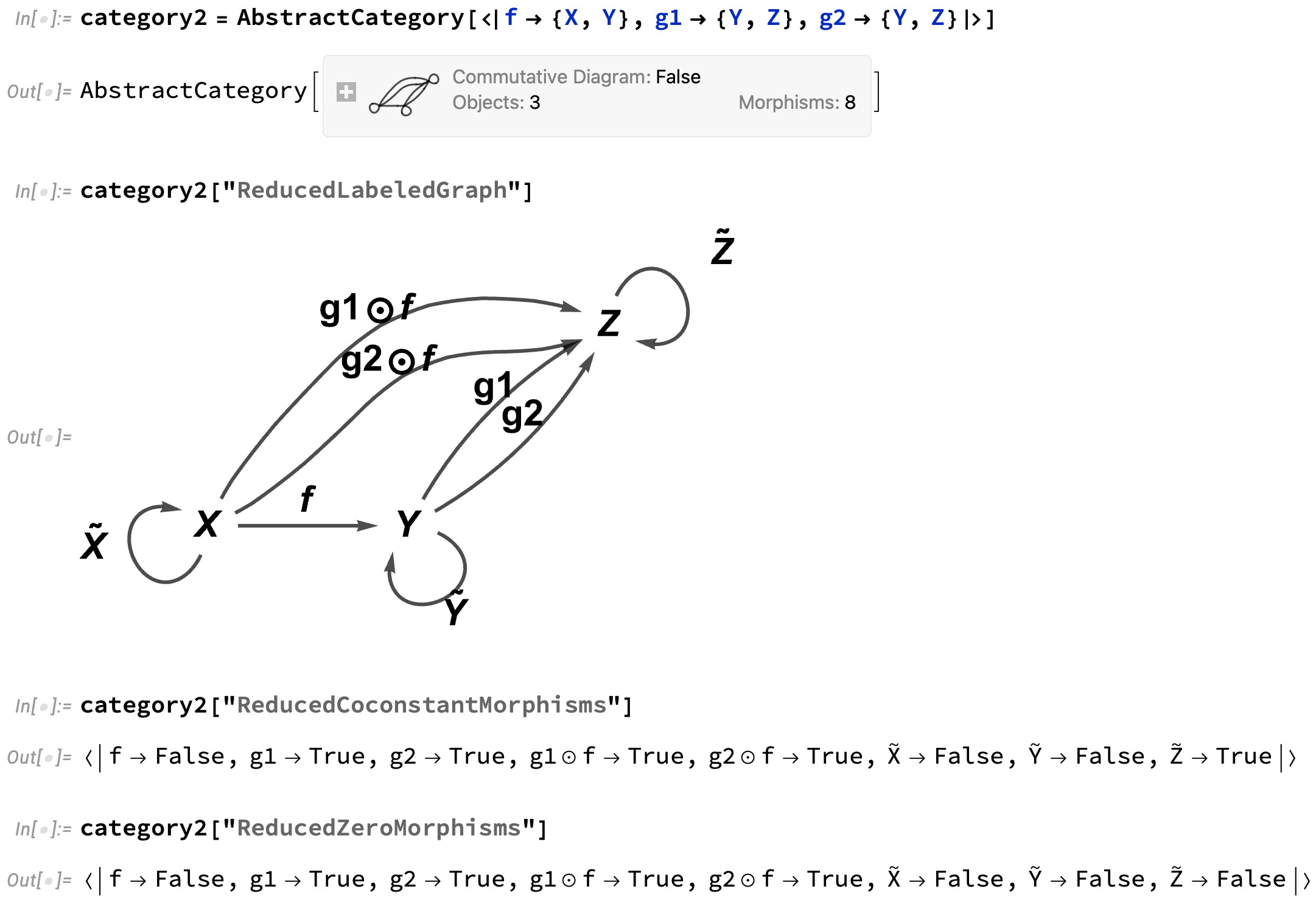}
\end{framed}
\caption{On the left, the \texttt{AbstractCategory} object corresponding to the basic diagrammatic setup of a constant morphism, showing that initially morphism $f$ is not a constant morphism, and hence also not a zero morphism. On the right, the \texttt{AbstractCategory} object corresponding to the basic diagrammatic setup of a coconstant morphism, showing that initially morphism $f$ is not a coconstant morphism, and hence also not a zero morphism.}
\label{fig:Figure23}
\end{figure}

\begin{figure}[ht]
\centering
\begin{framed}
\includegraphics[width=0.495\textwidth]{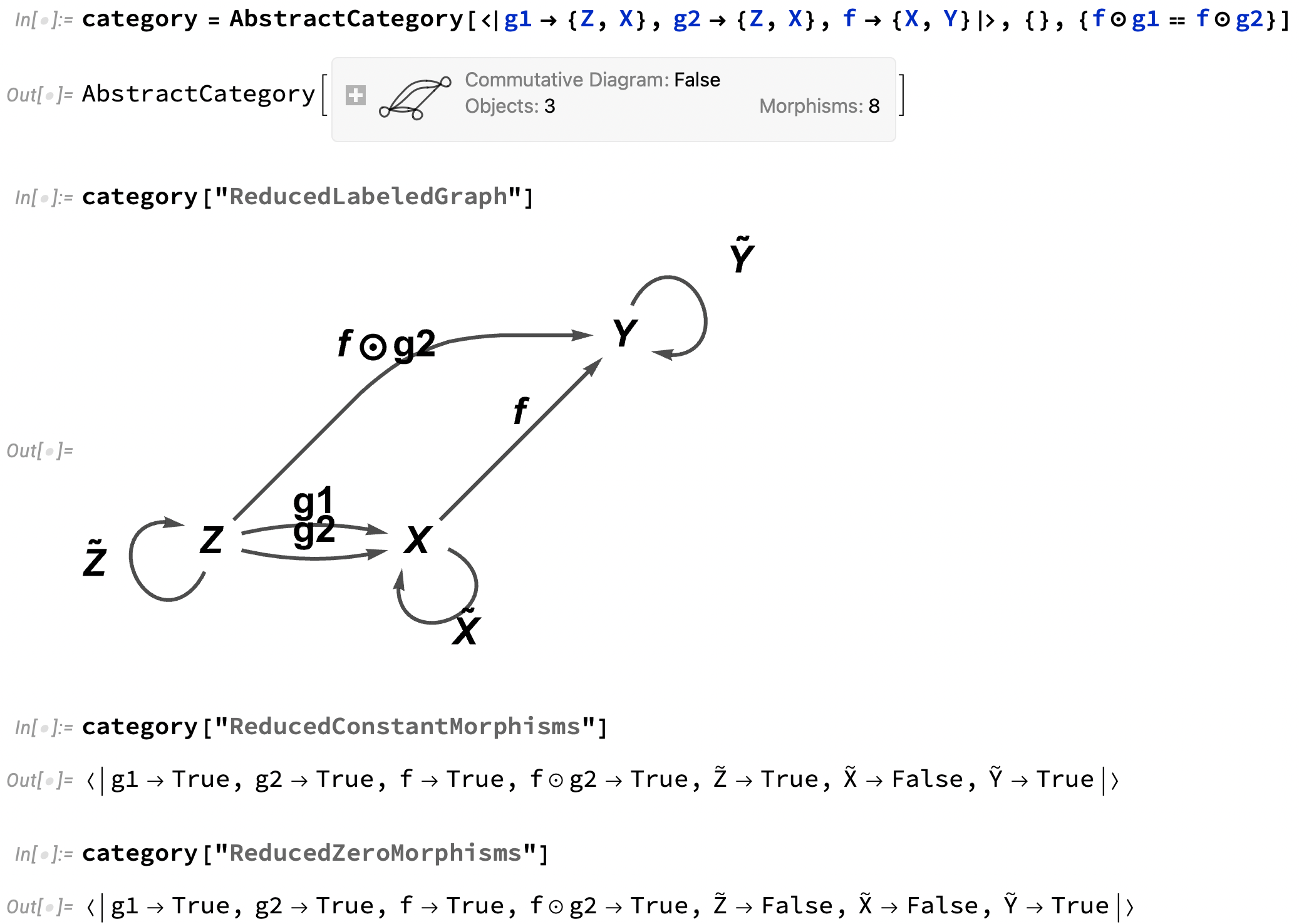}
\vrule
\includegraphics[width=0.495\textwidth]{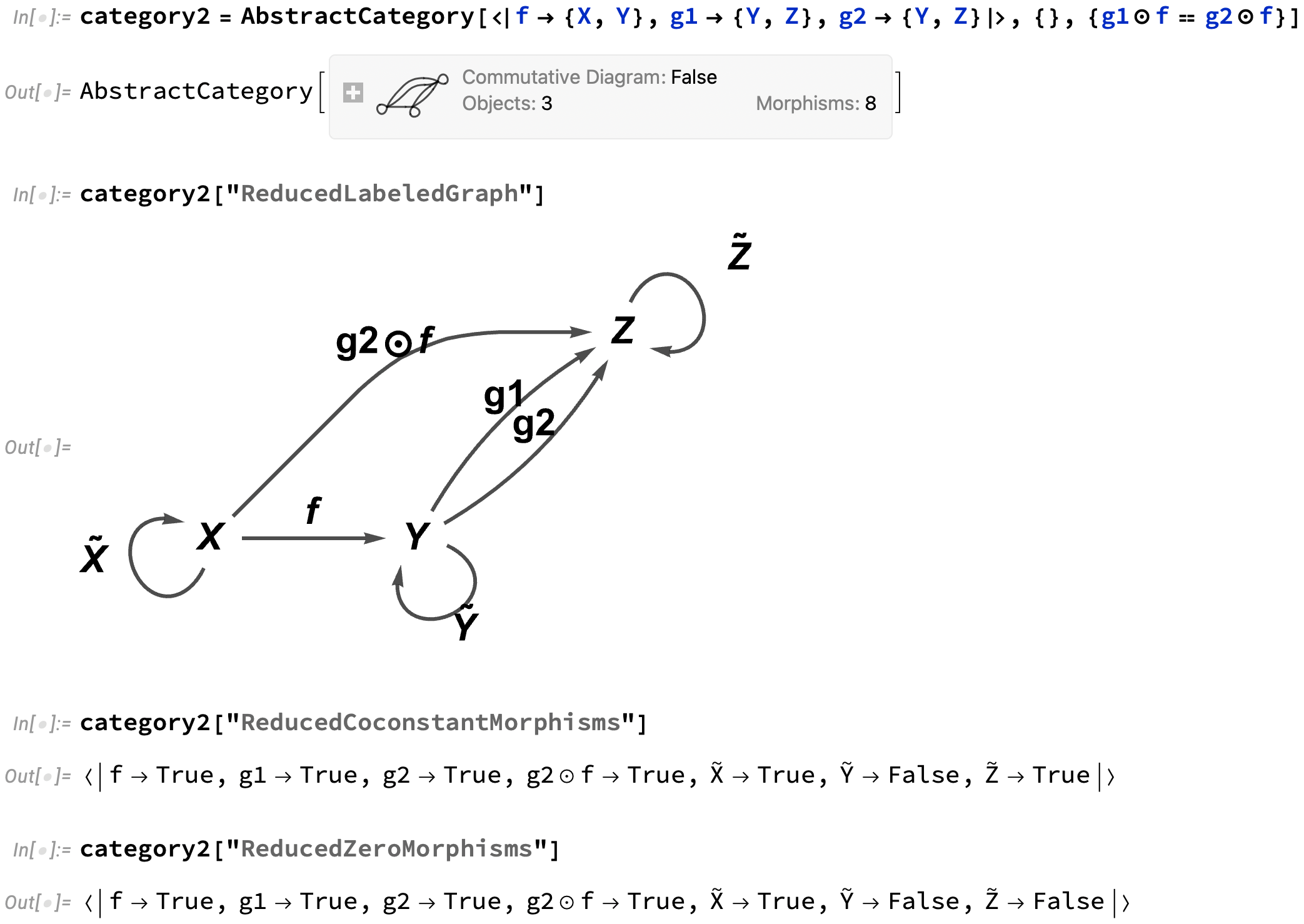}
\end{framed}
\caption{On the left, the \texttt{AbstractCategory} object corresponding to the basic diagrammatic setup of a constant morphism, with the additional algebraic equivalence ${f \circ g_1 = f \circ g_2}$ imposed, showing that morphism $f$ has now become a constant morphism, and hence also a zero morphism. On the right, the \texttt{AbstractCategory} object corresponding to the basic diagrammatic setup of a coconstant morphism, with the additional algebraic equivalence ${g_1 \circ f = g_2 \circ f}$ imposed, showing that morphism $f$ has now become a coconstant morphism, and hence also a zero morphism.}
\label{fig:Figure24}
\end{figure}

The final example of a dual construction that we shall discuss within the context of the present section, due to its close relationship with the monomorphism vs. epimorphism duality discussed previously, is that of \textit{constant} vs. \textit{coconstant} morphisms. If the morphism ${f : X \to Y}$ in the category ${\mathcal{C}}$ is such that, for all objects $Z$ and all pairs of morphisms ${g_1 : Z \to X}$ and ${g_2 : Z \to X}$, one necessarily has that ${\left( f \circ g_1 : Z \to Y \right) = \left( f \circ g_2 : Z \to Y \right)}$, then ${f : X \to Y}$ is a constant morphism:

\begin{multline}
\forall \left( f : X \to Y \right) \in \mathrm{hom} \left( \mathcal{C} \right), \qquad \left( f : X \to Y \right) \text{ is a constant morphism},\\
\iff \qquad \forall Z \in \mathrm{ob} \left( \mathcal{C} \right), \qquad \forall \left( g_1 : Z \to X \right), \left( g_2 : Z \to X \right) \in \mathrm{hom} \left( \mathcal{C} \right),\\
\text{ one has } \qquad \left( f \circ g_1 : Z \to Y \right) = \left( f \circ g_2 : Z \to Y \right),
\end{multline}
or, illustrated diagrammatically, one has:

\begin{equation}
\begin{tikzcd}
& X \arrow[dr, "f"] &\\
Z \arrow[ur, swap, "g_1"] \arrow[ur, bend left, "g_2"] \arrow[rr, swap, "f \circ g_1"] \arrow[rr, bend right, swap, "f \circ g_2"] & & Y
\end{tikzcd} \qquad \mapsto \qquad
\begin{tikzcd}
& X \arrow[dr, "f"] &\\
Z \arrow[ur, swap, "g_1"] \arrow[ur, bend left, "g_2"] \arrow[rr, swap, "f \circ g_1 = f \circ g_2"] & & Y
\end{tikzcd}.
\end{equation}
Note that this is identical to the diagram that appears in the \textit{hypothesis} of the definition of a monomorphism (and therefore we may say, slightly loosely, that a monomorphism is a constant morphism for which one necessarily has that ${\left( g_1 : Z \to X \right) = \left( g_2 : Z \to X \right)}$). Dually, if the morphism ${f : X \to Y}$ in the category ${\mathcal{C}}$ is such that, for all objects $Z$ and all pairs of morphisms ${g_1 : Y \to Z}$ and ${g_2 : Y \to Z}$, one necessarily has that ${\left( g_1 \circ f : X \to Z \right) = \left( g_2 \circ f : X \to Z \right)}$, then ${f : X \to Y}$ is a coconstant morphism:

\begin{multline}
\forall \left( f : X \to Y \right) \in \mathrm{hom} \left( \mathcal{C} \right), \qquad \left( f : X \to Y \right) \text{ is a coconstant morphism},\\
\iff \qquad \forall Z \in \mathrm{ob} \left( \mathcal{C} \right), \qquad \forall \left( g_1 : Y \to Z \right), \left( g_2 : Y \to Z \right) \in \mathrm{hom} \left( \mathcal{C} \right),\\
\text{ one has } \qquad \left( g_1 \circ f : X \to Z \right) = \left( g_2 \circ f : X \to Z \right),
\end{multline}
or, illustrated diagrammatically, one has:

\begin{equation}
\begin{tikzcd}
& Y \arrow[dr, swap, "g_1"] \arrow[dr, bend left, "g_2"] &\\
X \arrow[ur, "f"] \arrow[rr, swap, "g_1 \circ f"] \arrow[rr, bend right, swap, "g_2 \circ f"] & & Z
\end{tikzcd} \qquad \mapsto \qquad
\begin{tikzcd}
& Y \arrow[dr, swap, "g_1"] \arrow[dr, bend left, "g_2"] &\\
X \arrow[ur, "f"] \arrow[rr, swap, "g_1 \circ f = g_2 \circ f"] & & Z
\end{tikzcd}.
\end{equation}
Note that this is, likewise, identical to the diagram that appears in the \textit{hypothesis} of the definition of an epimorphism (and therefore we may also say, again slightly loosely, that an epimorphism is a coconstant morphism for which one necessarily has that ${\left( g_1 : Y \to Z \right) = \left( g_2 : Y \to Z \right)}$). Constant morphisms are so-named because they generalize the notion of constant functions in mathematical analysis (and, accordingly, coconstant morphisms may be thought of as generalizing the notion of zero-maps in analysis); any morphism that is both a constant morphism and a coconstant morphism is known as a 
\textit{zero} morphism, since such morphisms share certain key algebraic properties with morphisms mapping into and out of zero objects (by the same token, constant morphisms share certain algebraic properties with morphisms mapping into terminal objects, and coconstant morphisms share certain algebraic properties with morphisms mapping out of initial objects)\cite{pareigis}. Figure \ref{fig:Figure23} illustrates the basic diagrammatic setup for both constant and coconstant morphisms in \textsc{Categorica}, and demonstrates that initially (in the absence of any further algebraic equivalences) the morphism ${f : X \to Y}$ is not a constant morphism (in the former case), or not a coconstant morphism (in the latter case), and hence also not a zero morphism (in either case). Figure \ref{fig:Figure24} shows that, by imposing the algebraic equivalence ${\left( f \circ g_1 : Z \to Y \right) = \left( f \circ g_2 : Z \to Y \right)}$ in the former case and ${\left( g_1 \circ f : X \to Z \right) = \left( g_2 \circ f : X \to Z \right)}$ in the latter case, one is able to force the morphism ${f : X \to Y}$ to be a constant morphism in the former example, and to be a coconstant morphism in the latter example, and therefore a zero morphism in both examples.

The \texttt{AbstractCategory} function in \textsc{Categorica} also contains functionality for detecting and manipulating many other special types of morphism, including \textit{endomorphisms} (i.e. morphisms mapping an object $X$ to itself) through the \textit{``Endomorphisms''} property and its variants, \textit{automorphisms} (i.e. isomorphisms between an object $X$ and itself) through the \textit{``Automorphisms''} property and its variants, etc. There are also many other special types of category (beyond simply the groupoids and commutative diagrams discussed thus far) that can be automatically detected and represented using \textsc{Categorica}, with \texttt{AbstractCategory} in many cases including specialized functionality for manipulating them, including \textit{discrete categories} (i.e. categories in which the only morphisms are identity morphisms) through the \textit{``DiscreteCategoryQ''} property, \textit{indiscrete categories} (i.e. categories in which there exists a unique morphism between every pair of objects) through the \textit{``IndiscreteCategoryQ''} property, \textit{balanced categories} (i.e. categories in which every bimorphism is also necessarily an isomorphism) through the \textit{``BalancedCategoryQ''} property, etc. Needless to say, many additional properties and capabilities are also planned for future development. In many cases these properties have rather elegant mathematical interpretations; for instance, discrete categories may be thought as being a natural category-theoretic generalization of discrete topological spaces (i.e. topological spaces equipped with a \textit{discrete topology}, in which every subset is open), since in a discrete topological space the only allowable paths are the identity paths from a point to itself, etc. However, in the interests of brevity, we will not cover these additional capabilities in any detail here.

\clearpage

\section{Functors and Fibrations}
\label{sec:Section3}

\begin{figure}[ht]
\centering
\begin{framed}
\includegraphics[width=0.495\textwidth]{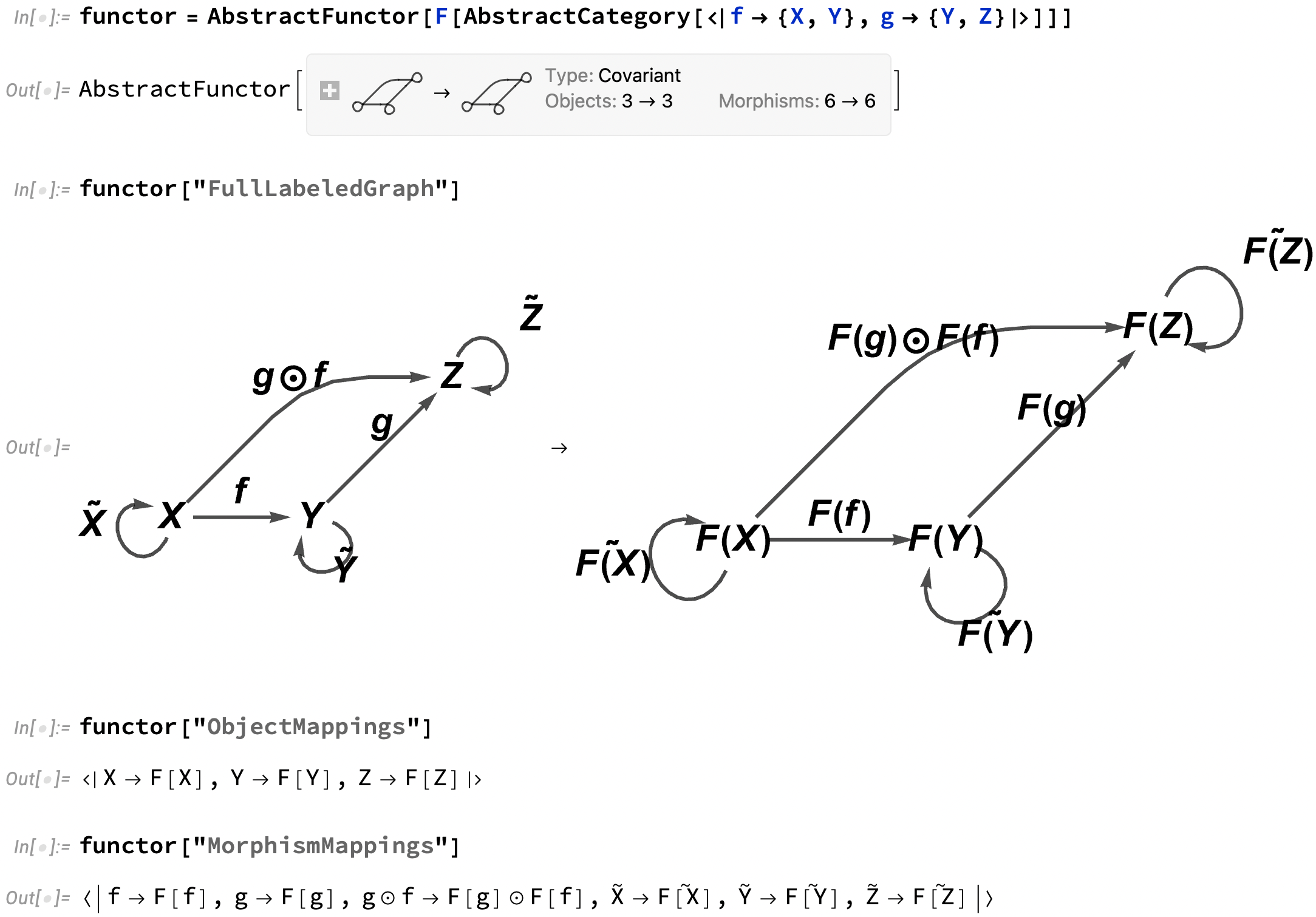}
\vrule
\includegraphics[width=0.495\textwidth]{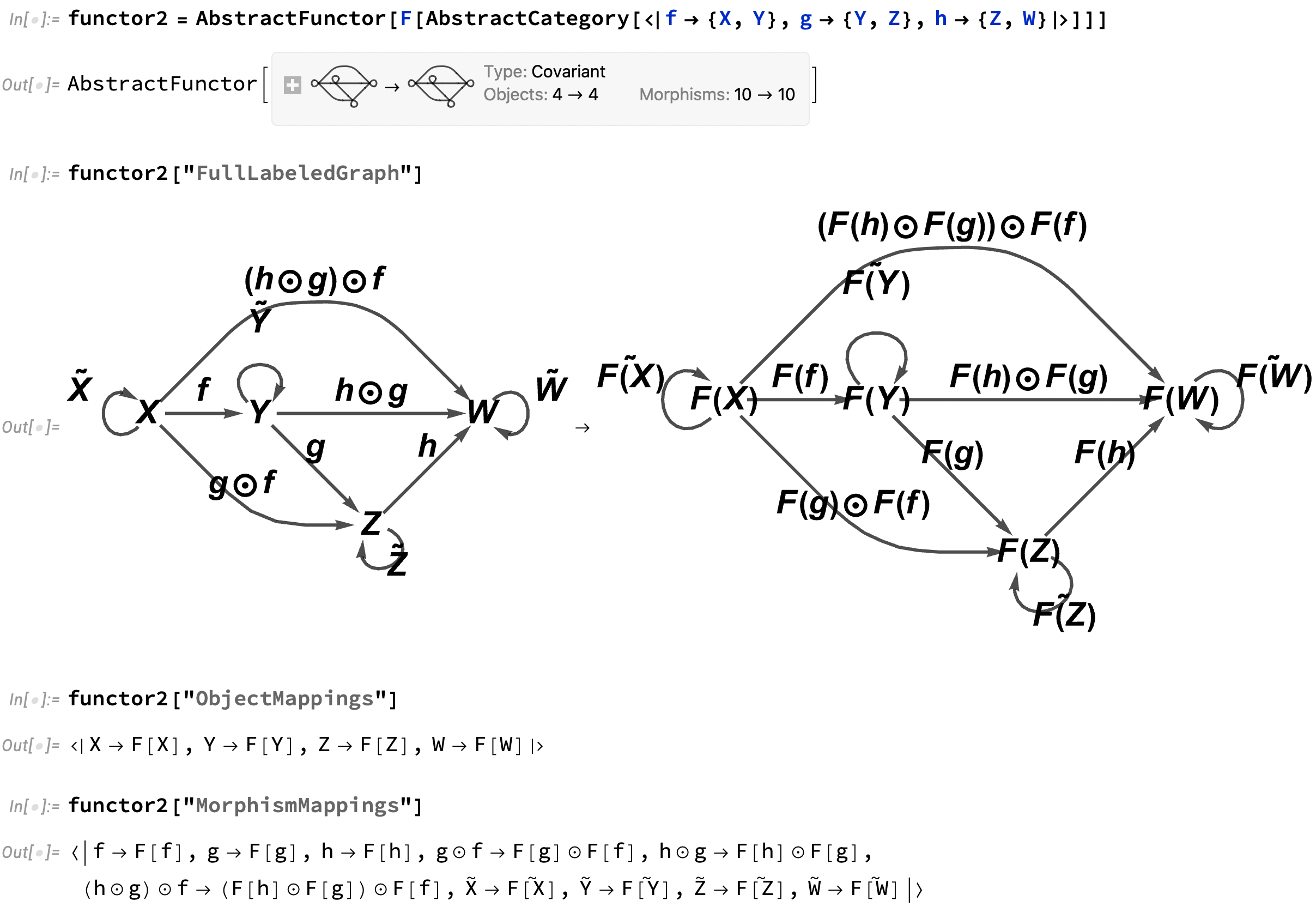}
\end{framed}
\caption{On the left, the \texttt{AbstractFunctor} object for a simple (three object, six-morphism) domain category, showing the explicit mappings between objects and morphisms in the domain and codomain categories. On the right, the \texttt{AbstractFunctor} object for a slightly larger (four-object, ten-morphism) domain category, showing the explicit mappings between objects and morphisms in the domain and codomain categories.}
\label{fig:Figure12}
\end{figure}

A homomorphism (i.e. a structure-preserving map) between two categories is known as a \textit{functor}\cite{maclane}\cite{jacobson}. More precisely, the map ${F : \mathcal{C} \to \mathcal{D}}$ from a category ${\mathcal{C}}$ to a category ${\mathcal{D}}$ is a functor if and only if it associates every object $X$ in category ${\mathcal{C}}$ to a corresponding object ${F \left( X \right)}$ in category ${\mathcal{D}}$, i.e:

\begin{equation}
\forall X \in \mathrm{ob} \left( \mathcal{C} \right), \qquad \exists F \left( X \right) \in \mathrm{ob} \left( \mathcal{D} \right),
\end{equation}
and every morphism ${f : X \to Y}$ in category ${\mathcal{C}}$ to a corresponding morphism ${F \left( f \right) : F \left( X \right) \to F \left( Y \right)}$ in category ${\mathcal{D}}$, i.e:

\begin{equation}
\forall \left( f : X \to Y \right) \in \mathrm{hom} \left( \mathcal{C} \right), \qquad \exists \left( F \left( f \right) : F \left( X \right) \to F \left( Y \right) \right) \in \mathrm{hom} \left( \mathcal{D} \right),
\end{equation}
in such a way that the identity morphisms are all preserved, such that the identity morphism ${id_X : X \to X}$ on object $X$ in category ${\mathcal{C}}$ is always mapped to the identity morphism ${id_{F \left( X \right)} : F \left( X \right) \to F \left( X \right)}$ on object ${F \left( X \right)}$ in category ${\mathcal{D}}$, i.e:

\begin{equation}
\forall X \in \mathrm{ob} \left( \mathcal{C} \right), \qquad \left( F \left( id_X \right) : F \left( X \right) \to F \left( X \right) \right) = \left( id_{F \left( X \right)} : F \left( X \right) \to F \left( X \right) \right),
\end{equation}
and the composition of morphisms is also preserved, such that the composite morphism ${g \circ f : X \to Z}$ in category ${\mathcal{C}}$ (where ${f : X \to Y}$ and ${g : Y \to Z}$ are also morphisms in ${\mathcal{C}}$) is always mapped to the composite morphism ${F \left( g \right) \circ F \left( f \right) : F \left( X \right) \to F \left( Z \right)}$ in category ${\mathcal{D}}$, i.e:

\begin{multline}
\forall \left( f : X \to Y \right), \left( g : Y \to Z \right) \in \mathrm{hom} \left( \mathcal{C} \right),\\
\left( F \left( g \circ f \right) : F \left( X \right) \to F \left( Z \right) \right) = \left( F \left( g \right) \circ F \left( f \right) : F \left( X \right) \to F \left( Z \right) \right).
\end{multline}
The functor construction may be illustrated diagrammatically by means of the following minimal example:

\begin{equation}
\begin{tikzcd}
& Y \arrow[dr, "g"] \arrow[loop above, "id_Y"] &\\
X \arrow[ur, "f"] \arrow[rr, swap, "g \circ f"] \arrow[loop below, "id_X"] & & Z \arrow[loop below, "id_Z"]
\end{tikzcd} \qquad \mapsto \qquad
\begin{tikzcd}
& F \left( Y \right) \arrow[dr, "F \left( g \right)"] \arrow[loop above, "id_{F \left( Y \right)}"] &\\
F \left( X \right) \arrow[ur, "F \left( f \right)"] \arrow[rr, swap, "F \left( g \right) \circ F \left( f \right)"] \arrow[loop below, "id_{F \left( X \right)}"] & & F \left( Z \right) \arrow[loop below, "id_{F \left( Z \right)}"]
\end{tikzcd}.
\end{equation}
\textit{Functoriality} turns out to be an extremely powerful algebraic condition, and functors may consequently be thought of as generalizing many central constructions in pure mathematics, including group actions in the context of group theory (wherein one considers functors from single-object groupoids, i.e. groups, to the category ${\mathbf{Set}}$ whose objects are sets and whose morphisms are set-valued functions), group \textit{representations} in the context of representation theory (wherein one again considers functors from groups/single-object groupoids, but now to the category ${\mathbf{Vect}}$ whose objects are vector spaces and whose morphisms are linear maps), and tangent bundles in differential geometry (wherein one considers functors from the category ${\textbf{Man}^p}$ whose objects are differentiable manifolds of smoothness class ${C^p}$ and whose morphisms are $p$-times continuously differentiable maps, to the category ${\mathbf{Vect} \left( \textbf{Man}^p \right)}$ whose objects are vector bundles over the manifolds in ${\mathbf{Man}^p}$ and whose morphisms are vector bundle homomorphisms). Indeed, as mentioned previously, even categorical diagrams themselves may be formalized as functors ${D : \mathcal{J} \to \mathcal{C}}$ from some index category/scheme ${\mathcal{J}}$ to an arbitrary category ${\mathcal{C}}$. Power sets in set theory (i.e. functors from the category ${\mathbf{Set}}$ of sets and set-valued functions to itself) and fundamental groups in algebraic topology (i.e. functors from the category ${\mathbf{Top_{*}}}$ of pointed topological spaces and continuous functions between them, to the category ${\mathbf{Grp}}$ of groups and group homomorphisms) constitute further examples of key functorial constructions in mathematics that have previously been alluded to within this article. Figure \ref{fig:Figure12} shows the implementation of the minimal example presented above, along with a slightly larger example, in the \textsc{Categorica} framework using the \texttt{AbstractFunctor} function; much like with the \texttt{AbstractCategory} objects themselves, \textsc{Categorica} automatically keeps track of, and enforces, all necessary algebraic equivalences between objects and morphisms within the codomain category of an \texttt{AbstractFunctor} object that must be imposed in order to maintain consistency with the functor axioms described above.

\begin{figure}[ht]
\centering
\begin{framed}
\includegraphics[width=0.495\textwidth]{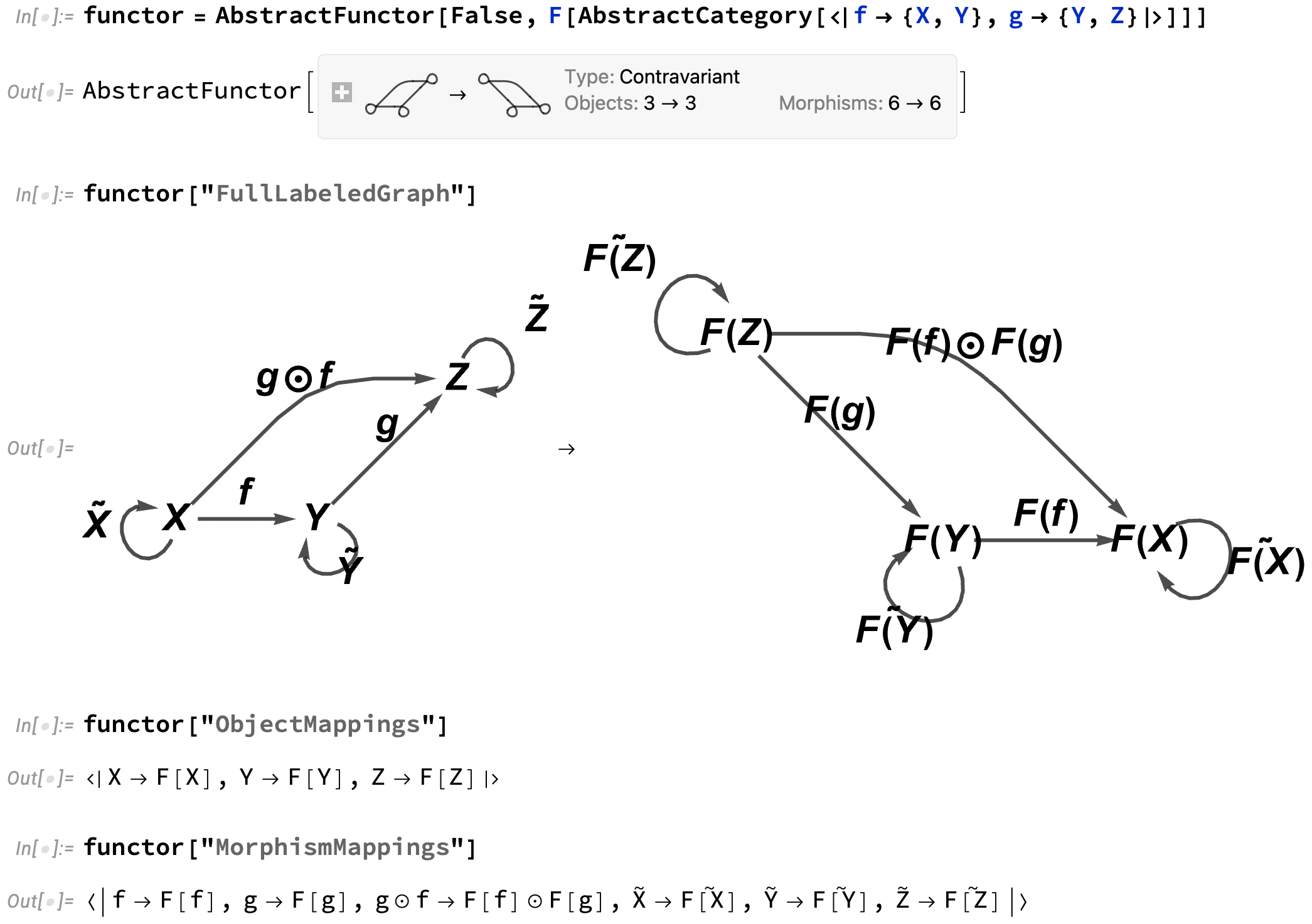}
\vrule
\includegraphics[width=0.495\textwidth]{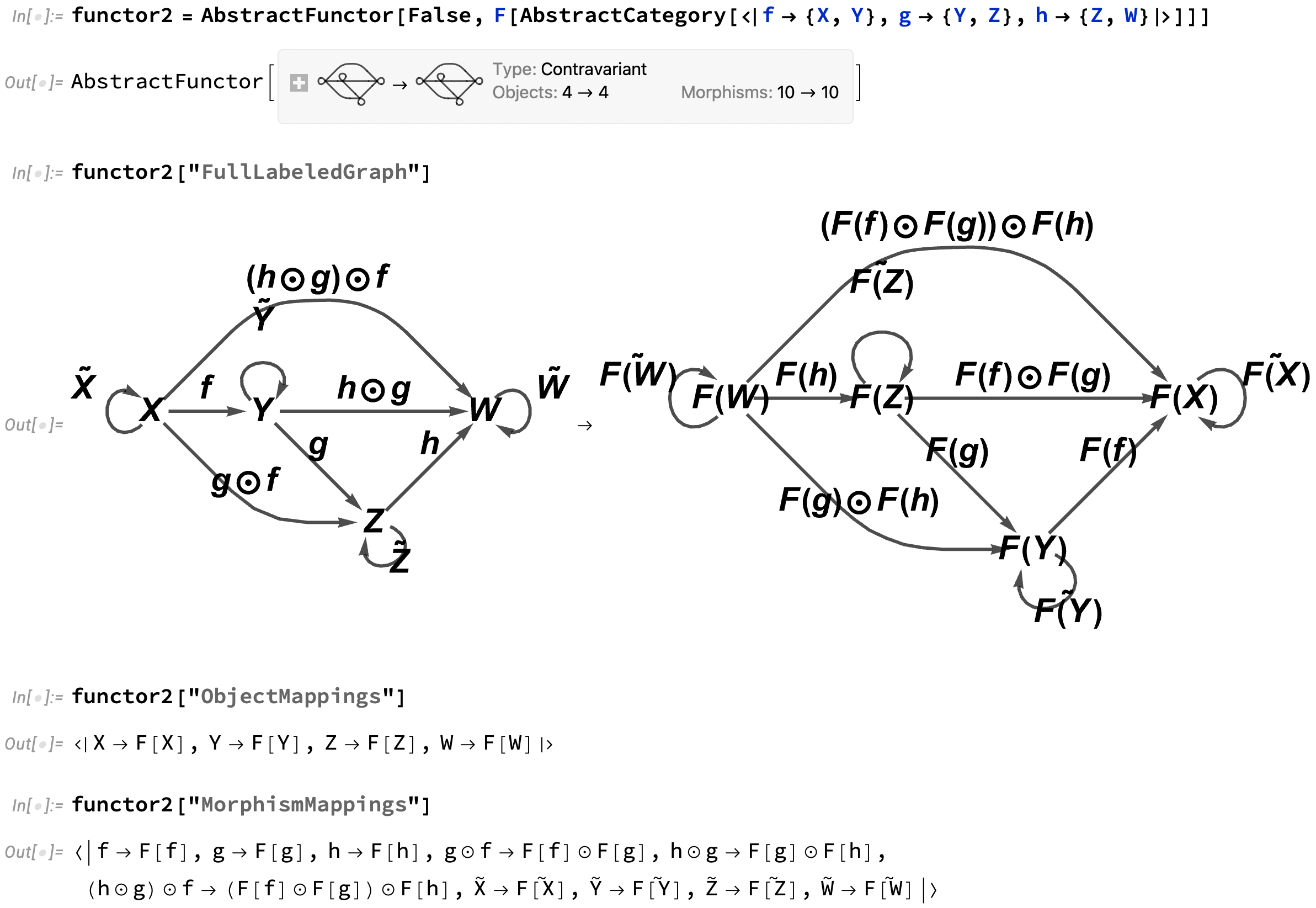}
\end{framed}
\caption{On the left, the contravariant \texttt{AbstractFunctor} object for a simple (three-object, six-morphism) domain category, showing the explicit mappings between objects and morphisms in the (dual) domain and codomain categories, illustrating reversal of the order of morphism composition. On the left, the contravariant \texttt{AbstractFunctor} object for a slightly larger (four-object, ten-morphism) domain category, showing the explicit mappings between objects and morphisms in the (dual) domain and codomain categories, illustrating reversal of the order of morphism composition.}
\label{fig:Figure13}
\end{figure}

Although the formal definition presented above corresponds to the case of \textit{covariant} functors (namely functors which preserve the directions of morphisms, and hence which also preserve the order of composition of morphisms), there also exists a dual notion of \textit{contravariant} functors, which obey analogous axioms but which reverse the directions of morphisms, and hence also reverse the order of morphism composition. The definition of the mapping on objects remains unchanged from the covariant case, but now every morphism ${f : X \to Y}$ in category ${\mathcal{C}}$ is mapped to a corresponding morphism ${F \left( f \right) : F \left( Y \right) \to F \left( X \right)}$ in category ${\mathcal{D}}$, i.e:

\begin{equation}
\forall \left( f : X \to Y \right) \in \mathrm{hom} \left( \mathcal{C} \right), \qquad \exists \left( F \left( f \right) : F \left( Y \right) \to F \left( X \right) \right) \in \mathrm{hom} \left( \mathcal{D} \right).
\end{equation}
The condition on identity morphisms is also unchanged from the covariant case, but now the composite morphism ${g \circ f : X \to Z}$ in category ${\mathcal{C}}$ (where ${f : X \to Y}$ and ${g : Y \to Z}$ are also morphisms in ${\mathcal{C}}$) is always mapped to the composite morphism ${F \left( f \right) \circ F \left( g \right) : F \left( Z \right) \to F \left( X \right)}$ in category ${\mathcal{D}}$, i.e:

\begin{multline}
\forall \left( f : X \to Y \right), \left( g : Y \to Z \right) \in \mathrm{hom} \left( \mathcal{C} \right),\\
\left( F \left( g \circ f \right) : F \left( Z \right) \to F \left( X \right) \right) = \left( F \left( f \right) \circ F \left( g \right) : F \left( Z \right) \to F \left( X \right) \right).
\end{multline}
The contravariant functor construction may, just as before, be illustrated diagrammatically by means of the following minimal example:

\begin{equation}
\begin{tikzcd}
& Y \arrow[dr, "g"] \arrow[loop above, "id_Y"] &\\
X \arrow[ur, "f"] \arrow[rr, swap, "g \circ f"] \arrow[loop below, "id_X"] & & Z \arrow[loop below, "id_Z"]
\end{tikzcd} \qquad \mapsto \qquad
\begin{tikzcd}
& F \left( Y \right) \arrow[dr, "F \left( f \right)"] \arrow[loop above, "id_{F \left( Y \right)}"] &\\
F \left( Z \right) \arrow[ur, "F \left( g \right)"] \arrow[rr, swap, "F \left( f \right) \circ F \left( g \right)"] \arrow[loop below, "id_{F \left( Z \right)}"] & & F \left( X \right) \arrow[loop below, "id_{F \left( X \right)}"]
\end{tikzcd}.
\end{equation}
Just as covariant functors from the category ${\mathbf{Man}^p}$ of differentiable manifolds and differentiable maps to the category ${\mathbf{Vect} \left( \mathbf{Man}^p \right)}$ of vector bundles and vector bundle homomorphisms may be used to generalize the construction of tangent bundles in differential geometry, the construction of \textit{cotangent} bundles may correspondingly be generalized using the corresponding \textit{contravariant} functor from ${\mathbf{Man}^p}$ (or, more precisely, from its dual/opposite category ${\left( \mathbf{Man}^p \right)^{op}}$) to ${\mathbf{Vect} \left( \mathbf{Man}^p \right)}$. Likewise, the operation of taking dual spaces in linear algebra corresponds to the application of a certain contravariant functor from the category ${\mathbf{Vect} \left( K \right)}$ of vector spaces and linear maps over some fixed field $K$ (or, more precisely, from its dual/opposite category ${\left( \mathbf{Vect} \left( K \right) \right)^{op}}$), to itself. Indeed, the power set construction in set theory may in fact be formalized as \textit{either} a covariant functor ${\mathcal{P} : \mathbf{Set} \to \mathbf{Set}}$ (as above), \textit{or} as a contravariant functor ${\mathcal{P} : \mathbf{Set}^{op} \to \mathbf{Set}}$. The notion of presheaves on a space $X$ in algebraic topology may be formalized as a functor from the dual/opposite category ${\mathcal{C}^{op}}$ of the category ${\mathcal{C}}$ whose objects are open sets ${U \subseteq X}$ and whose morphisms are inclusion maps ${U \subseteq V}$, to ${\mathbf{Set}}$; when appropriately extended to the case where the category ${\mathcal{C}}$ (and hence the domain category ${\mathcal{C}^{op}}$ of the contravariant functor) is arbitrary, this yields the corresponding category-theoretic notion of a presheaf, which is far more general. Figure \ref{fig:Figure13} demonstrates a \textsc{Categorica} implementation of the corresponding contravariant versions of the two functors previously shown in Figure \ref{fig:Figure12}, obtained in each case by setting the first argument of \texttt{AbstractFunctor} (i.e. the covariance argument) to \texttt{False}. At any time, the dual version of any \texttt{AbstractFunctor} object may be computed directly using the \textit{``SwapVariance''} property (analogous to the \textit{``DualCategory''} property of \texttt{AbstractCategory} objects), with \textsc{Categorica} automatically making any required modifications to the order of morphism composition within all morphism equivalence lists, etc.

\begin{figure}[ht]
\centering
\begin{framed}
\includegraphics[width=0.545\textwidth]{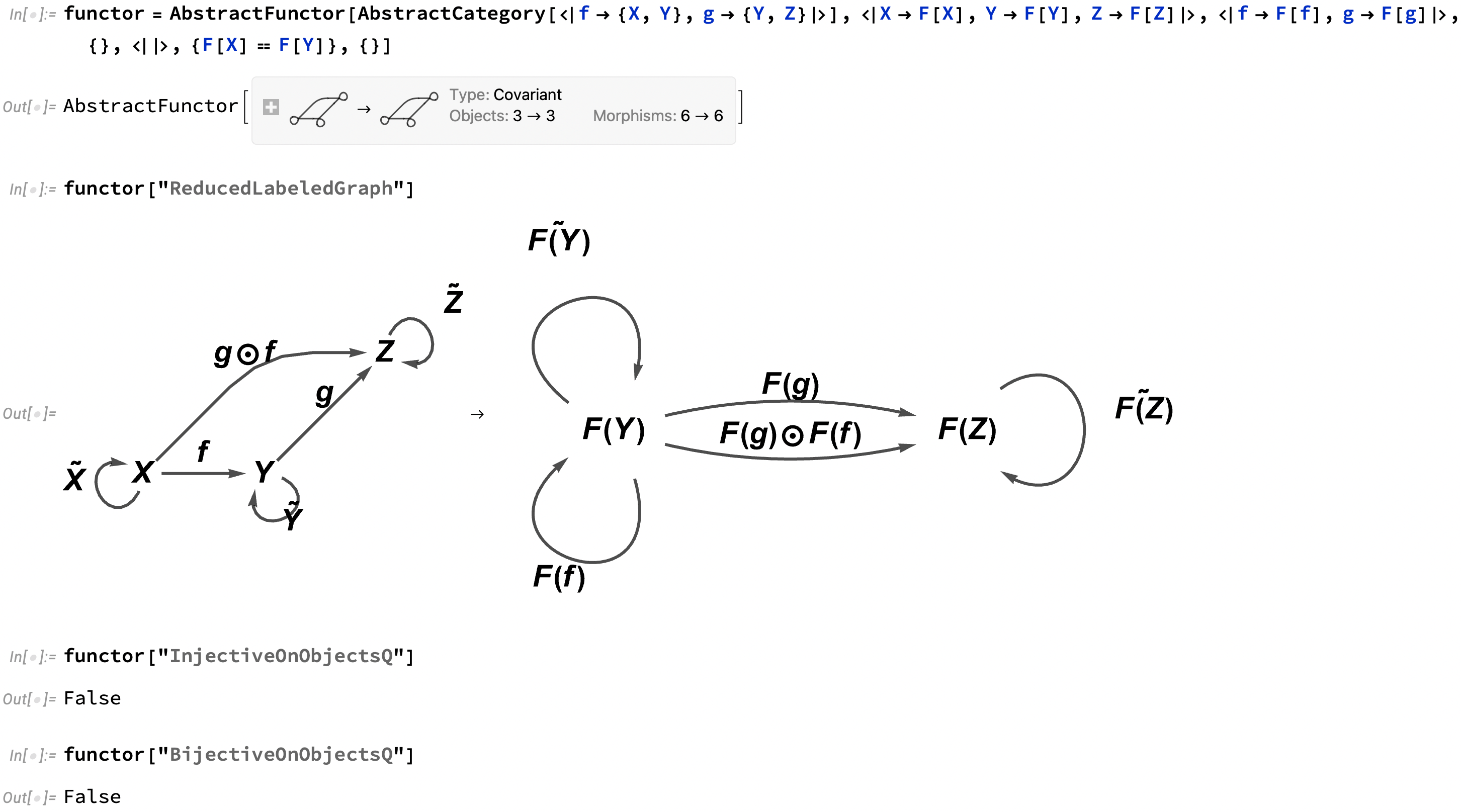}
\vrule
\includegraphics[width=0.445\textwidth]{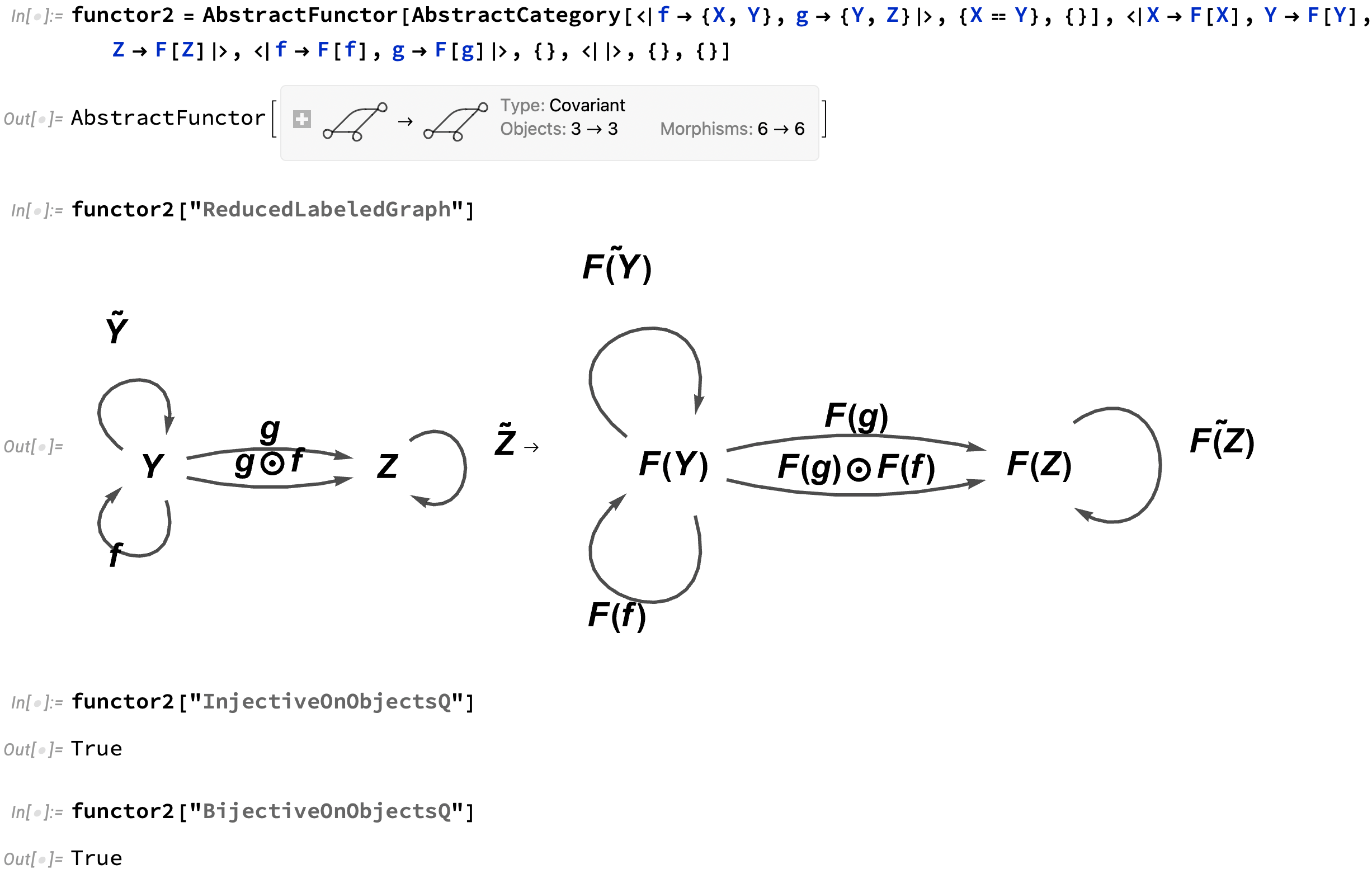}
\end{framed}
\caption{On the left, an \texttt{AbstractFunctor} object which is not injective, and hence not bijective, on objects because of the additional algebraic equivalence ${F \left( X \right) = F \left( Y \right)}$ being imposed on the objects of the codomain category. On the right, the corresponding \texttt{AbstractFunctor} object with injectivity, and hence bijectivity, on objects restored, by imposing the additional algebraic equivalence ${X = Y}$ on the objects of the domain category.}
\label{fig:Figure14}
\end{figure}

\begin{figure}[ht]
\centering
\begin{framed}
\includegraphics[width=0.445\textwidth]{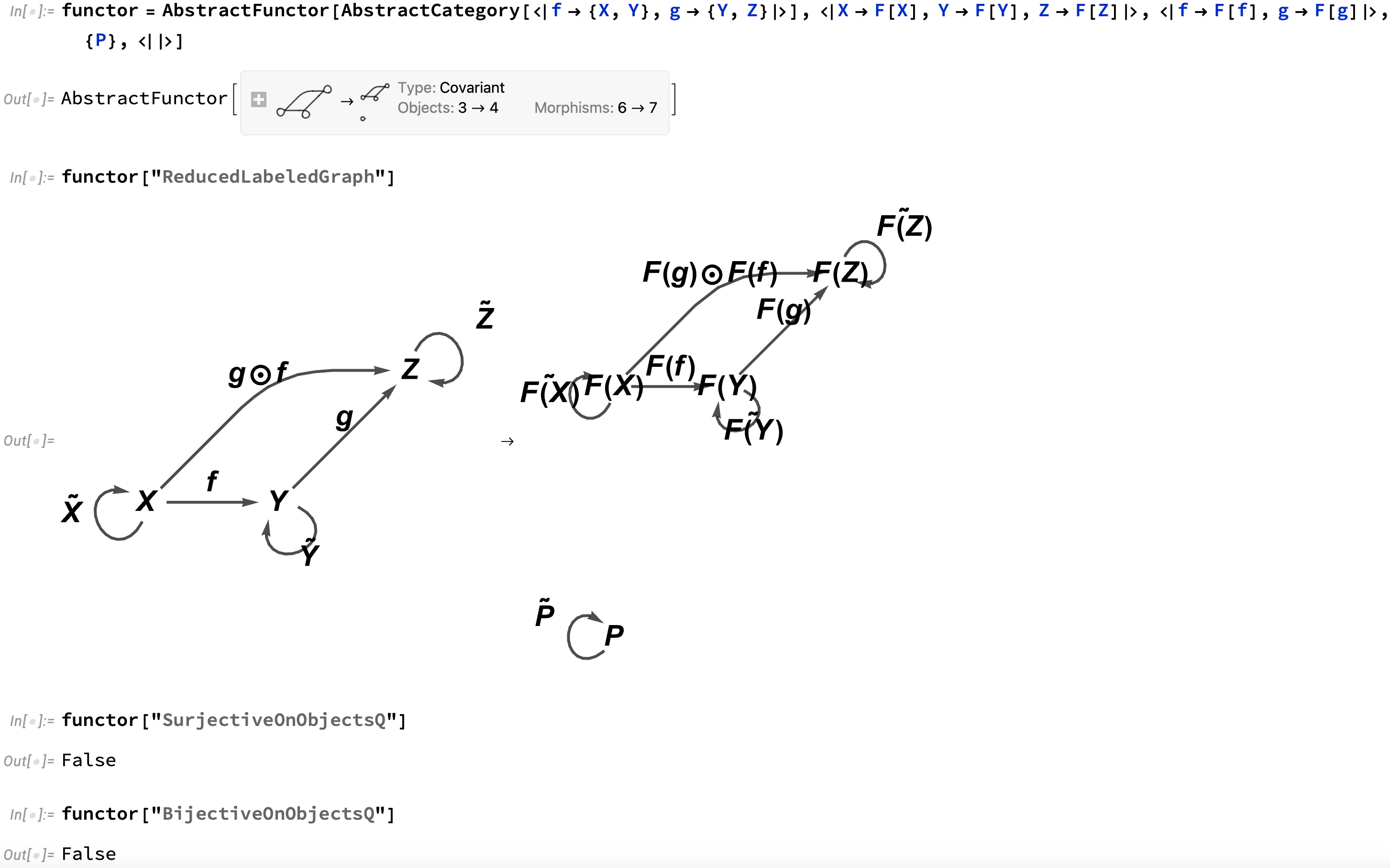}
\vrule
\includegraphics[width=0.545\textwidth]{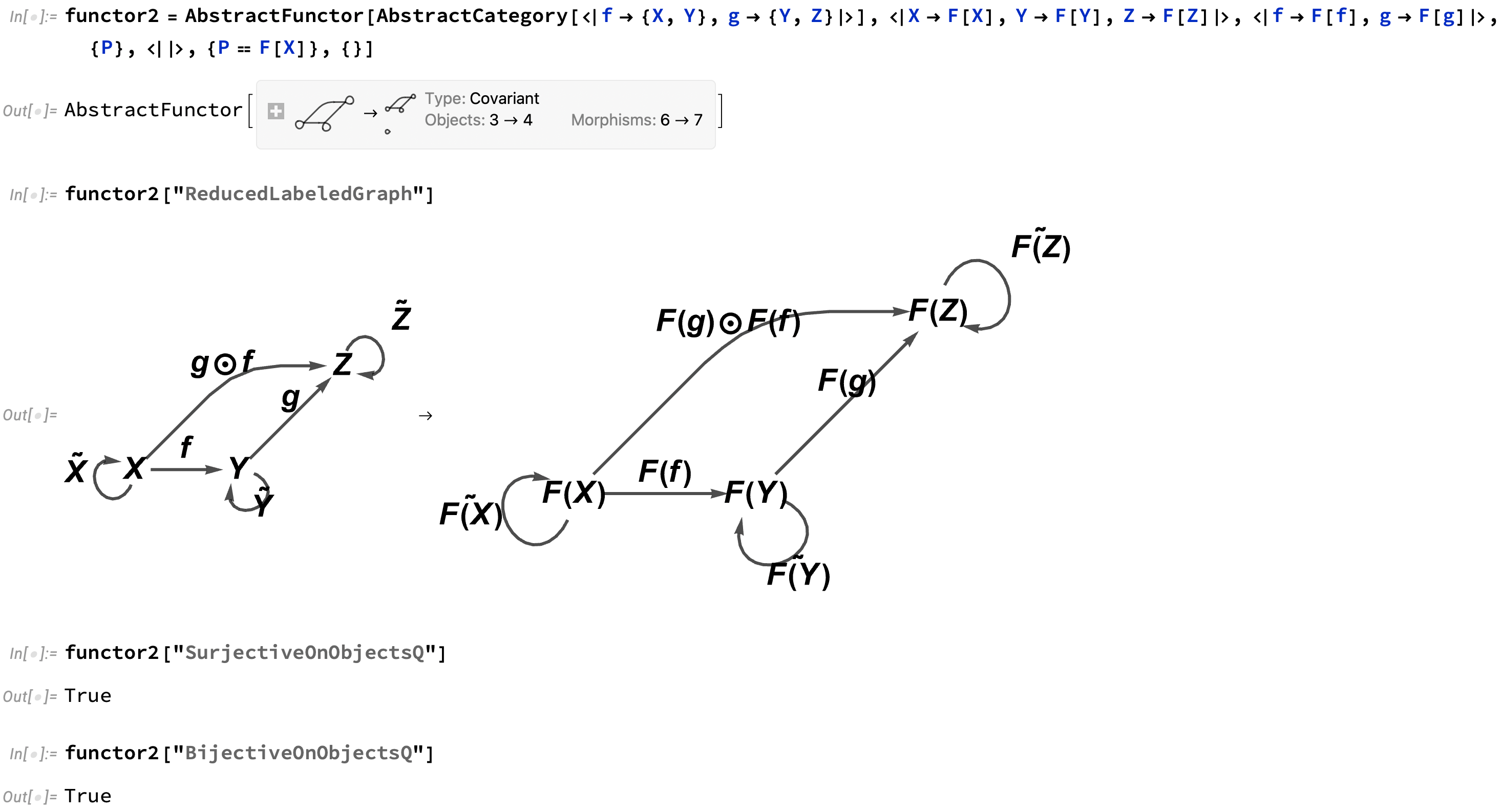}
\end{framed}
\caption{On the left, an \texttt{AbstractFunctor} object which is not surjective, and hence not bijective, on objects because of the additional object $P$ being introduced within the codomain category. On the right, the corresponding \texttt{AbstractFunctor} object with surjectivity, and hence bijectivity, on objects restored, by imposing the additional algebraic equivalence ${P = F \left( X \right)}$ on the objects of the codomain category.}
\label{fig:Figure15}
\end{figure}

In its full form, and in the absence of any additional variance information, an \texttt{AbstractFunctor} object is specified by the \texttt{AbstractCategory} object for its domain category, an association of object mappings, an association of morphism (or, more precisely, arrow) mappings, a list of new objects in the codomain category, a list of new morphisms (or, more precisely, arrows) in the codomain category, a list of new object equivalences in the codomain category, and a list of new morphism equivalences in the codomain category (as well as optional additional information regarding choices of composition and identity symbols, as with \texttt{AbstractCategory}). Note, in particular, that all mappings take place at the level of the underlying \texttt{AbstractQuiver} object, rather than on the \texttt{AbstractCategory} object itself, so as to guarantee consistency with the composition/functoriality axioms. All algebraic equivalences between objects and arrows/morphisms in the domain \texttt{AbstractCategory} object (or its underlying \texttt{AbstractQuiver} object) are automatically translated into corresponding algebraic equivalences in the codomain \texttt{AbstractCategory} object. We may therefore regard a functor ${F : \mathcal{C} \to \mathcal{D}}$ as corresponding to a pair of (set-valued, since we are dealing here with small categories) functions, ${F_{ob}}$ and ${F_{hom}}$, with the former mapping the set of objects in category ${\mathcal{C}}$ to the set of objects in category ${\mathcal{D}}$ and the latter mapping the set of morphisms in category ${\mathcal{C}}$ to the set of morphisms in category ${\mathcal{D}}$:

\begin{equation}
F_{ob} : \mathrm{ob} \left( \mathcal{C} \right) \to \mathrm{ob} \left( \mathcal{D} \right), \qquad \text{ and } \qquad F_{hom} : \mathrm{hom} \left( \mathcal{C} \right) \to \mathrm{hom} \left( \mathcal{D} \right).
\end{equation}
This, in turn, presents (at least) two distinct ways in which functors may be classified as either injective, surjective or bijective; namely, functors may be injective/surjective/bijective on objects, or on morphisms, or on both, or on neither (depending upon the injectivity/surjectivity/bijectivity properties of the functions ${F_{ob}}$ and ${F_{hom}}$). A simple way in which a functor may fail to be injective on objects is by having two objects in the codomain category be equivalent which were not previously equivalent in the domain category; for instance, consider the functor:

\begin{equation}
\begin{tikzcd}
& Y \arrow[dr, "g"] \arrow[loop above, "id_Y"] &\\
X \arrow[ur, "f"] \arrow[rr, swap, "g \circ f"] \arrow[loop below, "id_X"] & & Z \arrow[loop below, "id_Z"]
\end{tikzcd} \qquad \mapsto \qquad
\begin{tikzcd}
& F \left( Y \right) \arrow[dr, "F \left( g \right)"] \arrow[loop above, "id_{F \left( Y \right)}"] &\\
F \left( X \right) \arrow[ur, "F \left( f \right)"] \arrow[rr, swap, "F \left( g \right) \circ F \left( f \right)"] \arrow[loop below, "id_{F \left( X \right)}"] & & F \left( Z \right) \arrow[loop below, "id_{F \left( Z \right)}"]
\end{tikzcd},
\end{equation}
plus the additional algebraic condition that ${F \left( X \right) = F \left( Y \right)}$, imposed on the objects in the codomain category ${\mathcal{D}}$, such that one instead has:

\begin{equation}
\begin{tikzcd}
& Y \arrow[dr, "g"] \arrow[loop above, "id_Y"] &\\
X \arrow[ur, "f"] \arrow[rr, swap, "g \circ f"] \arrow[loop below, "id_X"] & & Z \arrow[loop below, "id_Z"]
\end{tikzcd} \qquad \mapsto \qquad
\begin{tikzcd}
F \left( X \right) = F \left( Y \right) \arrow[rr, bend left, "F \left( g \right)"] \arrow[rr, bend right, swap, "F \left( g \right) \circ F \left( f \right)"] \arrow[loop above, "id_{F \left( X \right)} = id_{F \left( Y \right)}"] \arrow[loop below, "F \left( f \right)"] & & F \left( Z \right) \arrow[loop right, "id_{F \left( Z \right)}"]
\end{tikzcd}.
\end{equation}
Such a functor would no longer be injective, and hence also no longer bijective, on objects, although its injectivity (and therefore bijectivity) on objects could nevertheless be restored by introducing a corresponding algebraic condition ${X = Y}$ on the objects in the domain category ${\mathcal{C}}$, such that one instead has:

\begin{equation}
\begin{tikzcd}
X = Y \arrow[rr, bend left, "g"] \arrow[rr, bend right, swap, "g \circ f"] \arrow[loop above, "id_X = id_Y"] \arrow[loop below, "f"] & & Z \arrow[loop right, "id_Z"]
\end{tikzcd} \qquad \mapsto \qquad
\begin{tikzcd}
F \left( X \right) = F \left( Y \right) \arrow[rr, bend left, "F \left( g \right)"] \arrow[rr, bend right, swap, "F \left( g \right) \circ F \left( f \right)"] \arrow[loop above, "id_{F \left( X \right)} = id_{F \left( Y \right)}"] \arrow[loop below, "F \left( f \right)"] & & F \left( Z \right) \arrow[loop right, "id_{F \left( Z \right)}"]
\end{tikzcd}.
\end{equation}
Figure \ref{fig:Figure14} implements this basic example in \textsc{Categorica}, demonstrating how functor injectivity on objects may be removed (by imposing the algebraic equivalence ${F \left( X \right) = F \left( Y \right)}$ on objects in the codomain category) and then subsequently reinstated (by imposing the additional algebraic equivalence ${X = Y}$ on objects in the domain category). Dual to this construction, a simple way in which a functor may fail to be surjective on objects is by having a new object be introduced in the codomain category that did not previously exist in the domain category; for instance, consider now the functor:

\begin{equation}
\begin{tikzcd}
& Y \arrow[dr, "g"] \arrow[loop above, "id_Y"] &\\
X \arrow[ur, "f"] \arrow[rr, swap, "g \circ f"] \arrow[loop below, "id_X"] & & Z \arrow[loop below, "id_Z"]
\end{tikzcd} \qquad \mapsto \qquad
\begin{tikzcd}
& F \left( Y \right) \arrow[dr, "F \left( g \right)"] \arrow[loop above, "id_{F \left( Y \right)}"] & & P \arrow[loop above, "id_P"]\\
F \left( X \right) \arrow[ur, "F \left( f \right)"] \arrow[rr, swap, "F \left( g \right) \circ F \left( f \right)"] \arrow[loop below, "id_{F \left( X \right)}"] & & F \left( Z \right) \arrow[loop below, "id_{F \left( Z \right)}"] &
\end{tikzcd}.
\end{equation}
Such a functor would no longer be surjective, and hence also no longer bijective, on objects, although its surjectivity (and therefore bijectivity) on objects could nevertheless be restored by introducing a new algebraic condition ${P = F \left( X \right)}$ on the objects in the codomain category ${\mathcal{D}}$, such that one instead has:

\begin{equation}
\begin{tikzcd}
& Y \arrow[dr, "g"] \arrow[loop above, "id_Y"] &\\
X \arrow[ur, "f"] \arrow[rr, swap, "g \circ f"] \arrow[loop below, "id_X"] & & Z \arrow[loop below, "id_Z"]
\end{tikzcd} \qquad \mapsto \qquad
\begin{tikzcd}
& F \left( Y \right) \arrow[dr, "F \left( g \right)"] \arrow[loop above, "id_{F \left( Y \right)}"] &\\
F \left( X \right) = P \arrow[ur, "F \left( f \right)"] \arrow[rr, swap, "F \left( g \right) \circ F \left( f \right)"] \arrow[loop below, "id_{F \left( X \right)} = id_P"] & & F \left( Z \right) \arrow[loop below, "id_{F \left( Z \right)}"]
\end{tikzcd}.
\end{equation}
Figure \ref{fig:Figure15} implements this dual example in \textsc{Categorica}, demonstrating now how functor surjectivity on objects may be removed (by introducing the new object $P$ in the codomain category) and then subsequently reinstated (by imposing the algebraic equivalence ${P = F \left( X \right)}$ on objects in the codomain category), in much the same way.

\begin{figure}[ht]
\centering
\begin{framed}
\includegraphics[width=0.495\textwidth]{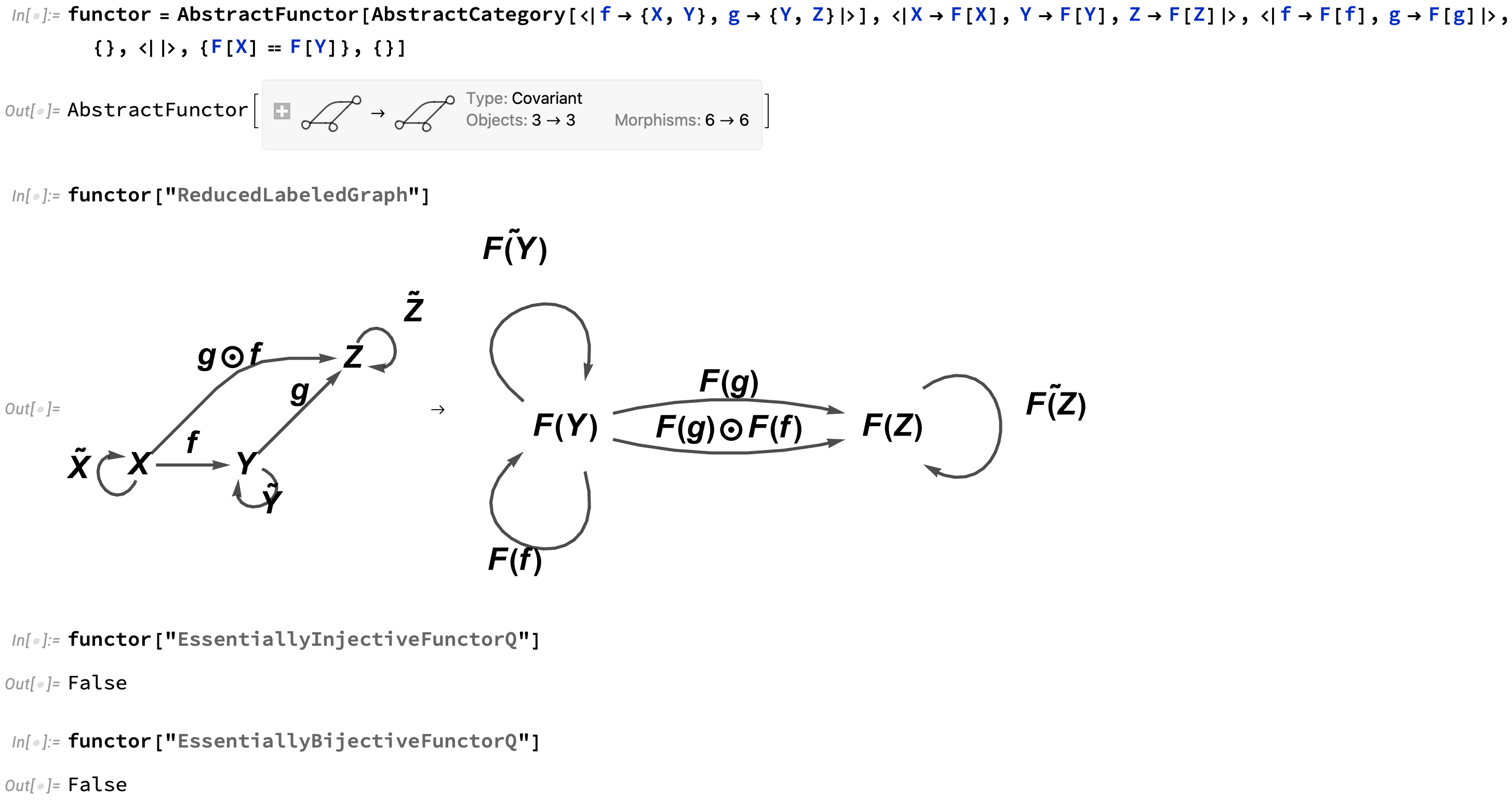}
\vrule
\includegraphics[width=0.495\textwidth]{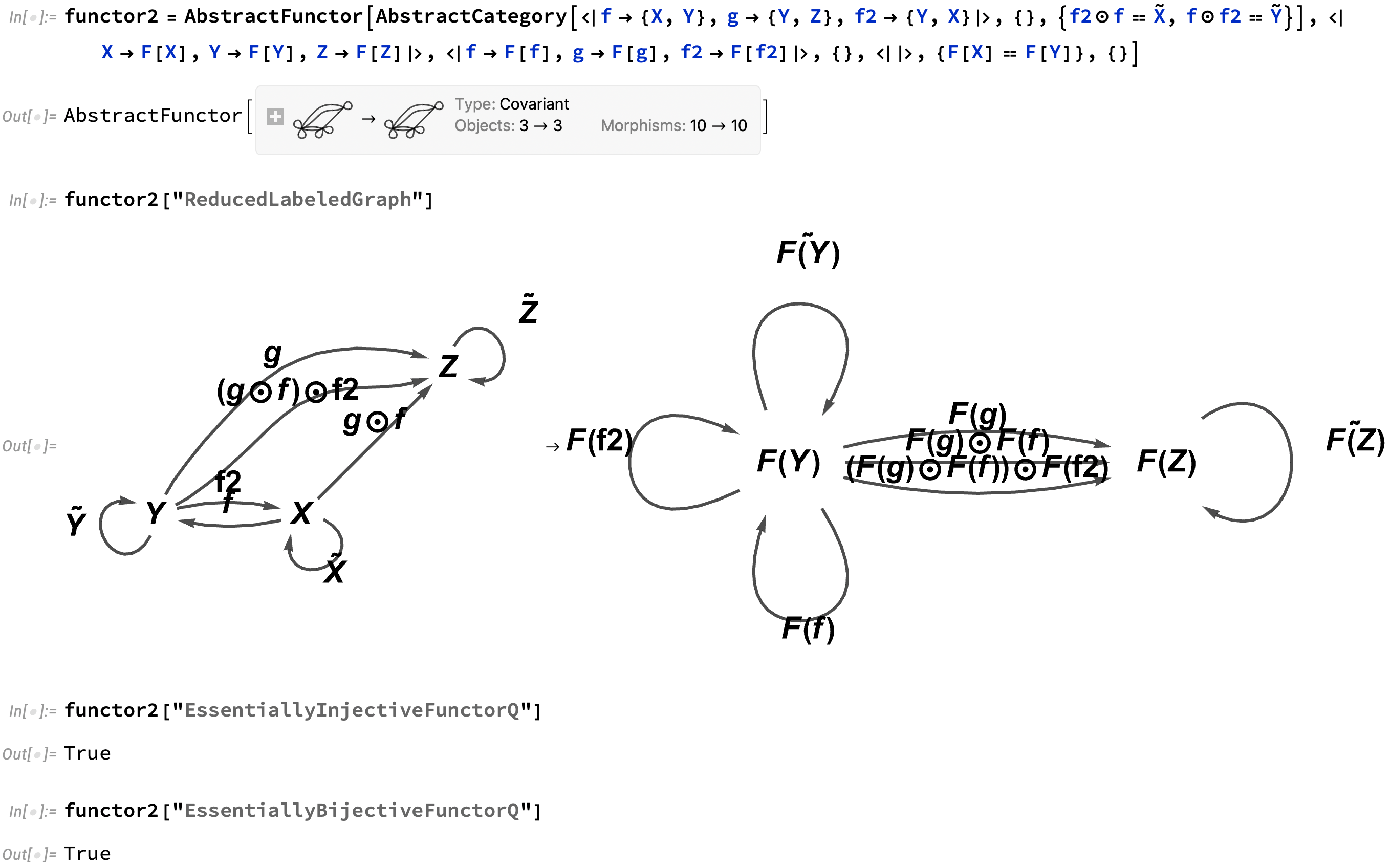}
\end{framed}
\caption{On the left, an \texttt{AbstractFunctor} object which is not essentially injective, and hence not essentially bijective, because of the additional algebraic equivalence ${F \left( X \right) = F \left( Y \right)}$ being imposed on the objects of the codomain category. On the right, the corresponding \texttt{AbstractFunctor} object with essential injectivity, and hence essential bijectivity, restored, by introducing a new isomorphism ${X \cong Y}$ between objects of the domain category.}
\label{fig:Figure16}
\end{figure}

\begin{figure}[ht]
\centering
\begin{framed}
\includegraphics[width=0.495\textwidth]{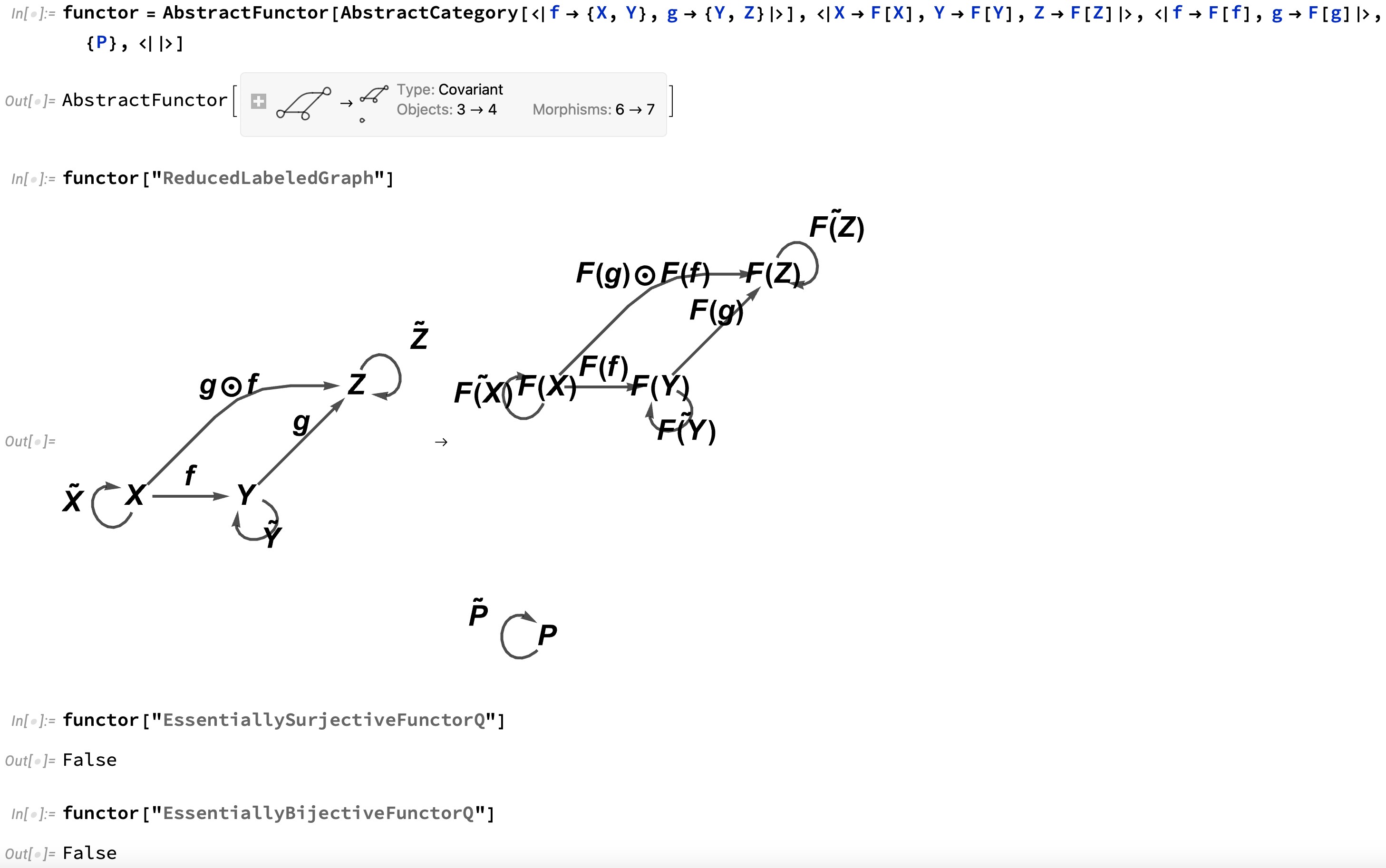}
\vrule
\includegraphics[width=0.495\textwidth]{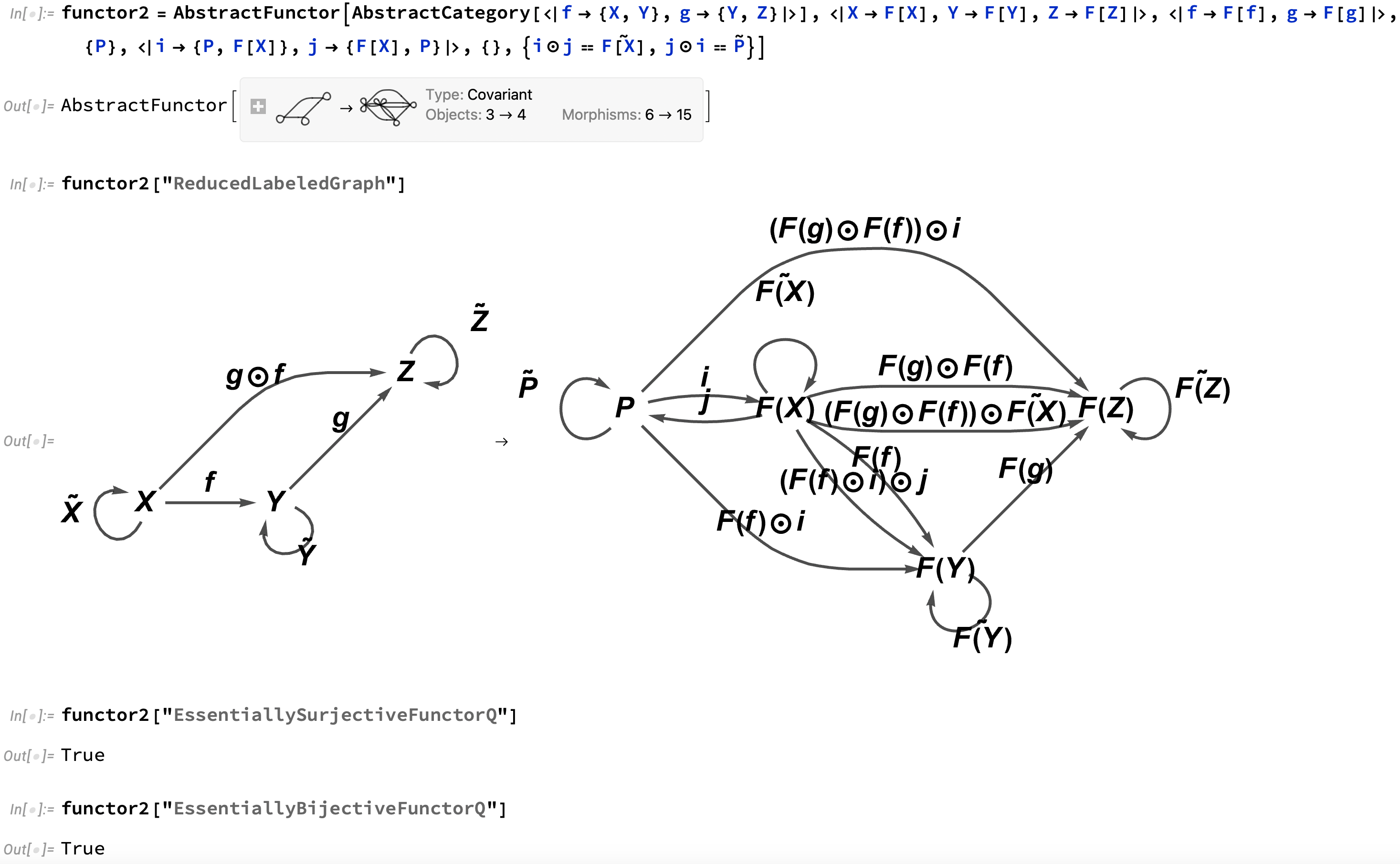}
\end{framed}
\caption{On the left, an \texttt{AbstractFunctor} object which is not essentially surjective, and hence not essentially bijective, because of the additional object $P$ being introduced within the codomain category. On the right, the corresponding \texttt{AbstractFunctor} object with essential sujrectivity, and hence essential bijectivity, restored, by introducing a new isomorphism ${P \cong F \left( X \right)}$ between objects of the codomain category.}
\label{fig:Figure17}
\end{figure}

However, there may nevertheless be cases in which a functor is only injective, surjective or bijective on objects \textit{up to isomorphism}, which in turn motivates the concepts of \textit{essential injectivity}, \textit{essential surjectivity} and \textit{essential bijectivity}. For instance, consider again the example presented above, in which a functor fails to be injective on objects because of an additional algebraic condition ${F \left( X \right) = F \left( Y \right)}$ that has been imposed on the objects in the codomain category ${\mathcal{D}}$:

\begin{equation}
\begin{tikzcd}
& Y \arrow[dr, "g"] \arrow[loop above, "id_Y"] &\\
X \arrow[ur, "f"] \arrow[rr, swap, "g \circ f"] \arrow[loop below, "id_X"] & & Z \arrow[loop below, "id_Z"]
\end{tikzcd} \qquad \mapsto \qquad
\begin{tikzcd}
F \left( X \right) = F \left( Y \right) \arrow[rr, bend left, "F \left( g \right)"] \arrow[rr, bend right, swap, "F \left( g \right) \circ F \left( f \right)"] \arrow[loop above, "id_{F \left( X \right)} = id_{F \left( Y \right)}"] \arrow[loop below, "F \left( f \right)"] & & F \left( Z \right) \arrow[loop right, "id_{F \left( Z \right)}"]
\end{tikzcd}.
\end{equation}
Now, rather than imposing the strict algebraic equality ${X = Y}$ on the objects in the domain category ${\mathcal{C}}$, suppose instead that one introduces a new inverse morphism ${f^{-1} : Y \to X}$ in ${\mathcal{C}}$, such that:

\begin{equation}
\left( f \circ f^{-1} : Y \to Y \right) = \left( id_Y : Y \to Y \right), \qquad \text{ and } \qquad \left( f^{-1} \circ f : X \to X \right) = \left( id_X : X \to X \right),
\end{equation}
and therefore the objects $X$ and $Y$ become isomorphic, thus yielding:

\begin{equation}
\begin{tikzcd}
& Y \arrow[dr, "g"] \arrow[dl, bend right, swap, "f^{-1}"] \arrow[loop above, "f \circ f^{-1} = id_Y"] &\\
X \arrow[ur, swap, "f"] \arrow[rr, swap, "g \circ f"] \arrow[loop below, "f^{-1} \circ f = id_X"] & & Z \arrow[loop below, "id_Z"]
\end{tikzcd},
\end{equation}
in the domain category, which maps to:

\begin{equation}
\mapsto \qquad \begin{tikzcd}
F \left( X \right) = F \left( Y \right) \arrow[rr, bend left, "F \left( g \right)"] \arrow[rr, bend right, swap, "F \left( g \right) \circ F \left( f \right)"] \arrow[loop above, "id_{F \left( X \right)} = id_{F \left( Y \right)}"] \arrow[loop below, "F \left( f \right)"] \arrow[loop left, "F \left( f^{-1} \right)"] & & F \left( Z \right) \arrow[loop right, "id_{F \left( Z \right)}"]
\end{tikzcd},
\end{equation}
in the codomain category. Then, the resulting functor would still not be strictly injective, nor strictly bijective, on objects, but it would be \textit{essentially} injective, and hence \textit{essentially} bijective (i.e. injective/bijective up to isomorphism). Figure \ref{fig:Figure16} implements this example in \textsc{Categorica}, showing how essential injectivity may be reinstated by imposing the isomorphism ${X \cong Y}$ on objects in the domain category. Dually, consider again the previous example of a functor that fails to be surjective on objects because of a new object $P$ that has been introduced in the codomain category:

\begin{equation}
\begin{tikzcd}
& Y \arrow[dr, "g"] \arrow[loop above, "id_Y"] &\\
X \arrow[ur, "f"] \arrow[rr, swap, "g \circ f"] \arrow[loop below, "id_X"] & & Z \arrow[loop below, "id_Z"]
\end{tikzcd} \qquad \mapsto \qquad
\begin{tikzcd}
& F \left( Y \right) \arrow[dr, "F \left( g \right)"] \arrow[loop above, "id_{F \left( Y \right)}"] & & P \arrow[loop above, "id_P"]\\
F \left( X \right) \arrow[ur, "F \left( f \right)"] \arrow[rr, swap, "F \left( g \right) \circ F \left( f \right)"] \arrow[loop below, "id_{F \left( X \right)}"] & & F \left( Z \right) \arrow[loop below, "id_{F \left( Z \right)}"]
\end{tikzcd}.
\end{equation}
Similarly, rather than imposing the strict algebraic equality ${P = F \left( X \right)}$  on the objects in the codomain category ${\mathcal{D}}$, we can instead introduce a new pair of morphisms ${i : P \to F \left( X \right)}$ and ${j : F \left( X \right) \to P}$, such that:

\begin{multline}
\left( i \circ j : F \left( X \right) \to F \left( X \right) \right) = \left( id_{F \left( X \right)} : F \left( X \right) \to F \left( X \right) \right),\\
\text{ and } \qquad \left( j \circ i : P \to P \right) = \left( id_P : P \to P \right),
\end{multline}
and therefore force the objects $P$ and ${F \left( X \right)}$ to be isomorphic, hence yielding:

\begin{equation}
\begin{tikzcd}
& Y \arrow[dr, "g"] \arrow[loop above, "id_Y"] &\\
X \arrow[ur, "f"] \arrow[rr, swap, "g \circ f"] \arrow[loop below, "id_X"] & & Y \arrow[loop below, "id_Z"]
\end{tikzcd} \qquad \mapsto \qquad
\begin{tikzcd}
P \arrow[dd, "i"] \arrow[dr, bend left = 40, "F \left( f \right) \circ i"] \arrow[ddrr, bend left = 60, "\left( F \left( g \right) \circ F \left( f \right) \right) \circ i"] \arrow[loop above, "j \circ i = id_P"] & & &\\
& F \left( Y \right) \arrow[dr, "F \left( g \right)"] \arrow[loop above, "id_{F \left( Y \right)}"] &\\
F \left( X \right) \arrow[ur, "F \left( f \right)"] \arrow[rr, swap, "F \left( g \right) \circ F \left( f \right)"] \arrow[uu, bend left, "j"] \arrow[loop below, "i \circ j = id_{F \left( X \right)}"] & & F \left( Z \right) \arrow[loop below, "id_{F \left( Z \right)}"]
\end{tikzcd}.
\end{equation}
Now, the resulting functor is still not strictly surjective, nor strictly bijective, on objects, but it is \textit{essentially} surjective, and hence \textit{essentially} bijective (i.e. surjective/bijective up to isomorphism). Figure \ref{fig:Figure17} implements this dual example in \textsc{Categorica}, showing how essential surjectivity may be reinstated by imposing the isomorphism ${P \cong F \left( X \right)}$ on objects in the codomain category.

\begin{figure}[ht]
\centering
\begin{framed}
\includegraphics[width=0.495\textwidth]{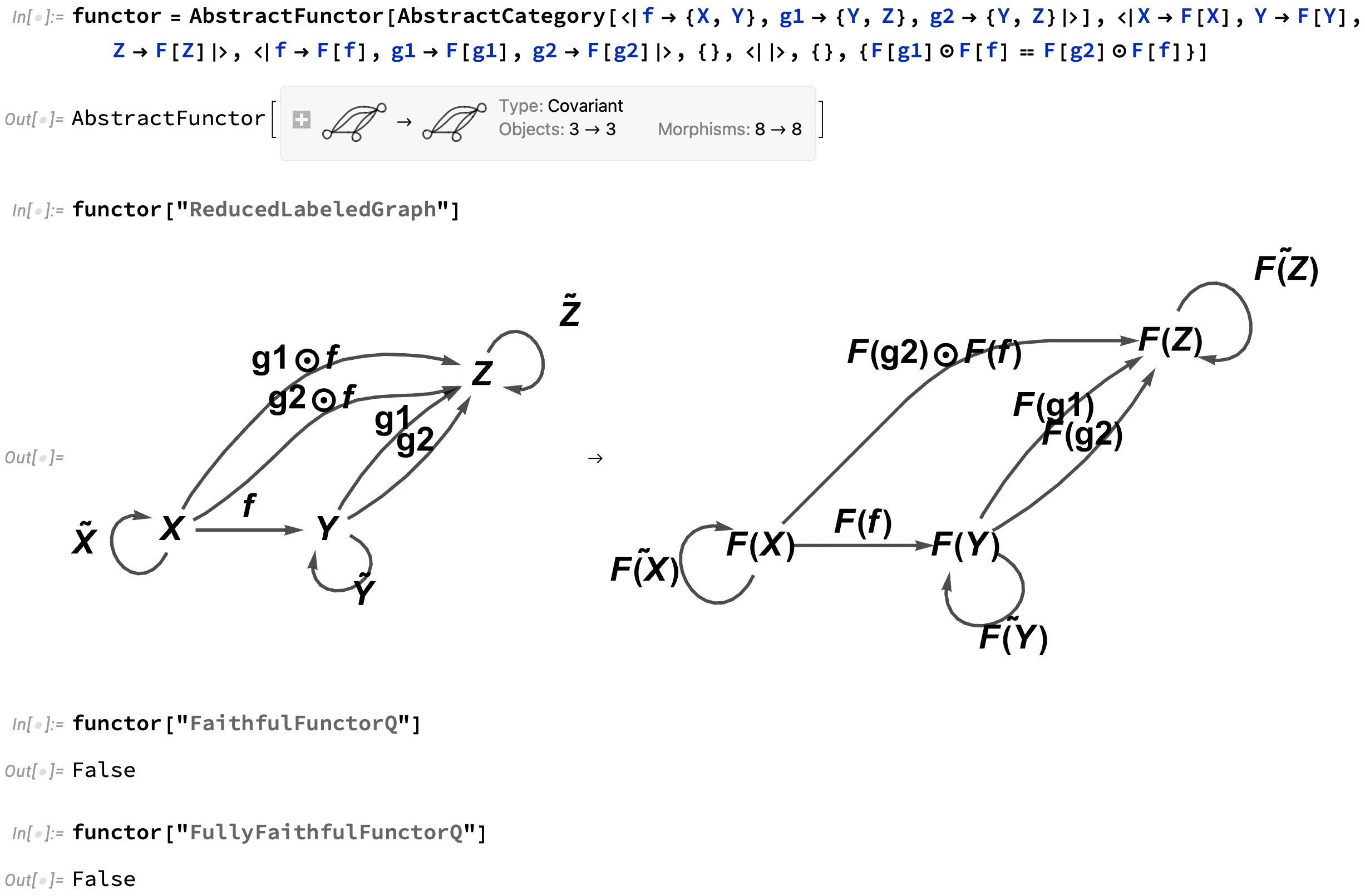}
\vrule
\includegraphics[width=0.495\textwidth]{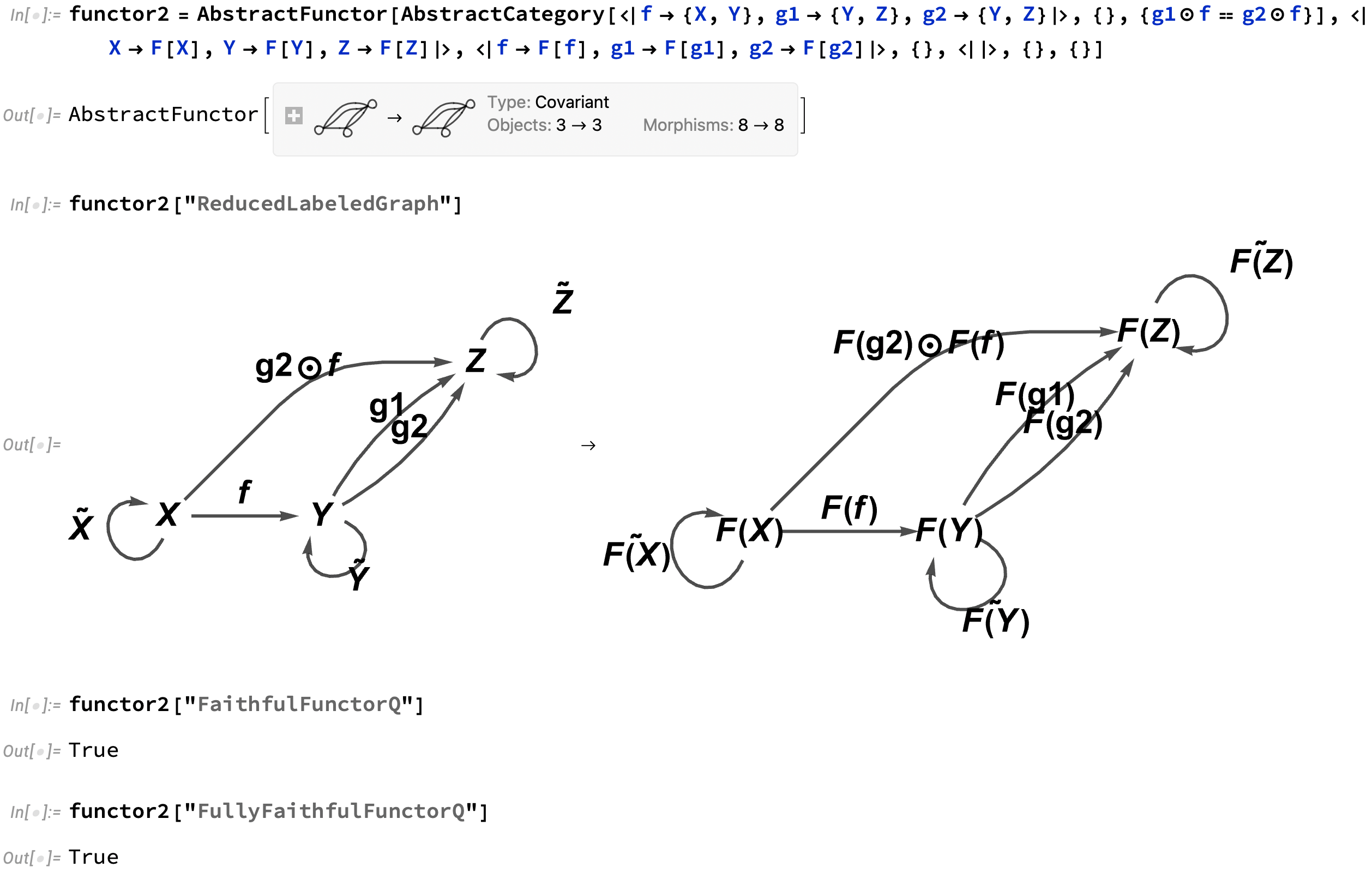}
\end{framed}
\caption{On the left, an \texttt{AbstractFunctor} object which is not faithful, and hence not fully faithful, because of the additional algebraic equivalence ${F \left( g_1 \right) \circ F \left( f \right) = F \left( g_2 \right) \circ F \left( f \right)}$ being imposed on the morphisms of the codomain category. On the right, the corresponding \texttt{AbstractFunctor} object with faithfulness, and hence full faithfulness, restored, by imposing the additional algebraic equivalence ${g_1 \circ f = g_2 \circ f}$ on the morphisms of the domain category.}
\label{fig:Figure18}
\end{figure}

\begin{figure}[ht]
\centering
\begin{framed}
\includegraphics[width=0.495\textwidth]{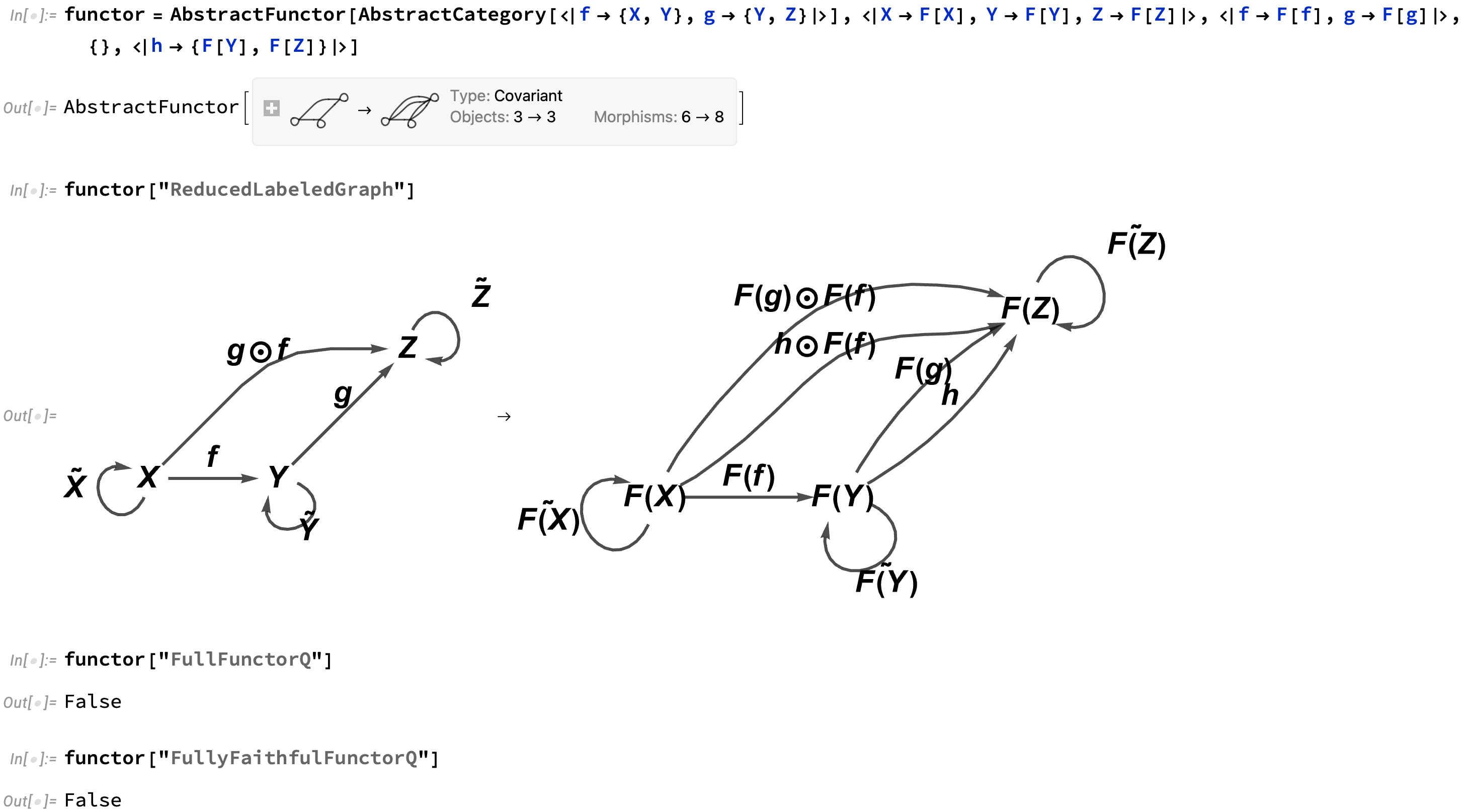}
\vrule
\includegraphics[width=0.495\textwidth]{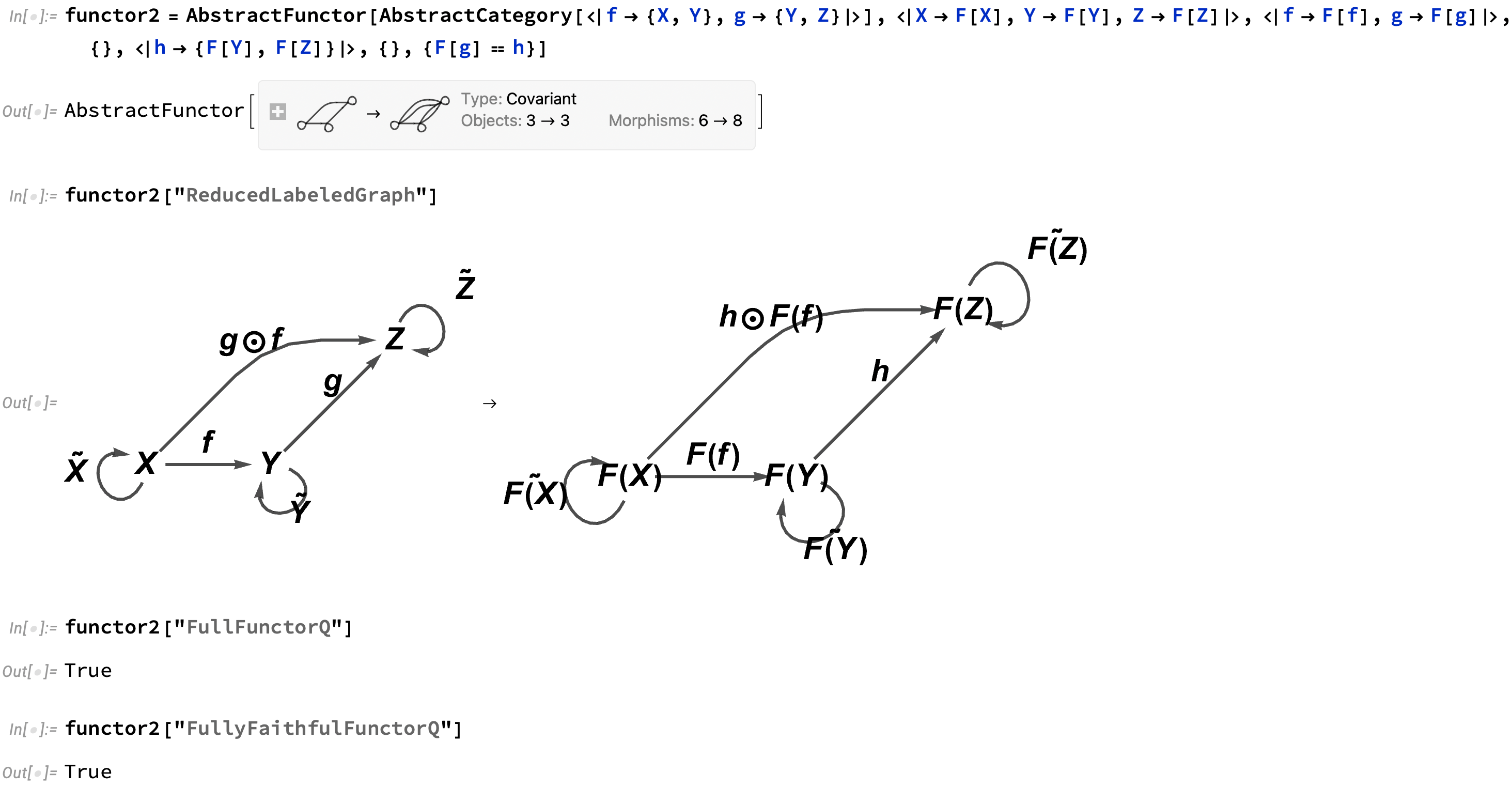}
\end{framed}
\caption{On the left, an \texttt{AbstractFunctor} object which is not full, and hence not fully faithful, because of the additional morphism $h$ being introduced within the codomain category. On the right, the corresponding \texttt{AbstractFunctor} object with fullness, and hence full faithfulness, restored, by imposing the additional algebraic equivalence ${F \left( g \right) = h}$ on the morphisms of the codomain category.}
\label{fig:Figure19}
\end{figure}

On the other hand, shifting now from considering the properties of the ${F_{ob}}$ function to those of the ${F_{hom}}$ function, functors that are injective, surjective or bijective on \textit{morphisms} are known as \textit{faithful} functors, \textit{full} functors or \textit{fully faithful} functors, respectively. A simple way in which a functor may fail to be faithful (i.e. fail to be injective on morphisms) is by having two morphisms in the codomain category be equivalent which were not previously equivalent in the domain category; for instance, consider the functor:

\begin{equation}
\begin{tikzcd}
& Y \arrow[dr, swap, "g_1"] \arrow[dr, bend left, "g_2"] \arrow[loop above, "id_Y"] &\\
X \arrow[ur, "f"] \arrow[rr, swap, "g_1 \circ f"] \arrow[rr, bend right, swap, "g_2 \circ f"] \arrow[loop below, "id_X"] & & Z \arrow[loop below, "id_Z"]
\end{tikzcd} \qquad \mapsto \qquad
\begin{tikzcd}
& F \left( Y \right) \arrow[dr, swap, "F \left( g_1 \right)"] \arrow[dr, bend left, "F \left( g_2 \right)"] \arrow[loop above, "id_{F \left( Y \right)}"] &\\
F \left( X \right) \arrow[ur, "F \left( f \right)"] \arrow[rr, swap, "F \left( g_1 \right) \circ F \left( f \right)"] \arrow[rr, bend right, swap, "F \left( g_2 \right) \circ F \left( f \right)"] \arrow[loop below, "id_{F \left( X \right)}"] & & F \left( Z \right) \arrow[loop below, "id_{F \left( Z \right)}"]
\end{tikzcd},
\end{equation}
plus the additional algebraic condition that:

\begin{equation}
\left( F \left( g_1 \right) \circ F \left( f \right) : F \left( X \right) \to F \left( Z \right) \right) = \left( F \left( g_2 \right) \circ F \left( f \right) : F \left( X \right) \to F \left( Z \right) \right),
\end{equation}
on the morphisms in the codomain category ${\mathcal{D}}$, such that one instead has:

\begin{equation}
\begin{tikzcd}
& Y \arrow[dr, swap, "g_1"] \arrow[dr, bend left, "g_2"] \arrow[loop above, "id_Y"] &\\
X \arrow[ur, "f"] \arrow[rr, swap, "g_1 \circ f"] \arrow[rr, bend right, swap, "g_2 \circ f"] \arrow[loop below, "id_X"] & & Z \arrow[loop below, "id_Z"]
\end{tikzcd} \qquad \mapsto \qquad
\begin{tikzcd}
& F \left( Y \right) \arrow[dr, swap, "F \left( g_1 \right)"] \arrow[dr, bend left, "F \left( g_2 \right)"] \arrow[loop above, "id_{F \left( Y \right)}"] &\\
F \left( X \right) \arrow[ur, "F \left( f \right)"] \arrow[rr, swap, "F \left( g_1 \right) \circ F \left( f \right) = F \left( g_2 \right) \circ F \left( f \right)"] \arrow[loop below, "id_{F \left( X \right)}"] & & F \left( Z \right) \arrow[loop below, "id_{F \left( Z \right)}"]
\end{tikzcd}.
\end{equation}
Such a functor would no longer be faithful (i.e. no longer injective on morphisms), and hence also no longer fully faithful (i.e. no longer bijective on morphisms), although its faithfulness, and therefore full faithfulness, could nevertheless be restored by introducing a corresponding algebraic condition:

\begin{equation}
\left( g_1 \circ f : X \to Z \right) = \left( g_2 \circ f : X \to Z \right),
\end{equation}
on the morphisms in the domain category ${\mathcal{C}}$, such that one now has:

\begin{equation}
\begin{tikzcd}
& Y \arrow[dr, swap, "g_1"] \arrow[dr, bend left, "g_2"] \arrow[loop above, "id_Y"] &\\
X \arrow[ur, "f"] \arrow[rr, swap, "g_1 \circ f = g_2 \circ f"] \arrow[loop below, "id_X"] & & Z \arrow[loop below, "id_Z"]
\end{tikzcd} \qquad \mapsto \qquad
\begin{tikzcd}
& F \left( Y \right) \arrow[dr, swap, "F \left( g_1 \right)"] \arrow[dr, bend left, "F \left( g_2 \right)"] \arrow[loop above, "id_{F \left( Y \right)}"] &\\
F \left( X \right) \arrow[ur, "F \left( f \right)"] \arrow[rr, swap, "F \left( g_1 \right) \circ F \left( f \right) = F \left( g_ 2 \right) \circ F \left( f \right)"] \arrow[loop below, "id_{F \left( Z \right)}"] & & F \left( Z \right) \arrow[loop below, "id_{F \left( Z \right)}"]
\end{tikzcd}.
\end{equation}
Figure \ref{fig:Figure18} implements this elementary example in \textsc{Categorica}, illustrating how the faithfulness of a functor may be removed (by imposing the algebraic equivalence ${F \left( g_1 \right) \circ F \left( f \right) = F \left( g_2 \right) \circ F \left( f \right)}$ on morphisms in the codomain category) and then subsequently reinstated (by imposing the additional algebraic equivalence ${g_1 \circ f = g_2 \circ f}$ on morphisms in the domain category). Dual to this construction, a simple way in which a functor may fail to be full (i.e. fail to be surjective on morphisms) is by having a new morphism be introduced in the codomain category that did not previously exist in the domain category; for instance, consider now the functor:

\begin{equation}
\begin{tikzcd}
& Y \arrow[dr, "g"] \arrow[loop above, "id_Y"] &\\
X \arrow[ur, "f"] \arrow[rr, swap, "g \circ f"] \arrow[loop below, "id_X"] & & Z \arrow[loop below, "id_Z"]
\end{tikzcd} \qquad \mapsto \qquad
\begin{tikzcd}
& F \left( Y \right) \arrow[dr, swap, "F \left( g \right)"] \arrow[dr, bend left, "h"] \arrow[loop above, "id_{F \left( Y \right)}"] &\\
F \left( X \right) \arrow[ur, "F \left( f \right)"] \arrow[rr, swap, "F \left( g \right) \circ F \left( f \right)"] \arrow[rr, bend right, swap, "h \circ F \left( f \right)"] \arrow[loop below, "id_{F \left( X \right)}"] & & F \left( Z \right) \arrow[loop below, "id_{F \left( Z \right)}"]
\end{tikzcd}.
\end{equation}
Such a functor would no longer be full (i.e. no longer surjective on morphisms), and hence also no longer fully faithful (i.e. no longer bijective on morphisms), although its fullness, and hence full faithfulness, could nevertheless be restored by introducing a new algebraic condition:

\begin{equation}
\left( F \left( g \right) : F \left( Y \right) \to F \left( Z \right) \right) = \left( h : F \left( Y \right) \to F \left( Z \right) \right),
\end{equation}
on morphisms in the codomain category ${\mathcal{D}}$, such that one now has:

\begin{equation}
\begin{tikzcd}
& Y \arrow[dr, "g"] \arrow[loop above, "id_Y"] &\\
X \arrow[ur, "f"] \arrow[rr, swap, "g \circ f"] \arrow[loop below, "id_X"] & & Z \arrow[loop below, "id_Z"]
\end{tikzcd} \qquad \mapsto \qquad
\begin{tikzcd}
& F \left( Y \right) \arrow[dr, "F \left( g \right) = h"] \arrow[loop above, "id_{F \left( Y \right)}"] &\\
F \left( X \right) \arrow[ur, "F \left( f \right)"] \arrow[rr, swap, "F \left( g \right) \circ F \left( f \right) = h \circ F \left( f \right)"] \arrow[loop below, "id_{F \left( X \right)}"] & & F \left( Z \right) \arrow[loop below, "id_{F \left( Z \right)}"]
\end{tikzcd}.
\end{equation}
Figure \ref{fig:Figure19} implements this dual example in \textsc{Categorica}, illustrating how the fullness of a functor may be removed (by introducing the new morphism ${h : F \left( Y \right) \to F \left( Z \right)}$ in the codomain category) and then subsequently reinstated (by imposing the algebraic equivalence ${F \left( g \right) = h}$ on morphisms in the codomain category).

The final example of a functorial construction that we shall cover within this section, due to its rather foundational significance in the fields of algebraic geometry, algebraic topology and (higher) homotopy theory, is that of a (Grothendieck) \textit{fibration}\cite{grothendieck2}\cite{gray}. This purely category-theoretic notion of a fibration is effectively a grand generalization of the notion of a fiber bundle (or topological fibration) in topology, wherein one considers a functor of the general form ${F : \mathcal{E} \to \mathcal{B}}$, with the domain category ${\mathcal{E}}$ being the \textit{total category} of the fibration (thus playing an analogous role to that of the \textit{total space} of a fiber bundle) and the codomain category ${\mathcal{B}}$ being the \textit{base category} (thus playing an analogous role to that of the \textit{base space} of the same fiber bundle). For each object $X$ in the base category ${\mathcal{B}}$, one can consider the \textit{fiber category} ${F_X}$ at $X$, which is a category whose object set ${\mathrm{ob} \left( F_X \right)}$ is the set of objects in the total category ${\mathcal{E}}$ that map to $X$ under the functor $F$ (i.e. the preimage of $X$ under the function ${F_{ob}}$) and whose morphism set ${\mathrm{hom} \left( F_X \right)}$ is the set of morphisms in the total category ${\mathcal{E}}$ that map to the identity morphism ${id_X : X \to X}$ on $X$ under $F$ (i.e. the preimage of ${id_X : X \to X}$ under the function ${F_{hom}}$), i.e:

\begin{multline}
\forall X \in \mathrm{ob} \left( \mathcal{B} \right), \qquad F_X \text{ is the fiber category of functor } F : \mathcal{E} \to \mathcal{B} \text{ at } X\\
\iff \mathrm{ob} \left( F_X \right) = \left\lbrace U \in \mathrm{ob} \left( \mathcal{E} \right) : F \left( U \right) = X \right\rbrace, \qquad \text{ and }\\
\mathrm{hom} \left( F_X \right) = \left\lbrace \left( f : U \to V \right) \in \mathrm{hom} \left( \mathcal{E} \right) : \left( F \left( f \right) : F \left( U \right) \to F \left( V \right) \right) = \left( id_X : X \to X \right) \right\rbrace.
\end{multline}
In order to qualify as a valid Grothendieck fibration, however, the functor ${F : \mathcal{E} \to \mathcal{B}}$ must satisfy an additional, and rather technical, condition known as \textit{contravariant pseudofunctoriality} (for the corresponding dual construction, known as a Grothendieck \textit{opfibration}, the relevant condition is \textit{covariant pseudofunctoriality}), which is typically expressed in terms of \textit{Cartesian morphisms}. The formal definition of a Cartesian morphism is unfortunately a little opaque; specifically, if the morphism ${f : X \to Y}$ in the total category ${\mathcal{E}}$ is such that, for all objects $Z$ and all morphisms ${h : Z \to Y}$ in ${\mathcal{E}}$, and all morphisms ${u : F \left( Z \right) \to F \left( X \right)}$ in the base category ${\mathcal{B}}$ such that ${\left( F \left( f \right) \circ u : F \left( Z \right) \to F \left( Y \right) \right) = \left( F \left( h \right) : F \left( Z \right) \to F \left( Y \right) \right)}$, there necessarily exists a unique morphism ${v : Z \to X}$ in the total category ${\mathcal{E}}$ such that ${\left( f \circ v : Z \to Y \right) = \left( h : Z \to Y \right)}$ in ${\mathcal{E}}$ and ${\left( u : F \left( Z \right) \to F \left( X \right) \right) = \left( F \left( v \right) : F \left( X \right) \to F \left( Y \right) \right)}$ in ${\mathcal{B}}$, then ${f : X \to Y}$ is Cartesian:

\begin{multline}
\forall \left( f : X \to Y \right) \in \mathrm{hom} \left( \mathcal{E} \right), \qquad \left( f : X \to Y \right) \text{ is a Cartesian morphism}\\
\iff \qquad \forall Z \in \mathrm{ob} \left( \mathcal{E} \right), \qquad \forall \left( h : Z \to Y \right) \in \mathrm{hom} \left( \mathcal{E} \right), \qquad \forall \left( u : F \left( Z \right) \to F \left( X \right) \right) \in \mathrm{hom} \left( \mathcal{B} \right),\\
\left( F \left( f \right) \circ u : F \left( Z \right) \to F \left( Y \right) \right) = \left( F \left( h \right) : F \left( Z \right) \to F \left( Y \right) \right) \qquad \implies \qquad \exists! \left( v : Z \to X \right) \in \mathrm{hom} \left( \mathcal{E} \right),\\
\text{ such that } \qquad \left( f \circ v : Z \to Y \right) = \left( h : Z \to Y \right), \qquad \text{ and }\\
\left( u : F \left( Z \right) \to F \left( X \right) \right) = \left( F \left( v \right) : F \left( Z \right) \to F \left( X \right) \right),
\end{multline}
or, illustrated diagrammatically, one has the following basic setup:

\begin{equation}
\begin{tikzcd}
\forall Z \arrow[dd, dashed, swap, "\exists! v"] \arrow[ddrr, swap, "\forall h"] \arrow[ddrr, bend left, "f \circ v"] & &\\ \\
X \arrow[rr, swap, "f"] & & Y
\end{tikzcd} \qquad \mapsto \qquad
\begin{tikzcd}
F \left( \forall Z \right) \arrow[dd, swap, "\forall u"] \arrow[ddrr, swap, "F \left( \forall h \right)"] \arrow[ddrr, bend left, "F \left( f \right) \circ u"] & &\\ \\
F \left( X \right) \arrow[rr, swap, "F \left( f \right)"] & & F \left( Y \right)
\end{tikzcd},
\end{equation}
which then collapses down to:

\begin{equation}
\begin{tikzcd}
\forall Z \arrow[dd, dashed, swap, "\exists! v"] \arrow[ddrr, "f \circ v = \forall h"] & &\\ \\
X \arrow[rr, swap, "f"] & & Y
\end{tikzcd} \qquad \mapsto \qquad
\begin{tikzcd}
F \left( \forall Z \right) \arrow[dd, swap, "u = F \left( \exists! v \right)"] \arrow[ddrr, "F \left( f \right) \circ u = F \left( \forall h \right)"] & &\\ \\
F \left( X \right) \arrow[rr, swap, "F \left( f \right)"] & & F \left( Y \right)
\end{tikzcd}.
\end{equation}
The contravariant pseudofunctoriality condition that characterizes Grothendieck fibrations is then the condition that, for every object $Y$ in the total category ${\mathcal{E}}$, and every morphism ${f_0 : X_0 \to F \left( Y \right)}$ mapping to the image ${F \left( Y \right)}$ of $Y$ in the base category ${\mathcal{B}}$, there must exist some Cartesian morphism ${f : X \to Y}$ in the total category ${\mathcal{E}}$ such that the morphism ${f_0 : X_0 \to F \left( Y \right)}$ is the image of the Cartesian morphism ${f : X \to Y}$ in the base category ${\mathcal{B}}$, i.e. ${\left( F \left( f \right) : F \left( X \right) \to F \left( Y \right) \right) = \left( f_0 : X_0 \to F \left( Y \right) \right)}$ (and therefore also ${F \left( X \right) = X_0}$):

\begin{multline}
\forall Y \in \mathrm{ob} \left( \mathcal{E} \right), \qquad \forall X_0 \in \mathrm{ob} \left( \mathcal{B} \right), \qquad \forall \left( f_0 : X_0 \to F \left( Y \right) \right) \in \mathrm{hom} \left( \mathcal{B} \right),\\
\exists X \in \mathrm{ob} \left( \mathcal{E} \right), \qquad \exists \left( f : X \to Y \right), \qquad \text{ such that } \qquad \left( f : X \to Y \right) \text{ is Cartesian}, \qquad \text{ and}\\
F \left( X \right) = X_0, \qquad \text{ and } \qquad \left( F \left( f \right) : X_0 \to F \left( Y \right) \right) = \left( f_0 : X_0 \to F \left( Y \right) \right).
\end{multline}
There exist many notable specializations of this highly abstract Grothendieck fibration definition, including the case of \textit{discrete fibrations} (in which every fiber category ${F_X}$ is a discrete category, consisting solely of objects and their identity morphisms) and \textit{groupoidal fibrations} (in which every fiber category ${F_X}$ is a groupoid). \textsc{Categorica} does not yet possess the functionality to detect and characterize all cases of Grothendieck fibrations with complete generality, although important special cases (such as discrete fibrations, via the \textit{``DiscreteFibrationQ''} property, etc.) have indeed been implemented fully.

\begin{figure}[ht]
\centering
\begin{framed}
\includegraphics[width=0.515\textwidth]{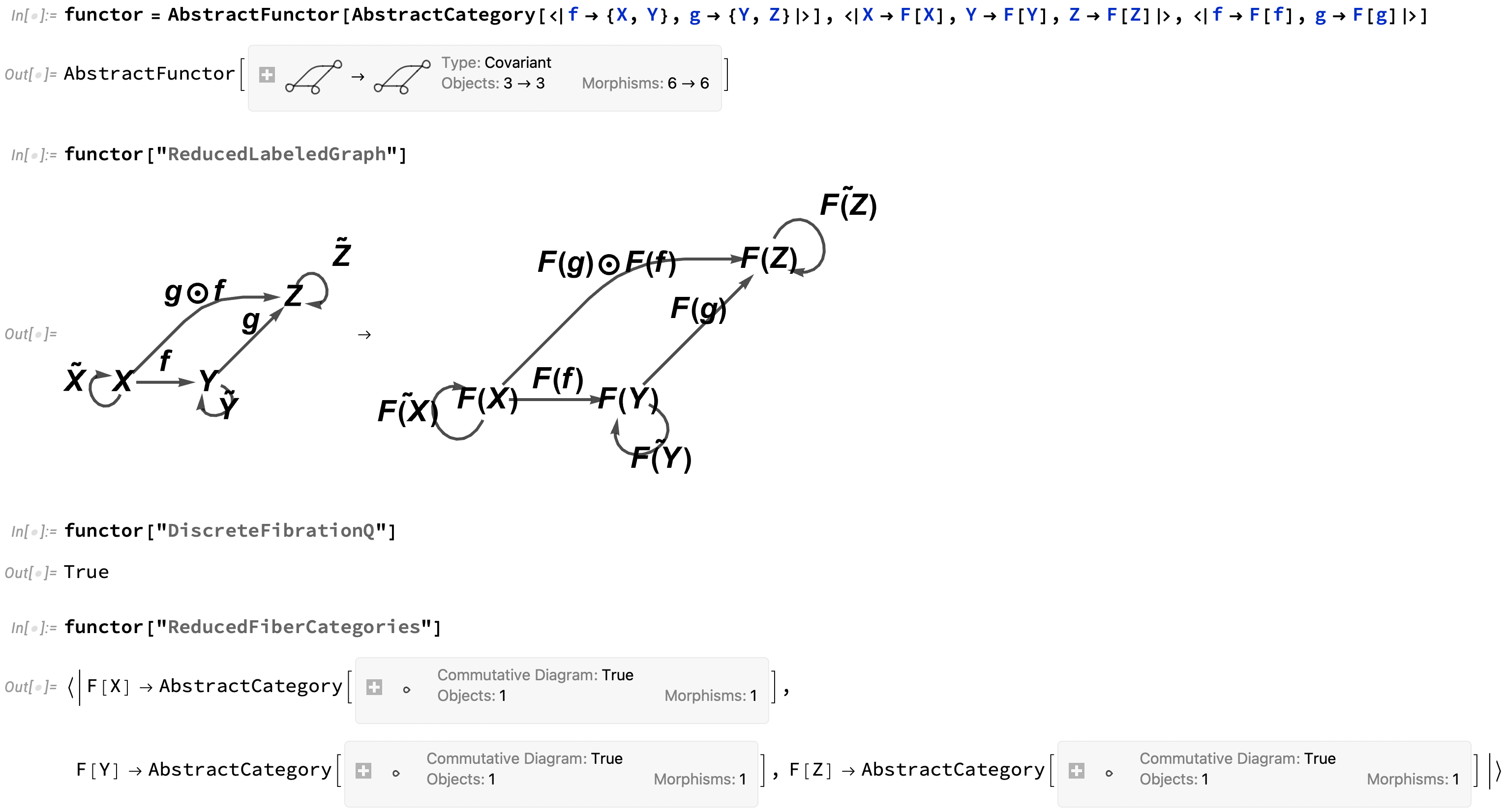}
\vrule
\includegraphics[width=0.475\textwidth]{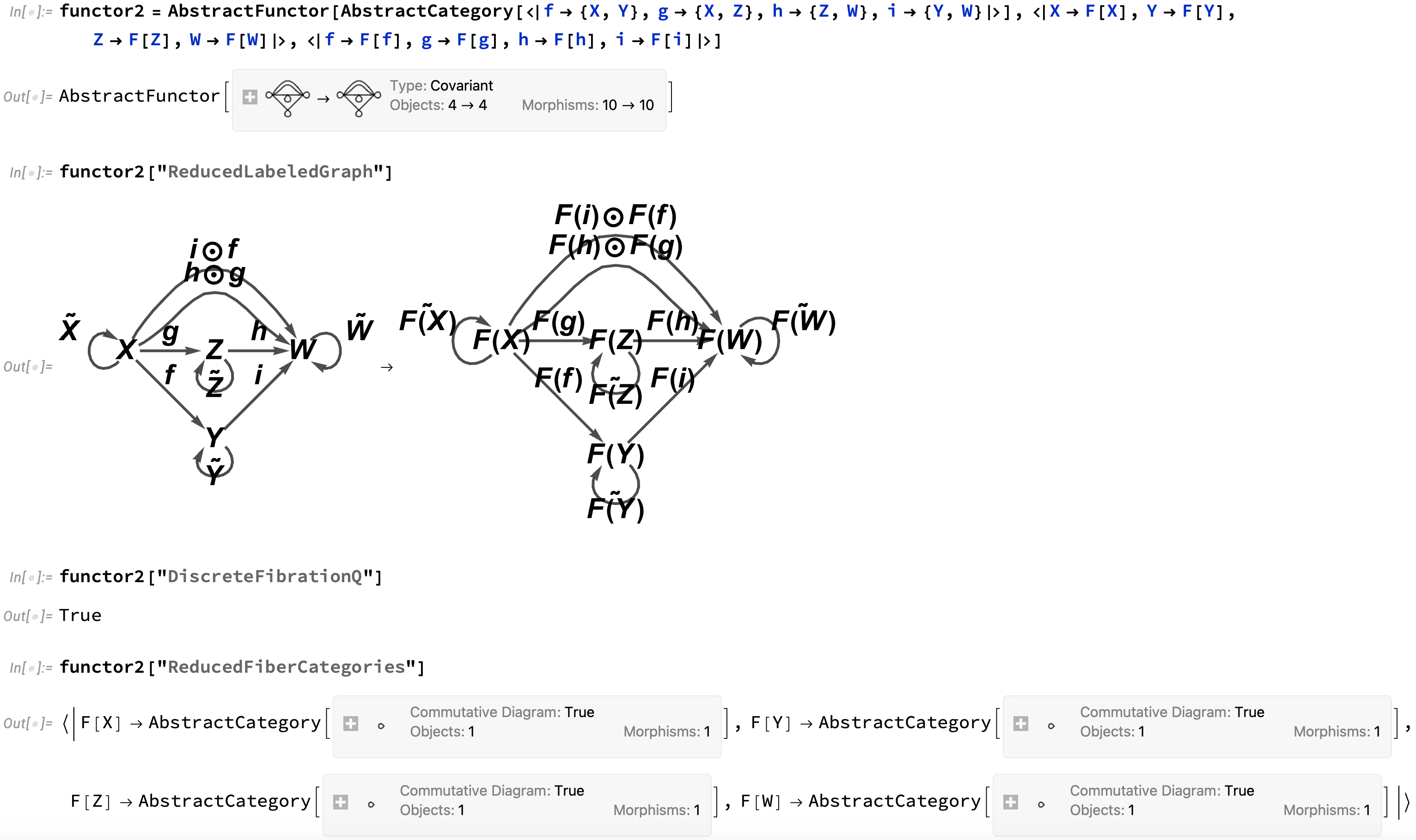}
\end{framed}
\caption{On the left, an \texttt{AbstractFunctor} object corresponding to a discrete fibration of a simple three-object, six-morphism total category over a three-object, six-morphism base category, showing the three discrete (single-object) fiber categories. On the right, an \texttt{AbstractFunctor} object corresponding to a discrete fibration of a slightly larger four-object, ten-morphism total category over a four-object, ten-morphism base category, showing the four discrete (single-object) fiber categories.}
\label{fig:Figure25}
\end{figure}

\begin{figure}[ht]
\centering
\begin{framed}
\includegraphics[width=0.495\textwidth]{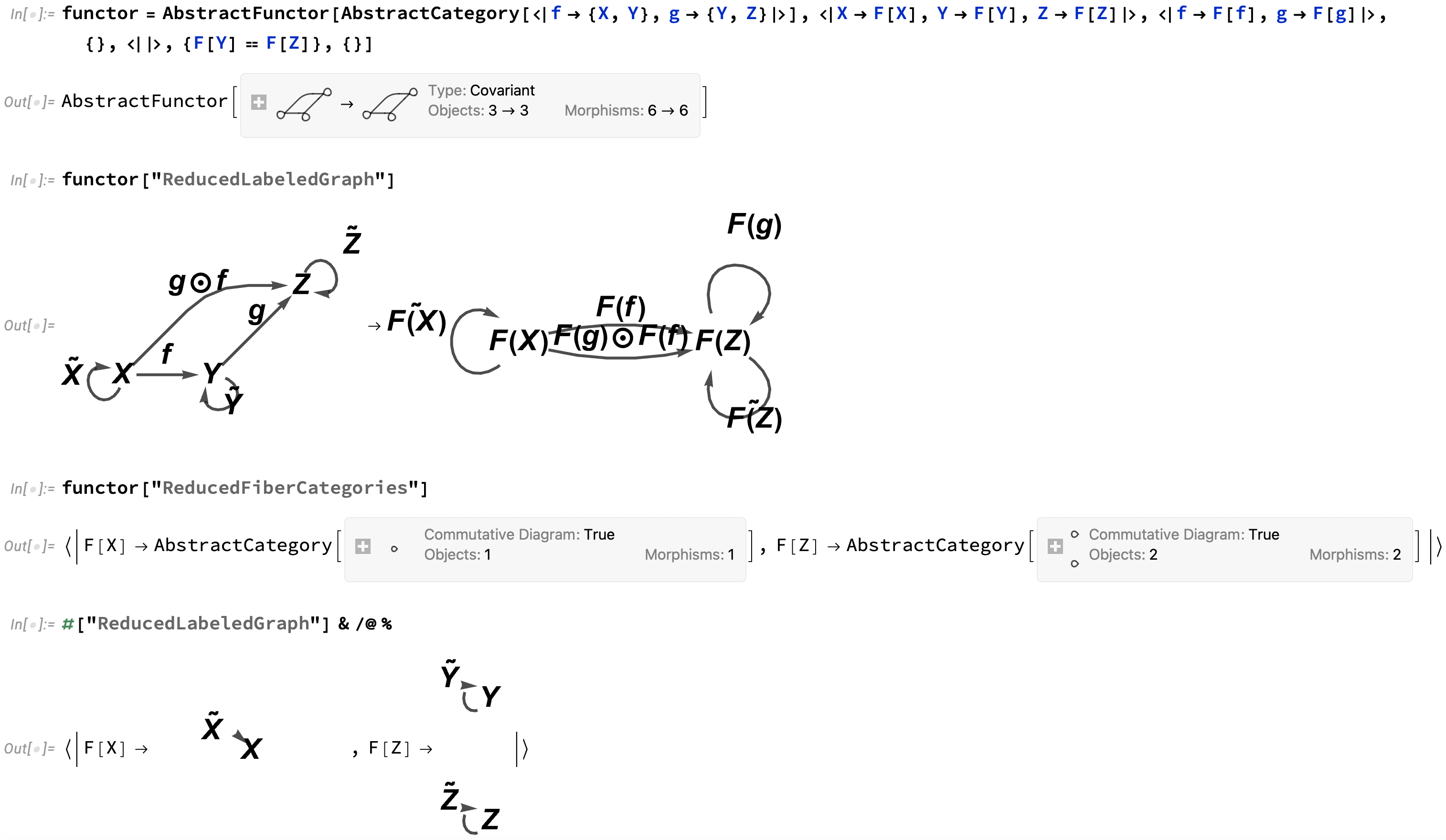}
\vrule
\includegraphics[width=0.495\textwidth]{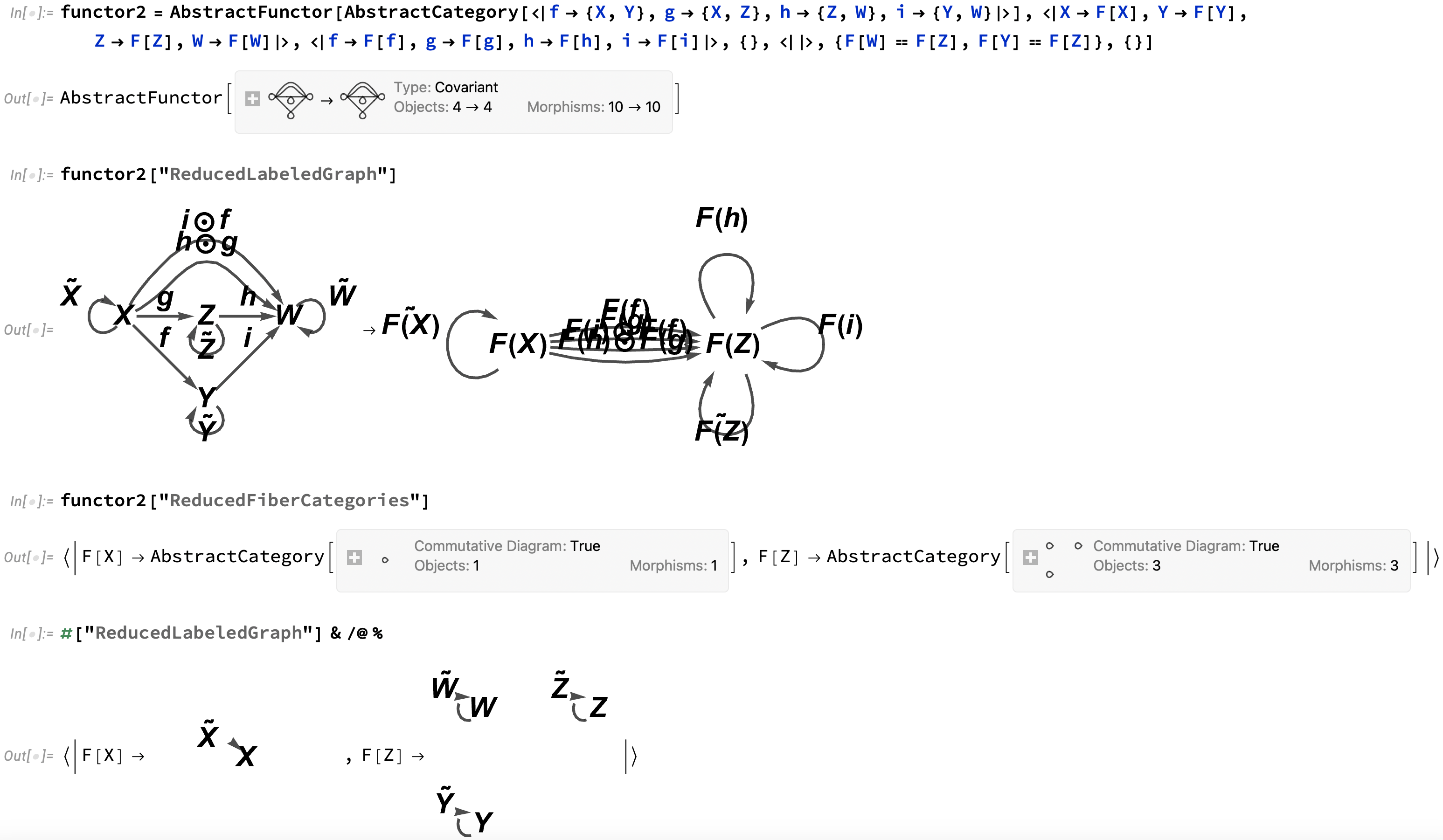}
\end{framed}
\caption{On the left, an \texttt{AbstractFunctor} object corresponding to a more general fibration of a simple three-object, six-morphism total category over a two-object, five-morphism base category, obtained by imposing the additional object equivalence ${F \left( Y \right) = F \left( Z \right)}$ in the codomain category, showing the two resulting discrete fiber categories. On the right, an \texttt{AbstractFunctor} object corresponding to a more general fibration of a slightly larger four-object, ten-morphism total category over a two-object, eight-morphism base category, obtained by impoosing the additional object equivalences ${F \left( W \right) = F \left( Z \right)}$ and ${F \left( Y \right) = F \left( Z \right)}$ in the codomain category, showing the two resulting discrete fiber categories.}
\label{fig:Figure26}
\end{figure}

\begin{figure}[ht]
\centering
\begin{framed}
\includegraphics[width=0.515\textwidth]{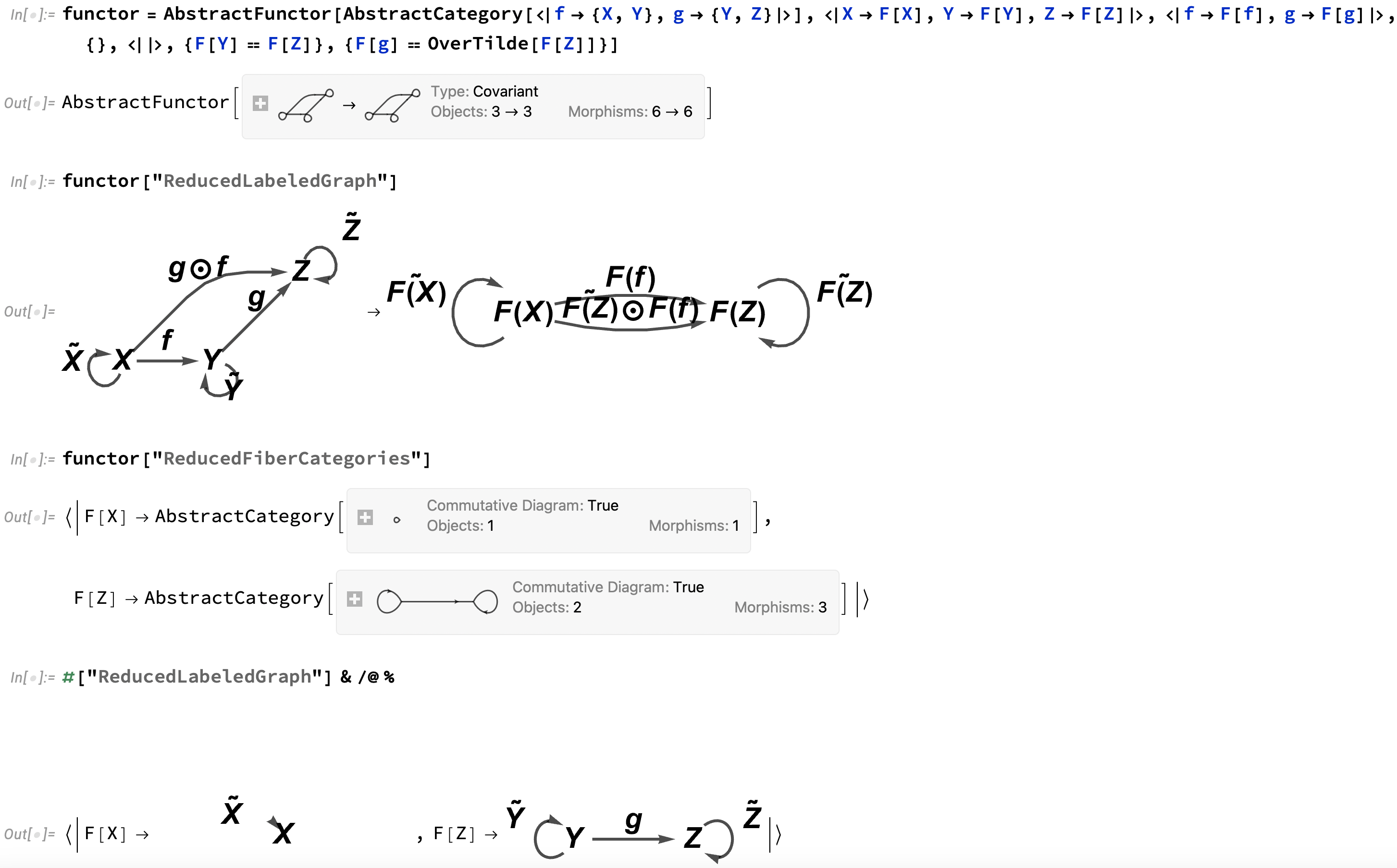}
\vrule
\includegraphics[width=0.475\textwidth]{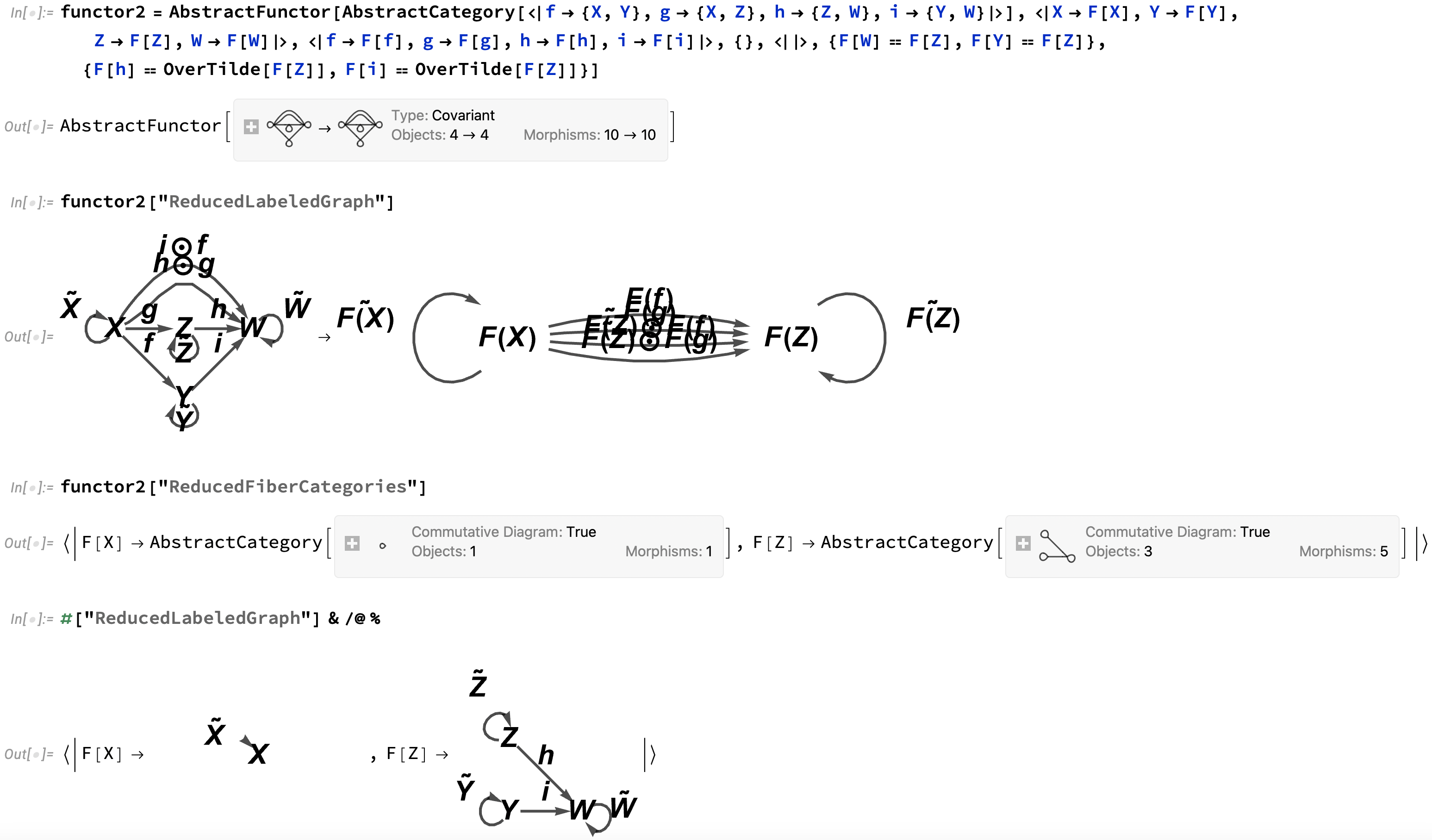}
\end{framed}
\caption{On the left, an \texttt{AbstractFunctor} object corresponding to a more general fibration of a simple three-object, six-morphism total category over a two-object, four-morphism base category, obtained by imposing the additional object equivalence ${F \left( Y  \right) = F \left( Z \right)}$ and additional morphism equivalence ${F \left( g \right) = id_{F \left( Z \right)}}$ in the codomain category, showing the two resulting non-discrete fiber categories. On the right, an \texttt{AbstractFunctor} object corresponding to a more general fibration of a slightly larger four-object, ten-morphism total category over a two-object, six-morphism base category, obtained by imposing the additional object equivalences ${F \left( W \right) = F \left( Z \right)}$ and ${F \left( Y \right) = F \left( Z \right)}$, and additional morphism equivalences ${F \left( h \right) = id_{F \left( Z \right)}}$ and ${F \left( i \right) = id_{F \left( Z \right)}}$, in the codomain category, showing the two resulting non-discrete fiber categories.}
\label{fig:Figure27}
\end{figure}

For instance, if we consider the following simple functor from a three-object, six-morphism \textit{total category} ${\mathcal{E}}$ to another three-object, six-morphism \textit{base category} ${\mathcal{B}}$:

\begin{equation}
\begin{tikzcd}
& Y \arrow[dr, "g"] \arrow[loop above, "id_Y"] &\\
X \arrow[ur, "f"] \arrow[rr, swap, "g \circ f"] \arrow[loop below, "id_X"] & & Z \arrow[loop below, "id_Z"]
\end{tikzcd} \qquad \mapsto \qquad
\begin{tikzcd}
& F \left( Y \right) \arrow[dr, "F \left( g \right)"] \arrow[loop above, "id_{F \left( Y \right)}"] &\\
F \left( X \right) \arrow[ur, "F \left( f \right)"] \arrow[rr, swap, "F \left( g \right) \circ F \left( f \right)"] \arrow[loop below, "id_{F \left( X \right)}"] & & F \left( Z \right) \arrow[loop below, "id_{F \left( Z \right)}"]
\end{tikzcd},
\end{equation}
or, alternatively, the following slightly more complex functor from a four-object, ten-morphism \textit{total category} ${\mathcal{E}}$ to another four-object, ten-morphism \textit{base category} ${\mathcal{B}}$:

\begin{equation}
\begin{tikzcd}
X \arrow[rr, "f"] \arrow[dd, swap, "g"] \arrow[ddrr, bend left = 10, "h \circ g"] \arrow[ddrr, bend right = 10, swap, "i \circ f"] \arrow[loop above, "id_X"] & & Y \arrow[dd, "i"] \arrow[loop above, "id_Y"] \\ \\
Z \arrow[rr, swap, "h"] \arrow[loop below, "id_Z"] & & W \arrow[loop below, "id_W"]
\end{tikzcd} \qquad \mapsto \qquad
\begin{tikzcd}
F \left( X \right) \arrow[rrr, "F \left( f \right)"] \arrow[ddd, swap, "F \left( g \right)"] \arrow[dddrrr, bend left = 15, "F \left( h \right) \circ F \left( g \right)"] \arrow[dddrrr, bend right = 15, swap, "F \left( i \right) \circ F \left( f \right)"] \arrow[loop above, "id_{F \left( X \right)}"] & & & F \left( Y \right) \arrow[ddd, "F \left( i \right)"] \arrow[loop above, "id_{F \left( Y \right)}"] \\ \\ \\
F \left( Z \right) \arrow[rrr, swap, "F \left( h \right)"] \arrow[loop below, "id_{F \left( Z \right)}"] & & & F \left( W \right) \arrow[loop below, "id_{F \left( W \right)}"]
\end{tikzcd},
\end{equation}
then we see that these functors both constitute more-or-less trivial cases of (discrete) fibrations. In the former case, there are three (discrete) fiber categories, one for each object ${F \left( X \right)}$, ${F \left( Y \right)}$ and ${F \left( Z \right)}$ in the base category ${\mathcal{B}}$, and each consisting of the single object $X$, $Y$ or $Z$ from the total category ${\mathcal{E}}$ (along with the corresponding identity morphism in each case), respectively; in the latter case, there are instead four (discrete) fiber categories, one for each object ${F \left( X \right)}$, ${F \left( Y \right)}$, ${F \left( Z \right)}$ and ${F \left( W \right)}$ in the base category ${\mathcal{B}}$, and each consisting of the single object $X$, $Y$, $Z$ or $W$ from the total category ${\mathcal{E}}$ (along with the corresponding identity morphism in each case), respectively. This is illustrated in Figure \ref{fig:Figure25}, in which \textsc{Categorica} correctly identifies that the relevant \texttt{AbstractFunctor} objects are both discrete fibrations, and proceeds to compute the relevant fiber categories (represented as an association of \texttt{AbstractCategory} objects - one for each object in the codomain/base category of the fibration). In the first instance, imposing the additional algebraic equivalence ${F \left( Y \right) = F \left( Z \right)}$ on objects in the codomain category now yields a fibration over a two-object, five-morphism base category:

\begin{equation}
\begin{tikzcd}
& Y \arrow[dr, "g"] \arrow[loop above, "id_Y"] &\\
X \arrow[ur, "f"] \arrow[rr, swap, "g \circ f"] \arrow[loop below, "id_X"] &  & Z \arrow[loop below, "id_Z"]
\end{tikzcd} \qquad \mapsto \qquad
\begin{tikzcd}
F \left( X \right) \arrow[loop left, "id_{F \left( X \right)}"] \arrow[rr, bend left, "F \left( f \right)"] \arrow[rr, bend right, swap, "F \left( g \right) \circ F \left( f \right)"] & & F \left( Y  \right) = F \left( Z \right) \arrow[loop above, "F \left( g \right)"] \arrow[loop below, "id_{F \left( Y\right)} = id_{F \left( Z \right)}"]
\end{tikzcd},
\end{equation}
while, in the second instance, imposing the additional algebraic equivalences ${F \left( W \right) = F \left( Z \right)}$ and ${F \left( Y \right) = F \left( Z \right)}$ on objects in the codomain category yields instead a fibration over a two-object, eight-morphism base category:

\begin{equation}
\begin{tikzcd}
X \arrow[rr, "f"] \arrow[dd, swap, "g"] \arrow[ddrr, bend left = 10, "h \circ g"] \arrow[ddrr, bend right = 10, swap, "i \circ f"] \arrow[loop above, "id_X"] & & Y \arrow[dd, "i"] \arrow[loop above, "id_Y"] \\ \\
Z \arrow[rr, swap, "h"] \arrow[loop below, "id_Z"] & & W \arrow[loop below, "id_W"]
\end{tikzcd} \qquad \mapsto \qquad
\begin{tikzcd}
F \left( X \right) \arrow[loop left, "id_{F \left( X \right)}"] \arrow[rrr, bend left = 30, "F \left( f \right)"] \arrow[rrr, bend left = 10, "F \left( g \right)"] \arrow[rrr, bend right = 10, swap, "F \left( i \right) \circ F \left( f \right)"] \arrow[rrr, bend right = 30, swap, "F \left( h \right) \circ F \left( g \right)"] & & & {\renewcommand{\arraystretch}{0.5}
\begin{tabular}{c}
${F \left( Y \right)}$\\
${= F \left( W \right)}$\\
${= F \left( Z \right)}$
\end{tabular}} \arrow[loop above, "F \left( h \right)"] \arrow[loop right, "F \left( i \right)"] \arrow[loop below, "id_{F \left( Y \right)} = id_{F \left( W \right)} = id_{F \left( Z \right)}"]
\end{tikzcd}.
\end{equation}
These fibrations are still discrete, but now in the former case there are only two discrete fiber categories: one for object ${F \left( X \right)}$ in the base category ${\mathcal{B}}$, consisting of the single object $X$ from the total category ${\mathcal{E}}$ (along with its identity morphism ${id_X : X \to X}$), and one for object ${F \left( Y \right) = F \left( Z \right)}$ in the base category ${\mathcal{B}}$, consisting of the pair of objects $Y$ and $Z$ from the total category ${\mathcal{E}}$ (along with their respective identity morphisms ${id_Y : Y \to Y}$ and ${id_Z : Z \to Z}$). In the latter case there are now also only two discrete fiber categories: one for object ${F \left( X \right)}$ in the base category ${\mathcal{B}}$, consisting of the single object $X$ from the total category ${\mathcal{E}}$ (along with its identity morphism ${id_X : X \to X}$), and one for object ${F \left( Y \right) = F \left( W \right) = F \left( Z \right)}$ in the base category ${\mathcal{B}}$, consisting of the triple of objects $Y$, $Z$ and $W$ from the total category ${\mathcal{E}}$ (along with their respective identity morphisms ${id_Y : Y \to Y}$, ${id_Z : Z \to Z}$ and ${id_W : W \to W}$). Once again, these discrete fiber categories can be computed automatically, directly from the \texttt{AbstractFunctor} objects in \textsc{Categorica}, as shown in Figure \ref{fig:Figure26}. Finally, we may consider applying, in the first instance, the additional algebraic equivalence ${\left(F \left( g \right) : F \left( Y \right) \to F \left( Z \right) \right) = \left( id_{F \left( Z \right)} : F \left( Z \right) \to F \left( Z \right) \right)}$ on morphisms in the codomain category, thus yielding the following fibration over a two-object, four-morphism base category:

\begin{equation}
\begin{tikzcd}
& Y \arrow[dr, "g"] \arrow[loop above, "id_Y"] &\\
X \arrow[ur, "f"] \arrow[rr, swap, "g \circ f"] \arrow[loop below, "id_X"] & & Z \arrow[loop below, "id_Z"]
\end{tikzcd} \qquad \mapsto \qquad
\begin{tikzcd}
F \left( X \right) \arrow[loop left, "id_{F \left( X \right)}"] \arrow[rr, bend left, "F \left( f \right)"] \arrow[rr, bend right, swap, "F \left( g \right) \circ F \left( f \right)"] & & F \left( Y \right) = F \left( Z \right) \arrow[loop below, looseness = 15, "F \left( g \right) = id_{F \left( Y \right)} = id_{F \left( Z \right)}"]
\end{tikzcd},
\end{equation}
or alternatively, in the second instance, the additional algebraic equivalences:

\begin{multline}
\left( F \left( h \right) : F \left( Z \right) \to F \left( W \right) \right) = \left( id_{F \left( Z \right)} : F \left( Z \right) \to F \left( Z \right) \right),\\
\text{ and } \qquad \left( F \left( i \right) : F \left( Y \right) \to F \left( W \right) \right) = \left( id_{F \left( Z \right)} : F \left( Z \right) \to F \left( Z \right) \right),
\end{multline}
on morphisms in the codomain category, thus yielding the following fibration over a two-object, six-morphism base category:

\begin{equation}
\begin{tikzcd}
X \arrow[rr, "f"] \arrow[dd, swap, "g"] \arrow[ddrr, bend left = 10, "h \circ g"] \arrow[ddrr, bend right = 10, swap, "i \circ f"] \arrow[loop above, "id_X"] & & Y \arrow[dd, "i"] \arrow[loop above, "id_Y"] \\ \\
Z \arrow[rr, swap, "h"] \arrow[loop below, "id_Z"] & & W \arrow[loop below, "id_W"]
\end{tikzcd} \qquad \mapsto \qquad
\begin{tikzcd}
F \left( X \right) \arrow[loop left, "id_{F \left( X \right)}"] \arrow[rrr, bend left = 30, "F \left( f \right)"] \arrow[rrr, bend left = 10, "F \left( g \right)"] \arrow[rrr, bend right = 10, swap, "F \left( i \right) \circ F \left( f \right)"] \arrow[rrr, bend right = 30, swap, "F \left( h \right) \circ F \left( g \right)"] & & & {\renewcommand{\arraystretch}{0.5}
\begin{tabular}{c}
${F \left( Y \right)}$\\
${ = F \left( W \right)}$\\
${ = F \left( Z \right)}$
\end{tabular}} \arrow[loop below, "\substack{F \left( h \right) = F \left( i \right) \\ = id_{F \left( Y \right)} = id_{F \left( W \right)} = id_{F \left( Z \right)}}"]
\end{tikzcd}.
\end{equation}
The resulting fibrations are no longer discrete, since in the former case the fiber category for object ${F \left( Y \right) = F \left( Z \right)}$ in the base category ${\mathcal{B}}$ now contains a morphism ${g : Y \to Z}$ connecting the pair of objects $Y$ and $Z$ from the total category ${\mathcal{E}}$, and in the latter case the fiber category for object ${F \left( Y \right) = F \left( W \right) = F \left( Z \right)}$ in the base category ${\mathcal{B}}$ now contains a pair of morphisms ${h : Z \to W}$ and ${i : Y \to W}$ connecting the triple of objects $Y$, $Z$ and $W$ from the total category ${\mathcal{E}}$, as seen in the explicit \textsc{Categorica} implementation in Figure \ref{fig:Figure27}. These examples make manifest the precise relationship between category-theoretic fibrations and topological ones: if one considers categories whose objects are points and whose morphisms are paths between those points, then the abstract definition of a Grothendieck fibration reduces concretely to the classical definition a fiber bundle, with the total category ${\mathcal{E}}$ being the total space, the base category ${\mathcal{B}}$ being the base space, and the fiber categories ${F_X}$ being the fibers for each point $X$ in the base.

There exist many other standard features and properties of functors that are currently supported within \textsc{Categorica}, and that can be computed automatically using the \texttt{AbstractFunctor} function, including several that can be effectively assembled as ``composites'' of the various features/properties already discussed above. For instance, a functor yielding an \textit{equivalence} of categories (detectable through the \textit{``EquivalenceFunctorQ''} property) is one that is simultaneously fully faithful \textit{and} essentially surjective; an \textit{embedding} functor (detectable through the \textit{``EmbeddingFunctorQ''} property) is one that is simultaneously faithful \textit{and} injective on objects; and a \textit{full} embedding functor (detectable through the \textit{``FullEmbeddingFunctorQ''} property) is one that is both \textit{fully faithful} and injective on objects. An \textit{embedding} of a category ${\mathcal{C}}$ into a category ${\mathcal{D}}$ via the embedding functor ${F : \mathcal{C} \to \mathcal{D}}$ may be thought of as being a weaker/more lax version of an \textit{inclusion} of a \textit{subcategory} ${\mathcal{C}}$ into a \textit{supercategory} ${\mathcal{D}}$, i.e. a functor ${F : \mathcal{C} \to \mathcal{D}}$ such that the object and morphism sets of the domain category ${\mathcal{C}}$ form (potentially improper) subsets of the respective object and morphism sets of the codomain category ${\mathcal{D}}$:

\begin{equation}
\mathrm{ob} \left( \mathcal{C} \right) \subseteq \mathrm{ob} \left( \mathcal{D} \right), \qquad \text{ and } \qquad \mathrm{hom} \left( \mathcal{C} \right) \subseteq \mathrm{hom} \left( \mathcal{D} \right),
\end{equation}
just as an \textit{equivalence} between categories ${\mathcal{C}}$ and ${\mathcal{D}}$ may be thought of as being a weaker/more lax version of a strict \textit{equality} of categories (wherein the object and morphism sets of categories ${\mathcal{C}}$ and ${\mathcal{D}}$ are required to be strictly identical). Indeed, strict inclusions of subcategories into supercategories are also detectable, by means of the \textit{``InclusionFunctorQ''} property in \texttt{AbstractFunctor} (along with the slightly stronger case in which the inclusion functor is also required to be full, detectable through the \textit{``FullInclusionFunctorQ''} property). Other common types of functors that \textsc{Categorica} contains specialized functionality for handling as part of the \texttt{AbstractFunctor} framework include \textit{endofunctors} (i.e. functors mapping a domain category ${\mathcal{C}}$ to itself) via the \textit{``EndofunctorQ''} property, \textit{identity functors} (i.e. a trivial, much stronger case of an endofunctor in which every object and morphism in the domain category ${\mathcal{C}}$ is mapped to itself) via the \textit{``IdentityFunctorQ''} property, \textit{constant functors} (i.e. functors mapping every object in the domain category ${\mathcal{C}}$ to some fixed object $X$ in the codomain category ${\mathcal{D}}$, and mapping every morphism in ${\mathcal{C}}$ to the identity morphism ${id_X : X \to X}$ on that object) via the \textit{``ConstantFunctorQ''} property, and \textit{conservative functors} (i.e. functors for which a given morphism ${F \left( f \right) : F \left( X \right) \to F \left( Y \right)}$ in the codomain category ${\mathcal{D}}$ being an isomorphism necessarily implies that the corresponding morphism ${f : X \to Y}$ in the domain category ${\mathcal{C}}$ must also have been an isomorphism) via the \textit{``ConservativeFunctorQ''} property, etc.

\clearpage

\section{Natural Transformations and Algorithmic Details}
\label{sec:Section4}

Having considered above the case of homomorphisms/structure-preserving maps between categories, in the form of functors, it is now possible to consider the case of homomorphisms/structure-preserving maps between \textit{functors}, in the form of \textit{natural transformations}. More precisely, if ${F : \mathcal{C} \to \mathcal{D}}$ and ${G : \mathcal{C} \to \mathcal{D}}$ are both functors mapping between the same pair of domain/codomain categories ${\mathcal{C}}$ and ${\mathcal{D}}$, then the map ${\eta : F \Rightarrow G}$ is a natural transformation if and only if it associates every object $X$ in category ${\mathcal{C}}$ to a corresponding morphism ${\eta_X : F \left( X \right) \to G \left( X \right)}$ in category ${\mathcal{D}}$:

\begin{equation}
\forall X \in \mathrm{ob} \left( \mathcal{C} \right), \qquad \left( \eta_X : F \left( X \right) \to G \left( X \right) \right) \in \mathrm{hom} \left( \mathcal{D} \right),
\end{equation}
known as the \textit{component} of the natural transformation ${\eta}$ at $X$, in such a way that, for every morphism ${f : X \to Y}$ in category ${\mathcal{C}}$, one has that the composition ${\eta_Y \circ F \left( f \right)}$ yields the same morphism as the composition ${G \left( f \right) \circ \eta_X}$ in category ${\mathcal{D}}$, i.e:

\begin{equation}
\forall \left( f : X \to Y \right) \in \mathrm{hom} \left( \mathcal{C} \right), \qquad \left( \eta_Y \circ F \left( f \right) : F \left( X \right) \to G \left( Y \right) \right) = \left( G \left( f \right) \circ \eta_X : F \left( X \right) \to G \left( Y \right) \right),
\end{equation}
which can be illustrated simply by means of the following commutative diagram (which therefore commutes in category ${\mathcal{D}}$):

\begin{equation}
\begin{tikzcd}
F \left( X \right) \arrow[rr, "\eta_X"] \arrow[dd, swap, "F \left( f \right)"] \arrow[ddrr, bend left = 10, "G \left( f \right) \circ \eta_X"] \arrow[ddrr, bend right = 10, swap, "\eta_Y \circ F \left( f \right)"] & & G \left( X \right) \arrow[dd, "G \left( f \right)"] \\ \\
F \left( Y \right) \arrow[rr, swap, "\eta_Y"] & & G \left( Y \right)
\end{tikzcd} \qquad \mapsto \qquad
\begin{tikzcd}
F \left( X \right) \arrow[rr, "\eta_X"] \arrow[dd, swap, "F \left( f \right)"] \arrow[ddrr, "\substack{\eta_Y \circ F \left( f \right) \\ = G \left( f \right) \circ \eta_X}"] & & G \left( X \right) \arrow[dd, "G \left( f \right)"] \\ \\
F \left( Y \right) \arrow[rr, swap, "\eta_Y"] & & G \left( Y \right)
\end{tikzcd}.
\end{equation}
Natural transformations abstract and generalize many powerful results and constructions in mathematics, especially in algebraic topology, perhaps most notably the Hurewicz theorem/Hurewicz homomorphism relating homotopy groups and homology groups of spaces: both the construction of homotopy groups and the construction of homology groups over a (pointed) topological space may be formalized as functors from the category ${\mathbf{Top_{*}}}$ of pointed topological spaces and continuous functions between them, to the category ${\mathbf{Grp}}$ of groups and group homomorphisms, with the Hurewicz homomorphism being elegantly represented as a natural transformation between these two functors (and the Hurewicz theorem corresponding to the statement that such a natural transformation necessarily exists). In linear algebra and functional analysis, the relationship between the double-dual ${V^{**}}$ of a vector space $V$ and the original space constitutes another familiar example of a commonplace natural transformation in mathematics; the double-dual operation may be formalized as an endofunctor on the category ${\mathbf{Vect}}$ of vector spaces and linear maps (i.e. a functor from ${\mathbf{Vect}}$ to itself), and there necessarily exists a natural transformation between this double-dual functor and the identity functor on ${\mathbf{Vect}}$. For the case of finite-dimensional vector spaces, this natural transformation becomes a \textit{natural isomorphism} (a stronger case, to be defined momentarily).

\begin{figure}[ht]
\centering
\begin{framed}
\includegraphics[width=0.515\textwidth]{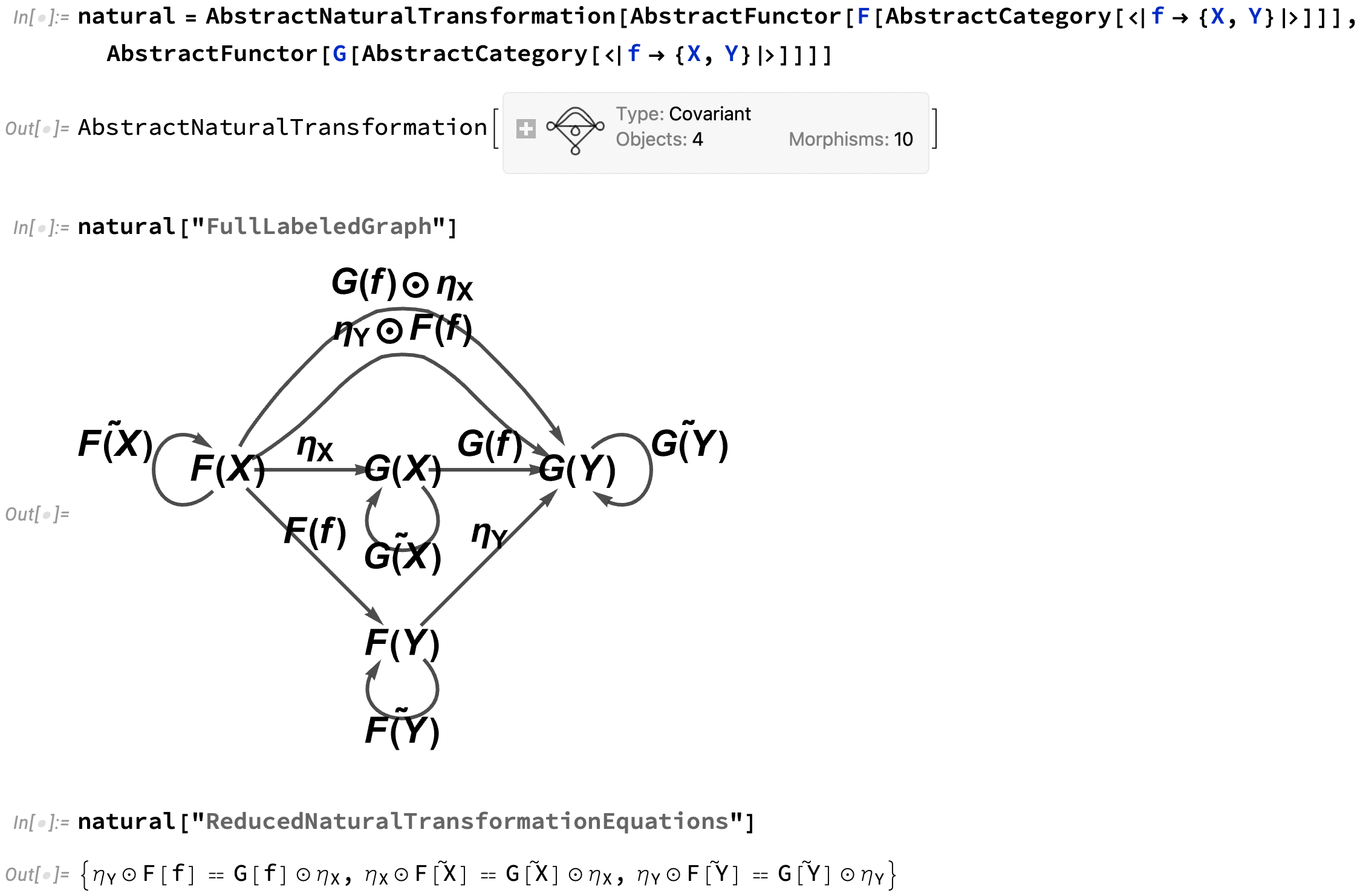}
\vrule
\includegraphics[width=0.475\textwidth]{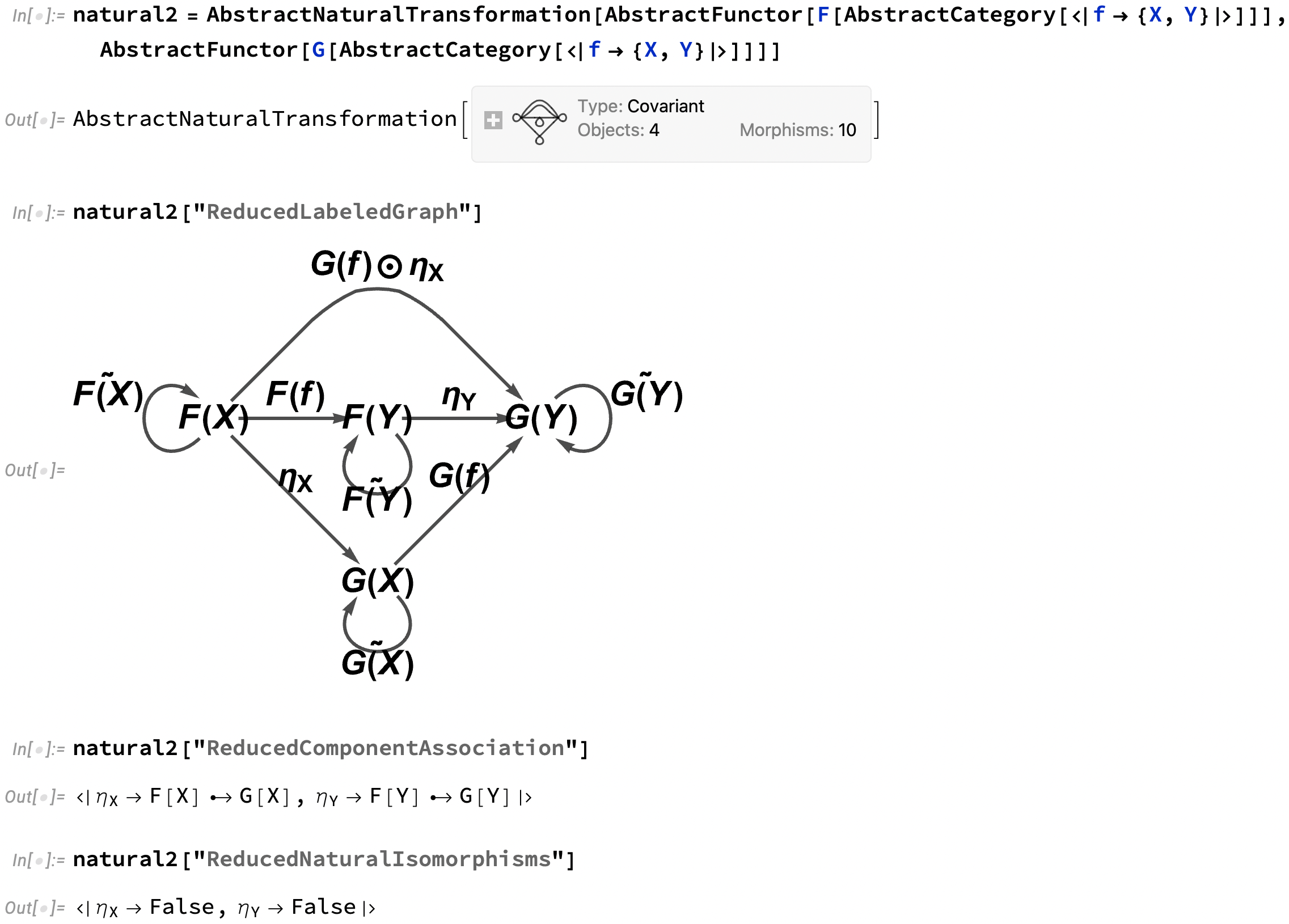}
\end{framed}
\caption{On the left, the \texttt{AbstractNaturalTransformation} object for a pair of functors ${F, G : \mathcal{C} \to \mathcal{D}}$ between two very simple (two-object, three-morphism) categories, showing the minimum set of algebraic equivalences on the component morphisms necessary for the natural transformation to be valid. On the right, the same \texttt{AbstractNaturalTransformation} object for the pair of functors ${F, G : \mathcal{C} \to \mathcal{D}}$ between the same two very simple (two-object, three-morphism) categories, showing the labeled graph representation with the aforementioned algebraic equivalences imposed, as well as the association of component morphisms, showing that neither component is a natural isomorphism.}
\label{fig:Figure28}
\end{figure}

\begin{figure}[ht]
\centering
\begin{framed}
\includegraphics[width=0.495\textwidth]{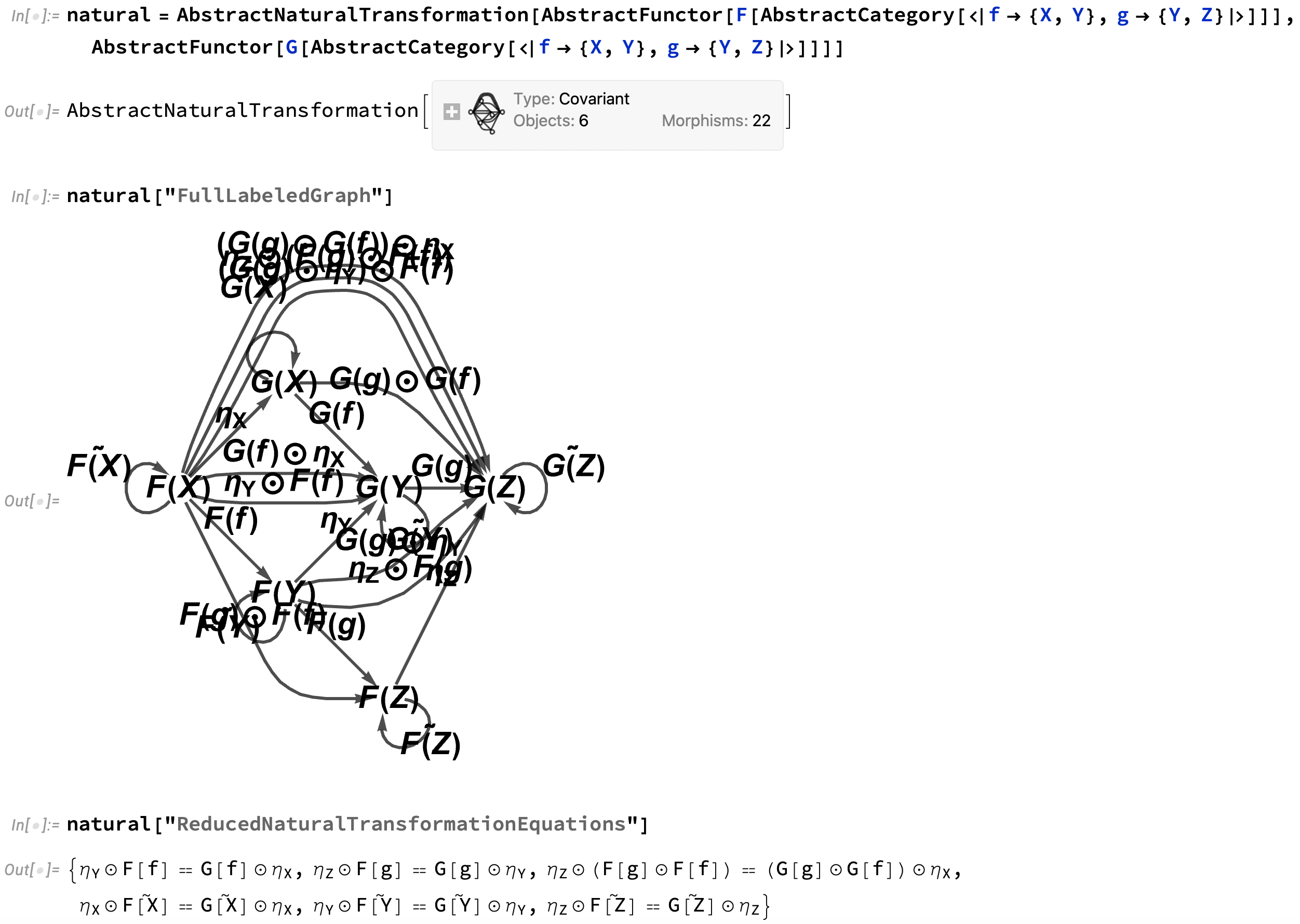}
\vrule
\includegraphics[width=0.495\textwidth]{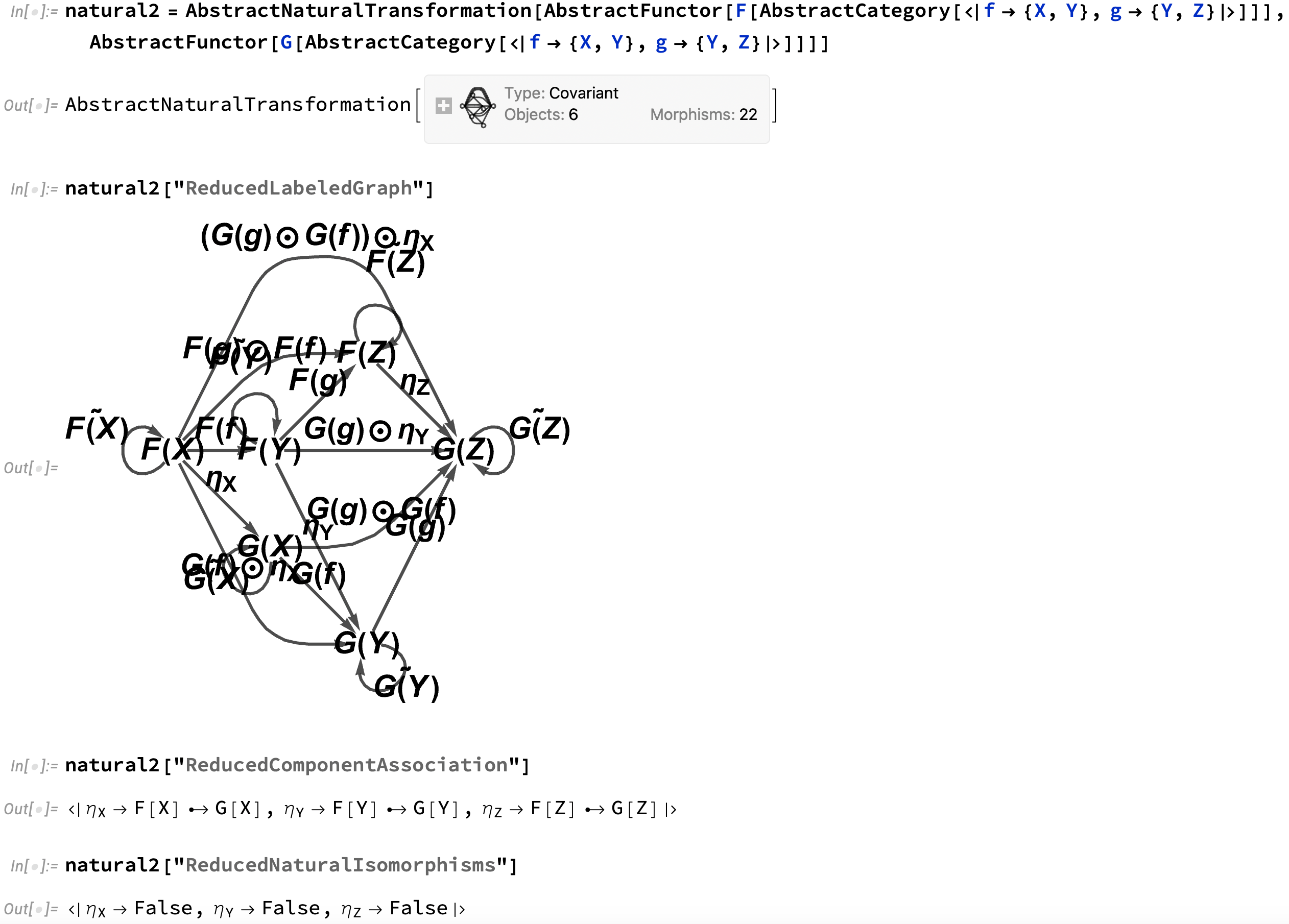}
\end{framed}
\caption{On the left, the \texttt{AbstractNaturalTransformation} object for a pair of functors ${F, G : \mathcal{C} \to \mathcal{D}}$ between two slightly more complex (three-object, six-morphism) categories, showing the minimum set of algebraic equivalences on the component morphisms necessary for the natural transformation to be valid. On the right, the same \texttt{AbstractNaturalTransformation} objects for the pair of functors ${F, G : \mathcal{C} \to \mathcal{D}}$ between the same two slightly more complex (three-object, six-morphism) categories, showing the labeled graph representation with the aforementioned algebraic equivalences imposed, as well as the association of component morphisms, showing that none of the components is a natural isomorphism.}
\label{fig:Figure29}
\end{figure}

\textsc{Categorica} is able to represent arbitrary natural transformations between \texttt{AbstractFunctor} objects using the \texttt{AbstractNaturalTransformation} function. In particular, it is able to compute the minimum set of algebraic equivalences between component morphisms that is necessary in order for a particular transformation between \texttt{AbstractFunctor} objects to be natural, and to represent these equivalences in both equational and (commutative) diagrammatic form, as illustrated in Figure \ref{fig:Figure28} for the minimal case presented above, in which one considers a natural transformation between a pair of functors ${F : \mathcal{C} \to \mathcal{D}}$ and ${G : \mathcal{C} \to \mathcal{D}}$ acting on a common domain category ${\mathcal{C}}$ consisting of a pair of objects $X$ and $Y$ and a single morphism ${f : X \to Y}$, i.e. one has:

\begin{equation}
\begin{tikzcd}
X \arrow[r, "f"] \arrow[loop left, "id_X"] & Y \arrow[loop right, "id_Y"]
\end{tikzcd} \qquad \mapsto \qquad
\begin{tikzcd}
F \left( X \right) \arrow[r, "F \left( f \right)"] \arrow[loop left, "id_{F \left( X \right)}"] & F \left( Y \right) \arrow[loop right, "id_{F \left( Y \right)}"]
\end{tikzcd},
\end{equation}
and:

\begin{equation}
\begin{tikzcd}
X \arrow[r, "f"] \arrow[loop left, "id_X"] & Y \arrow[loop right, "id_Y"]
\end{tikzcd} \qquad \mapsto \qquad
\begin{tikzcd}
G \left( X \right) \arrow[r, "G \left( f \right)"] \arrow[loop left, "id_{G \left( X \right)}"]& G \left( Y \right) \arrow[loop right, "id_{G \left( Y \right)}"]
\end{tikzcd},
\end{equation}
respectively, in which the minimum algebraic condition for naturality is precisely the one given above:

\begin{equation}
\left( \eta_Y \circ F \left( f \right) : F \left( X \right) \to G \left( Y \right) \right) = \left( G \left( f \right) \circ \eta_X : F \left( X \right) \to G \left( Y \right) \right),
\end{equation}
plus a couple of corollary conditions on the identity morphisms. However, although this particular example is especially simple, one does not need to increase the size or complexity of the underlying categories and functors involved very much before the resulting naturality conditions become rather unwieldy to manipulate by hand. For instance, let us consider the next most obvious example, in which the common domain category ${\mathcal{C}}$ of the functors ${F : \mathcal{C} \to \mathcal{D}}$ and ${G : \mathcal{C} \to \mathcal{D}}$ consists instead of three objects $X$, $Y$ and $Z$, along with a pair of composable morphisms ${f : X \to Y}$ and ${g : Y \to Z}$, i.e. one now has:

\begin{equation}
\begin{tikzcd}
& Y \arrow[dr, "g"] \arrow[loop above, "id_Y"] &\\
X \arrow[ur, "f"] \arrow[rr, swap, "g \circ f"] \arrow[loop below, "id_X"] & & Z \arrow[loop below, "id_Z"]
\end{tikzcd} \qquad \mapsto \qquad
\begin{tikzcd}
& F \left( Y \right) \arrow[dr, "F \left( g \right)"] \arrow[loop above, "id_{F \left( Y \right)}"] &\\
F \left( X \right) \arrow[ur, "F \left( f \right)"] \arrow[rr, swap, "F \left( g \right) \circ F \left( f \right)"] \arrow[loop below, "id_{F \left( X \right)}"] & & F \left( Z \right) \arrow[loop below, "id_{F \left( Z \right)}"]
\end{tikzcd},
\end{equation}
and:

\begin{equation}
\begin{tikzcd}
& Y \arrow[dr, "g"] \arrow[loop above, "id_Y"] &\\
X \arrow[ur, "f"] \arrow[rr, swap, "g \circ f"] \arrow[loop below, "id_X"] & & Z \arrow[loop below, "id_Z"]
\end{tikzcd} \qquad \mapsto \qquad
\begin{tikzcd}
& G \left( Y \right) \arrow[dr, "G \left( g \right)"] \arrow[loop above, "id_{G \left( Y \right)}"] &\\
G \left( X \right) \arrow[ur, "G \left( f \right)"] \arrow[rr, swap, "G \left( g \right) \circ G \left( f \right)"] \arrow[loop below, "id_{G \left( X \right)}"] & & G \left( Z \right) \arrow[loop below, "id_{G \left( Z \right)}"]
\end{tikzcd},
\end{equation}
respectively. In this case, the relevant commutative diagram characterizing the natural transformation ${\eta : F \Rightarrow G}$ becomes instead the commutative oblong:

\begin{equation}
\begin{tikzcd}
F \left( X \right) \arrow[rr, "F \left( f \right)"] \arrow[dd, swap, "\eta_X"] & & F \left( Y \right) \arrow[rr, "F \left( g \right)"] \arrow[dd, "\eta_Y"] & & F \left( Z \right) \arrow[dd, "\eta_Z"] \\ \\
G \left( X \right) \arrow[rr, swap, "G \left( f \right)"] & & G \left( Y \right) \arrow[rr, swap, "G \left( g \right)"] & & G \left( Z \right)
\end{tikzcd},
\end{equation}
wherein the two interior squares must both commute:

\begin{equation}
\begin{tikzcd}
F \left( X \right) \arrow[rr, "F \left( f \right)"] \arrow[dd, swap, "\eta_X"] \arrow[ddrr, bend left = 10, "\eta_Y \circ F \left( f \right)"] \arrow[ddrr, bend right = 10, swap, "G \left( f \right) \circ \eta_X"] & & F \left( Y \right) \arrow[rr, "F \left( g \right)"] \arrow[dd, "\eta_Y"] \arrow[ddrr, bend left = 10, "\eta_Z \circ F \left( g \right)"] \arrow[ddrr, bend right = 10, swap, "G \left( g \right) \circ \eta_Y"] & & F \left( Z \right) \arrow[dd, "\eta_Z"] \\ \\
G \left( X \right) \arrow[rr, swap, "G \left( f \right)"] & & G \left( Y \right) \arrow[rr, swap, "G \left( g \right)"] & & G \left( Z \right)
\end{tikzcd} \qquad \mapsto \qquad
\begin{tikzcd}
F \left( X \right) \arrow[rr, "F \left( f \right)"] \arrow[dd, swap, "\eta_X"] \arrow[ddrr, "\substack{\eta_Y \circ F \left( f \right)\\ = G \left( f \right) \circ \eta_X}"] & & F \left( Y \right) \arrow[rr, "F \left( g \right)"] \arrow[dd, "\eta_Y"] \arrow[ddrr, "\substack{\eta_Z \circ F \left( g \right)\\ = G \left( g \right) \circ \eta_Y}"] & & F \left( Z \right) \arrow[dd, "\eta_Z"] \\ \\
G \left( X \right) \arrow[rr, swap, "G \left( f \right)"] & & G \left( Y \right) \arrow[rr, swap, "G \left( g \right)"] & & G \left( Z \right)
\end{tikzcd},
\end{equation}
i.e:

\begin{multline}
\left( \eta_Y \circ F \left( f \right) : F \left( X \right) \to G \left( Y \right) \right) = \left( G \left( f \right) \circ \eta_X : F \left( X \right) \to G \left( Y \right) \right),\\
\text{ and } \qquad \left( \eta_Z \circ F \left( g \right) : F \left( Y \right) \to G \left( Z \right) \right) = \left( G \left( g \right) \circ \eta_Y : F \left( Y \right) \to G \left( Z \right) \right),
\end{multline}
in addition to the outer rectangle:

\begin{equation}
\begin{tikzcd}
F \left( X \right) \arrow[rr, "F \left( f \right)"] \arrow[dd, swap, "\eta_X"] \arrow[ddrrrr, bend left = 10, "\left( \eta_Z \circ F \left( g \right) \right) \circ F \left( f \right)"] \arrow[ddrrrr, bend right = 10, swap, "\left( G \left( g \right) \circ G \left( f \right) \right) \circ \eta_X"] & & F \left( Y \right) \arrow[rr, "F \left( g \right)"] & & F \left( Z \right) \arrow[dd, "\eta_Z"] \\ \\
G \left( X \right) \arrow[rr, swap, "G \left( f \right)"] & & G \left( Y \right) \arrow[rr, swap, "G \left( g \right)"] & & G \left( Z \right)
\end{tikzcd} \qquad \mapsto \qquad
\begin{tikzcd}
F \left( X \right) \arrow[rr, "F \left( f \right)"] \arrow[dd, swap, "\eta_X"] \arrow[ddrrrr, "\substack{\left( \eta_Z \circ F \left( g \right) \right) \circ F \left( f \right)\\ = \left( G \left( g \right) \circ G \left( f \right) \right) \circ \eta_X}"] & & F \left( Y \right) \arrow[rr, "F \left( g \right)"] & & F \left( Z \right) \arrow[dd, "\eta_Z"] \\ \\
G \left( X \right) \arrow[rr, swap, "G \left( f \right)"] & & G \left( Y \right) \arrow[rr, swap, "G \left( g \right)"] & & G \left( Z \right)
\end{tikzcd},
\end{equation}
i.e:

\begin{equation}
\left( \left( \eta_Z \circ F \left( g \right) \right) \circ F \left( f \right) : F \left( X \right) \to G \left( Z \right) \right) = \left( \left( G \left( g \right) \circ G \left( f \right) \right) \circ \eta_X : F \left( X \right) \to G \left( Z \right) \right),
\end{equation}
plus all corresponding conditions on the identity morphisms. As shown in Figure \ref{fig:Figure29}, \textsc{Categorica} is nevertheless easily able to represent, and to perform automatic computations involving, these more complex cases of natural transformations too.

Figures \ref{fig:Figure28} and \ref{fig:Figure29} also demonstrate \textsc{Categorica}'s ability to detect \textit{natural isomorphisms}: components ${\eta_X : F \left( X \right) \to G \left( X \right)}$ of the natural transformation ${\eta : F \Rightarrow G}$ that are also isomorphisms, i.e. where the inverse morphism ${\eta_{X}^{-1} : G \left( X \right) \to F \left( X \right)}$ also exists in the codomain category ${\mathcal{D}}$, such that ${\eta_{X}^{-1} \circ \eta_X : F \left( X \right) \to F \left( X \right)}$ is the identity morphism on ${F \left( X \right)}$ and ${\eta_X \circ \eta_{X}^{-1} : G \left( X \right) \to G \left( X \right)}$ is the identity morphism on ${G \left( X \right)}$:

\begin{equation}
\begin{tikzcd}
F \left( X \right) \arrow[rr, bend left, "\eta_X"] \arrow[loop above, "id_{F \left( X \right)}"] \arrow[loop below, "\eta_{X}^{-1} \circ \eta_X"] & & G \left( X \right) \arrow[ll, bend left, "\eta_{X}^{-1}"] \arrow[loop above, "id_{G \left( Y \right)}"] \arrow[loop below, "\eta_X \circ \eta_{X}^{-1}"]
\end{tikzcd} \qquad \mapsto \qquad
\begin{tikzcd}
F \left( X \right) \arrow[rr, bend left, "\eta_X"] \arrow[loop below, "\eta_{X}^{-1} \circ \eta_X = id_{F \left( X \right)}"] & & G \left( X \right) \arrow[ll, bend left, "\eta_{X}^{-1}"] \arrow[loop below, "\eta_{X} \circ \eta_{X}^{-1} = id_{G \left( X \right)}"]
\end{tikzcd},
\end{equation}
i.e:

\begin{multline}
\forall X \in \mathrm{ob} \left( \mathcal{C} \right), \qquad \text{ component } \eta_X : F \left( X \right) \to G \left( X \right) \qquad \text{ of } \qquad \eta : \left( F : \mathcal{C} \to \mathcal{D} \right) \Rightarrow \left( G : \mathcal{C} \to \mathcal{D} \right),\\
\text{ is a natural isomorphism } \qquad \iff \qquad \exists \left( \eta_{X}^{-1} : G \left( X \right) \to F \left( X \right) \right) \in \mathrm{hom} \left( \mathcal{D} \right),\\
\text{ such that } \qquad \left( \eta_{X}^{-1} \circ \eta_X : F \left( X \right) \to F \left( X \right) \right) = \left( id_{F \left( X \right)} : F \left( X \right) \to F \left( X \right) \right),\\
\text{ and } \qquad \left( \eta_{X} \circ \eta_{X}^{-1} : G \left( X \right) \to G \left( X \right) \right) = \left( id_{G \left( X \right)} : G \left( X \right) \to G \left( X \right) \right).
\end{multline}
In addition to individual objects in the codomain category ${\mathcal{D}}$ being naturally isomorphic, one can also speak of a pair of functors ${F : \mathcal{C} \to \mathcal{D}}$ and ${G : \mathcal{C} \to \mathcal{D}}$ themselves as being naturally isomorphic: if \textit{all} components ${\eta_X : F \left( X \right) \to G \left( X \right)}$ are natural isomorphisms (for every object $X$ in the domain category ${\mathcal{C}}$), then the natural transformation ${\eta : F \Rightarrow G}$ itself is known as a natural isomorphism between functors:

\begin{multline}
\eta : \left( F : \mathcal{C} \to \mathcal{D} \right) \Rightarrow \left( G : \mathcal{C} \to \mathcal{D} \right) \text{ is a natural isomorphism } \qquad \iff \qquad \forall X \in \mathrm{ob} \left( \mathcal{C} \right),\\
\text{ component } \eta_X : F \left( X \right) \to G \left( X \right) \text{ is a natural isomorphism}.
\end{multline}
The intuition underlying the concept of natural isomorphisms between functors can be precisely formalized by considering the \textit{functor category} ${\mathcal{D}^{\mathcal{C}}}$ (assuming that ${\mathcal{C}}$ and ${\mathcal{D}}$ are arbitrary small categories), whose objects are given by the functors ${F : \mathcal{C} \to \mathcal{D}}$ and whose morphisms are given by the natural transformations ${\eta : F \Rightarrow G}$ between such functors, wherein all isomorphisms between objects represent natural isomorphisms between functors. To determine whether a given \texttt{AbstractNaturalTransformation} object in \textsc{Categorica} represents a natural isomorphism between functors, one can use the \textit{``NaturalIsomorphismQ''} property. \textsc{Categorica} is also able to determine whether the codomain categories ${\mathcal{D}}$ between the two functors actually match (via the \textit{``MatchingCodomainsQ''} property), and therefore whether the \texttt{AbstractNaturalTransformation} object does indeed represent a valid natural transformation (via the \textit{``ValidNaturalTransformationQ''} property), meaning \textit{both} that the two codomain categories match \textit{and} that the defining algebraic conditions of a natural transformation hold within those codomain categories.

Finally, it is worth noting that all of \textsc{Categorica}'s core algebraic and diagrammatic reasoning algorithms are represented internally in terms of (hyper)graph rewriting systems over labeled graph representations of quivers and categories, whereby one formalizes (hyper)graph rewriting rules as \textit{spans} of monomorphisms\cite{ehrig}\cite{ehrig2}:

\begin{equation}
\begin{tikzcd}
L & K \arrow[l, swap, "l"] \arrow[r, "r"] & R
\end{tikzcd},
\end{equation}
i.e. pairs of monomorphisms which share a common domain object $K$. Here, the ambient category ${\mathcal{C}}$ in which these monomorphisms exist is the category whose objects are (hyper)graphs and whose morphisms are (hyper)graph inclusion relations. In the above, object $L$ represents the (hyper)graph pattern appearing on the left-hand/input side of the rule, object $R$ represents the (hyper)graph pattern appearing on the right-hand/output side of the rule (i.e. the pattern intended to replace pattern $L$ wherever it appears), and object $K$ represents the ``residual'' sub(hyper)graph pattern shared by both patterns $L$ and $R$ that remains invariant after the left-hand-side $L$ pattern has been ``extracted'', but before the right-hand-side pattern $R$ has been ``injected''. Such a rewriting rule may then be said to \textit{match} a given (hyper)graph object $G$ if there exists a morphism ${m : L \to G}$ in the category ${\mathcal{C}}$. The resulting (hyper)graph $H$ obtained via application of the rewriting rule at this particular match may consequently be computed explicitly by completing the following \textit{double-pushout} diagram\cite{habel}:

\begin{equation}
\begin{tikzcd}
L \arrow[dd, swap, "m"] & & K \arrow[ll, swap, "l"] \arrow[dd, "n"] \arrow[rr, "r"] \arrow[ddll, bend left = 10, "m \circ l"] \arrow[ddll, bend right = 10, swap, "g \circ n"] \arrow[ddrr, bend left = 10, "p \circ r"] \arrow[ddrr, bend right = 10, swap, "h \circ n"] & & R \arrow[dd, "p"]\\ \\
G & & D \arrow[ll, "g"] \arrow[rr, swap, "h"] & & H
\end{tikzcd} \qquad \mapsto \qquad
\begin{tikzcd}
L \arrow[dd, swap, "m"] & & K \arrow[ll, swap, "l"] \arrow[dd, "n"] \arrow[rr, "r"] \arrow[ddll, swap, "\substack{m \circ l\\ = g \circ n}"] \arrow[ddrr, "\substack{p \circ r\\ = h \circ n}"] & & R \arrow[dd, "p"]\\ \\
G & & D \arrow[ll, "g"] \arrow[rr, swap, "h"] & & H
\end{tikzcd},
\end{equation}
i.e. one must proceed to find objects $D$ and $H$, and morphisms $n$, $p$, $g$ and $h$, in the category ${\mathcal{C}}$, such that the diagram above commutes:

\begin{multline}
\exists D, H \in \mathrm{ob} \left( \mathcal{C} \right), \qquad \exists \left( n : L \to D \right), \left( p : R \to H \right), \left( g : D \to G \right), \left( h : D \to H \right) \in \mathrm{hom} \left( \mathcal{C} \right),\\
\text{ such that }  \qquad \left( m \circ l : K \to G \right) = \left( g \circ n : K \to G \right),\\
\text{ and } \qquad \left( p \circ r : K \to H \right) = \left( h \circ n : K \to H \right).
\end{multline}
Furthermore, these objects and morphisms (and therefore, in particular, the object $H$ representing the resulting (hyper)graph) are guaranteed to be unique, up to natural isomorphism, by virtue of the following \textit{universal property} of the double-pushout diagram:

\begin{equation}
\begin{tikzcd}
& & L \arrow[dd, swap, "m"] \arrow[ddddll, bend right = 20, "\forall m^{*}"] \arrow[ddddll, bend right = 40, swap, "u_1 \circ m"] & & K \arrow[ll, swap, "l"] \arrow[dd, "n"] \arrow[rr, "r"] \arrow[ddll, swap, "\substack{m \circ l\\ = g \circ n}"] \arrow[ddrr, "\substack{p \circ r\\ = h \circ n}"] & & R \arrow[dd, "p"] \arrow[ddddrr, bend left = 20, swap, "\forall p^{*}"] \arrow[ddddrr, bend left = 40, "u_2 \circ p"] & &\\ \\
& & G \arrow[ddll, dashed, "\exists! u_1"] & & D \arrow[ll, "g"] \arrow[rr, swap, "h"] \arrow[ddllll, bend left = 20, swap, "\forall g^{*}"] \arrow[ddllll, bend left = 40, "u_1 \circ g"] \arrow[ddrrrr, bend right = 20, "\forall h^{*}"] \arrow[ddrrrr, bend right = 40, swap, "u_2 \circ h"] & & H \arrow[ddrr, dashed, "\exists! u_2"] & &\\ \\
\forall G^{*} & & & & & & & & \forall H^{*}
\end{tikzcd},
\end{equation}
which then collapses down to:

\begin{equation}
\mapsto \qquad \begin{tikzcd}
& & L \arrow[dd, swap, "m"] \arrow[ddddll, bend right, swap, "\forall m^{*} = u_1 \circ m"] & & K \arrow[ll, swap, "l"] \arrow[dd, "n"] \arrow[rr, "r"] \arrow[ddll, swap, "\substack{m \circ l\\ = g \circ n}"] \arrow[ddrr, "\substack{p \circ r\\ = h \circ n}"] & & R \arrow[dd, "p"] \arrow[ddddrr, bend left, "\forall p^{*} = u_2 \circ p"] & &\\ \\
& & G \arrow[ddll, dashed, "\exists! u_1"] & & D \arrow[ll, "g"] \arrow[rr, swap, "h"] \arrow[ddllll, bend left, "\forall g^{*} = u_1 \circ g"] \arrow[ddrrrr, bend right, swap, "\forall h^{*} = u_2 \circ h"] & & H \arrow[ddrr, dashed, "\exists! u_2"] & &\\ \\
\forall G^{*} & & & & & & & & \forall H^{*}
\end{tikzcd},
\end{equation}
i.e. for any other objects ${G^{*}}$ and ${H^{*}}$, and morphisms ${m^{*}}$, ${g^{*}}$, ${p^{*}}$ and ${h^{*}}$, in the category ${\mathcal{C}}$, such that the resulting squares commute:

\begin{multline}
\forall G^{*}, H^{*} \in \mathrm{ob} \left( \mathcal{C} \right), \qquad \forall \left( m^{*} : L \to G^{*} \right), \left( g^{*} : D \to G^{*} \right), \left( p^{*} : R \to H^{*} \right), \left( h^{*} : D \to H^{*} \right) \in \mathrm{hom} \left( \mathcal{C} \right)\\
\text{ such that } \qquad \left( m^{*} \circ l : K \to G^{*} \right) = \left( g^{*} \circ n : K \to G^{*} \right),\\
\text{ and } \qquad \left( p^{*} \circ r : G \to H^{*} \right) = \left( h^{*} \circ n : K \to H^{*} \right),
\end{multline}
there necessarily exists a pair of unique morphisms ${u_1 : G \to G^{*}}$ and ${u_2 : H \to H^{*}}$ such that the diagram above commutes, i.e:

\begin{multline}
\exists! \left( u_1 : G \to G^{*} \right), \qquad \text{ such that } \qquad \left( m^{*} : L \to G^{*} \right) = \left( u_1 \circ m : L \to G^{*} \right),\\
\text{ and } \qquad \left( g^{*} : D \to G^{*} \right) = \left( u_1 \circ g : D \to G^{*} \right),
\end{multline}
and:

\begin{multline}
\exists! \left( u_2 : H \to H^{*} \right), \qquad \text{ such that } \qquad \left( p^{*} : R \to H^{*} \right) = \left( u_2 \circ p : R \to H^{*} \right),\\
\text{ and } \qquad \left( h^{*} : D \to H^{*} \right) = \left( u_2 \circ h : D \to H^{*} \right).
\end{multline}
\textsc{Categorica} is able to exploit the formalism of double-pushout (hyper)graph rewriting in order to convert automatically between potentially slow, tedious algebraic computations and fast, explicit diagrammatic ones. To this end, the core of \textsc{Categorica}'s hypergraph rewriting system has been built on top of the \textsc{Gravitas} computational general relativity framework\cite{gorard10}\cite{gorard11} (which has previously been applied to the study of black hole physics\cite{gorard12}\cite{gorard13}\cite{gorard14} and algebraic quantum field theory\cite{gorard15}), which employs a bespoke suite of highly optimized, and massively parallelized, algorithms for rewriting large hypergraphs subject to generic algebraic constraints. We note also that such rewriting systems, especially in the \textit{multiway}/non-deterministic case, admit an elegant mathematical description in terms of higher categories and (higher) homotopy types\cite{arsiwalla}\cite{arsiwalla2}.

\clearpage

\section{Concluding Remarks}
\label{sec:Section5}

In this article, we have intended to provide a reasonably comprehensive overview of several of the core data structures of the \textsc{Categorica} framework, especially the foundational \texttt{AbstractQuiver}, \texttt{AbstractCategory}, \texttt{AbstractFunctor} and \texttt{AbstractNaturalTransformation} constructions, as well as a number of its core design principles; in particular, it is hoped that this overview will lay the groundwork for the forthcoming second article in this series, which will showcase and outline \textsc{Categorica}'s more advanced functionality for handling algebraic computations and automated theorem-proving for abstract universal properties (e.g. products, coproducts, pullbacks, pushouts, limits and colimits), support for monoidal categories and computations involving string diagrams, preliminary support for elementary toposes, and initial support for higher categories. There exist many planned directions for future research and development involving the \textsc{Categorica} framework, including the extension of the \texttt{AbstractFunctor} function to support the detection and handling of Grothendieck fibrations in complete generality, the inclusion of adjoint functor relationships between \texttt{AbstractFunctor} objects (and therefore support for weak equivalences between \texttt{AbstractCategory} objects), and the incorporation of full support for both strict and weak 2-categories (i.e. bicategories) by extending the current support for single-object bicategories via the \texttt{AbstractStrictMonoidalCategory} functionality\cite{maclane2}. Similarly, all instances of \texttt{AbstractCategory} objects in the \textsc{Categorica} framework are currently assumed by default to be enriched over \textbf{Set} (since the object sets and morphism sets are assumed to be pure Wolfram Language list structures); it would therefore be instructive to extend the present functionality to accommodate \texttt{AbstractCategory} objects enriched directly over arbitrary \texttt{AbstractStrictMonoidalCategory} objects\cite{kelly}\cite{kelly2}. Indeed, although there exists much preliminary infrastructural support for general monoidal categories, there are also many important cases and constructions (including braided monoidal categories\cite{joyal}\cite{chari}, monoidal functors\cite{aguiar}, etc.) which are used throughout the field of applied category, and which are not currently supported, although appropriate functionality for many of these is presently under active development.

\section*{Acknowledgments}
The author would like to thank participants in the GReTA - Graph Transformation Theory and Applications - seminar series (especially Nicolas Behr and Aleks Kissinger), and members of the Topos Institute (especially Owen Lynch and David Spivak) for various discussions and comments on aspects of the design and implementation of the \textsc{Categorica} framework. The author would also like to thank Mohamed Barakat, Manojna Namuduri and Stephen Wolfram for their encouragement and design input at several stages of the development process.

\end{document}